# Confidence Intervals for Causal Effects with Invalid Instruments using Two-Stage Hard Thresholding with Voting

Zijian Guo

*Rutgers University, Piscataway, USA*

Hyunseung Kang

*University of Wisconsin-Madison, Madison, USA*

T. Tony Cai

*University of Pennsylvania, Philadelphia, USA*

Dylan S. Small

*University of Pennsylvania, Philadelphia, USA*

**Summary**.  A major challenge in instrumental variables (IV) analysis is to find instruments that are valid, or have no direct effect on the outcome and are ignorable. Typically one is unsure whether all of the putative IVs are in fact valid.  We propose a general inference procedure in the presence of invalid IVs, called Two-Stage Hard Thresholding (TSHT) with voting.  TSHT uses two hard thresholding steps to select strong instruments and generate candidate sets of valid IVs. Voting takes the candidate sets and uses majority and plurality rules to determine the true set of valid IVs. In low dimensions, if the sufficient and necessary identification condition under invalid instruments is met, which is more general than the so-called 50% rule or the majority rule, our proposal (i) correctly selects valid IVs, (ii) consistently estimates the causal effect, (iii) produces valid confidence intervals for the causal effect, and (iv) has oracle-optimal width.  In high dimensions, we establish nearly identical results without oracle-optimality. In simulations, our proposal outperforms traditional and recent methods in the invalid IV literature. We also apply our method to re-analyze the causal effect of education on earnings.

## 1. Introduction

### 1.1. *Motivation: Invalid Instruments*

Instrumental variables (IV) analysis is a popular method to deduce causal effects in the presence of unmeasured confounding.  Informally, an IV analysis requires instruments that (A1) are associated with the exposure, (A2) have no direct pathway to the outcome, and (A3) are not related to unmeasured variables that affect the exposure and the outcome; see Section 2.1 for details.  A major challenge in IV analysis is to find valid instruments, i.e. instruments that satisfy (A2) and (A3).

In practice, potential candidate instruments become more plausible as valid instruments after controlling for covariates [Hernán and Robins, 2006, Baiocchi et al., 2014]. For example, a long-standing interest in economics is studying the causal effect of education on earnings and often, IV analysis is used to deduce the effect [Angrist and Krueger, 1991, Card, 1993, 1999]. A popular instrument in this analysis is a person's proximity to a college when growing up [Card, 1999, 1993]. However, proximity to a college may be related to a person's socioeconomic status, characteristics of a person's high school and other covariates that may affect a person's earnings. Thus, these covariates are controlled for in order for proximity to college to be a valid instrument. With the growing trend toward collecting large data sets with many variables, this approach of creating plausibly valid instruments by conditioning on potentially many covariates has increasing promise [Hernán and Robins, 2006, Swanson and Hernán, 2013, Baiocchi et al., 2014, Varian, 2014, Imbens, 2014].

Yet, even after controlling for covariates, some IVs may still turn out to be invalid and subsequent analysis assuming that all the IVs are valid after conditioning can be misleading [Murray, 2006]. For example, suppose for studying the causal effect of education on earnings, we used proximity as an IV and we control for high school test scores, high school size, individual's genetic makeup, and parents' education and socioeconomic background. But, if living close to college had other benefits beyond getting more education, say by being exposed to many programs available to high school students for job preparation and employers who come to the area to discuss employment opportunities for college students, then the IV, proximity to college, can directly affect individual's earning potential and violate (A2) [Card, 1999]. This problem is also prevalent in other applications of instrumental variables, most notably in Mendelian Randomization [Davey Smith and Ebrahim, 2003, 2004] where the instruments are genetic in nature and some instruments are likely to be invalid due to having pleiotropic effects [Lawlor et al., 2008, Burgess et al., 2015].

This paper tackles the problem of constructing confidence intervals for causal effects when invalid instruments may be present. We consider two major cases. The first case is where the number of covariates is small and fixed relative to the sample size; this setting is typical in Mendelian Randomization studies and many traditional applied settings. The second case is where the number of covariates and/or instruments is growing and may exceed the sample size, which is becoming more prevalent with modern large data sets.

## 1.2. Prior Work

In non-IV settings with high dimensional covariates as controls, Zhang and Zhang [2014], Javanmard and Montanari [2014], van de Geer et al. [2014], Belloni et al. [2014] and Cai and Guo [2017] provide confidence intervals for a treatment effect. In IV settings with high dimensional covariates (or IVs), Gautier and Tsybakov [2011], Belloni et al. [2012], Fan and Liao [2014] and Chernozhukov et al. [2015] provide confidence intervals for a treatment effect, under the assumption that all the IVs are valid after controlling for said covariates. In invalid IV settings, Kolesár et al. [2015] and Bowden et al. [2015] provide inferential methods for treatment effects. How-ever, the method requires that the effects of the instruments on the treatment are orthogonal to their direct effects on the outcome, a stringent assumption. Bowden et al. [2016], Burgess et al. [2016], Kang et al. [2016b] and Windmeijer et al. [2016] also work on the invalid IV setting, but without making the stringent orthogonality assumption.

Unfortunately, all of these methods (i) work only in the low dimensional setting and (ii) rely on the sufficient condition in Han [2008] and Kang et al. [2016b] known as the "50% rule" or the "majority rule" where majority of the instruments must be valid to establish consistency or inferential guarantees (see Section 2.2 for details); to the best our knowledge, no method in this literature establishes consistency or inferential guarantees under the necessary and sufficient condition for identification established in Theorem 1 of Kang et al. [2016b], including the setting where the majority rule is violated. For example, while Kang et al. [2016b] establish sufficient and necessary condition for identification under invalid instruments, its proposed estimator for the treatment effect, sisVIVE, relies on more stringent conditions for estimation consistency and is only demonstrated to work in low dimensional settings. Also, sisVIVE has no theoretical guarantees on inference, i.e. a $1 - \alpha$ confidence interval derived from sisVIVE may not cover the true treatment effect at least $1 - \alpha$ of the time. Similarly, the median estimator of Bowden et al. [2016] and Burgess et al. [2016] is only shown to be a consistent estimator of the treatment effect and there are no theoretical guarantees on inference. Both papers also focus on the setting where the IVs are completely uncorrelated/orthogonal to each other, although modifications are possible, and the data is derived from two independent samples. Finally, the median method may not be optimal in the sense that the confidence interval derived from it may not have the same length as the oracle counterpart confidence interval, say the two-stage least squares (TSLS, see Section 5.1 for details), that knows which instruments are invalid, a priori. Work by Windmeijer et al. [2016], which uses the adaptive Lasso [Zou, 2006], does provide theoretical guarantees and can handle correlated IVs without any modification. But, Windmeijer et al. [2016]'s method relies heavily on the 50% rule to establish these theoretical guarantees and the method fails when this rule does not hold. Finally, work by Kang et al. [2016a] does not rely on the 50% rule to obtain inferential quantities, but it is conservative and only works in low dimensional settings. See Section 3.6 for additional discussion of these methods.

### 1.3. Our Contributions

Our work makes three major contributions in inferring treatment effects in the presence of possibly invalid instruments. First, we propose a novel two-stage hard thresholding (TSHT) with voting that works both in low and high dimensional settings. Second, in the low dimensional setting, our method is the first method to be complete; our method only relies on the sufficient and necessary condition for identification under invalid instruments to (i) correctly select valid IVs, (ii) consistently estimate the causal effect, (iii) produce confidence intervals with the desired level of coverage, and (iv) achieve oracle optimality in the sense that it performs as well asymptotically as the oracle procedure that knows which instruments are valid. In particular, our method can guarantee these properties even when more than 50% of the instruments are invalid. Third, in the high dimensional setting, our method achieves the same selection, estimation, and inferential guarantees without the oracle-optimality.

As the name suggests, the two key components of our method are TSHT and voting. Specifically, TSHT sequentially uses hard thresholding to estimate a set of non-redundant, individually strong instruments (see Definition 2 and Section 3.2 for details) and uses this set to estimate candidate sets of valid IVs (see Definition 1 and Section 3.3 for details). Voting takes the candidate sets of valid IVs from TSHT and estimates a single set of valid IVs by selecting IVs that satisfy the identification conditions for invalid instruments, the plurality rule and the majority rule (see Sections 2.2 and 3.4 for details).

The outline of the paper is as follows. After describing the model setup in Section 2, we describe our procedure TSHT with voting in Section 3 and provide theoretical justification for it in Section 4. In Section 5, we investigate the performance of our procedure in a simulation study and compare with existing methods, in particular the median method with bootstrapping by Bowden et al. [2016] and Burgess et al. [2016] and the adaptive Lasso method by Windmeijer et al. [2016]. We find that our method and Windmeijer et al. [2016] are comparable when the 50% rule holds, although the median estimator suffers from coverage and optimality issues. However, when the 50% rule fails, our method dominates all these methods. In Section 6, we present an empirical study where we revisit the question of the causal effect of years of schooling on income using data from the Wisconsin Longitudinal Study. We provide conclusions and discussions in Section 7.

## 2. Model

To define causal effects and instruments, the potential outcome approach [Neyman, 1923, Rubin, 1974] laid out in Holland [1988] is used. For each individual $i \in \{1, \ldots, n\}$, let $Y_i^{(d, \mathbf{z})} \in \mathbb{R}$ be the potential outcome if the individual were to have exposure/treatment $d \in \mathbb{R}$ and instruments $\mathbf{z} \in \mathbb{R}^{p_z}$. Let $D_i^{(\mathbf{z})} \in \mathbb{R}$ be the potential exposure if the individual had instruments $\mathbf{z} \in \mathbb{R}^{p_z}$. For each individual, only one possible realization of $Y_i^{(d, \mathbf{z})}$ and $D_i^{(\mathbf{z})}$ is observed, denoted as $Y_i$ and $D_i$, respectively, based on his/her observed instrument values $\mathbf{Z}_{i\cdot} \in \mathbb{R}^{p_z}$ and exposure value $D_i$. We also denote baseline covariates for each individual $i$ as $\mathbf{X}_{i\cdot} \in \mathbb{R}^{p_x}$. In total, $n$ sets of outcome, exposure, instruments and baseline covariates, denoted as $(Y_i, D_i, \mathbf{Z}_{i\cdot}, \mathbf{X}_{i\cdot})$, are observed in an i.i.d. fashion.

Let $\mathbf{Y} = (Y_1, \ldots, Y_n)$ to be an $n$-dimensional vector of observed outcomes, $\mathbf{D} = (D_1, \ldots, D_n)$ to be an $n$-dimensional vector of observed exposures, $\mathbf{Z}$ to be a $n$ by $p_z$ matrix of instruments where row $i$ consists of $\mathbf{Z}_{i\cdot}$, and $\mathbf{X}$ to be an $n$ by $p_x$ matrix of covariates where row $i$ consists of $\mathbf{X}_{i\cdot}$. Let $\mathbf{W}$ to be an $n$ by $p = p_z + p_x$ matrix where $\mathbf{W}$ is the result of concatenating the matrices $\mathbf{Z}$ and $\mathbf{X}$ and let $\mathbf{\Sigma}^* = \mathbf{E}\left(\mathbf{W}_{i\cdot} \mathbf{W}_{i\cdot}^\mathsf{T}\right)$ be the positive definite covariance matrix of the instrument-covariate matrix.

For any vector $\mathbf{v} \in \mathbb{R}^p$, let $v_j$ denote the $j$th element of $\mathbf{v}$. Let $\|\mathbf{v}\|_1$, $\|\mathbf{v}\|_2$, and $\|\mathbf{v}\|_\infty$ denote the usual $1, 2$ and $\infty$-norms, respectively. Let $\mathrm{supp}(\mathbf{v}) \subseteq \{1, \dots, p\}$ denote the support of the vector $\mathbf{v}$, $\mathrm{supp}(\mathbf{v}) = \{j : v_j \neq 0\}$ and $\|\mathbf{v}\|_0$ denote the size of the support of $\mathbf{v}$, or equivalently, the number of non-zero elements in $\mathbf{v}$. For a set $J$, let $|J|$ denote its cardinality and $J^C$ denote its complement. For an $n$ by $p$ matrix $\mathbf{M} \in \mathbb{R}^{n \times p}$, we denote the $(i, j)$ element of matrix $\mathbf{M}$ as $M_{ij}$, the $i$th row as $\mathbf{M}_{i.}$, and the $j$th column as $\mathbf{M}_{.j}$. For any sets $A \subseteq \{1, \dots, n\}$ and $B \subseteq \{1, \dots, p\}$, let $\mathbf{M}_{A,B}$ denote the submatrix formed by the rows specified by the set $A$ and the columns specified by the set $B$. Let $\mathbf{M}^\mathsf{T}$ be the transpose of $\mathbf{M}$ and $\|\mathbf{M}\|_\infty$ represent the element-wise matrix sup-norm of matrix $\mathbf{M}$. For a symmetric matrix $\mathbf{M}$, let $\lambda_{\min}(\mathbf{M})$ and $\lambda_{\max}(\mathbf{M})$ denote the smallest and largest eigenvalues of $\mathbf{M}$, respectively.

For a sequence of random variables $X_n$, let $X_n \xrightarrow{p} X$ and $X_n \xrightarrow{d} X$ denote that $X_n$ converges to $X$ in probability and in distribution, respectively. For any two sequences $a_n$ and $b_n$, let $a_n \gg b_n$ denote that $\limsup_{n \to \infty} b_n / a_n = 0$; similarly, let $a_n \ll b_n$ denote $b_n \gg a_n$.

### 2.1. Models and Instrumental Variables Assumptions

We consider the Additive LInear, Constant Effects (ALICE) model of Holland [1988] and extend it to allow for possibly invalid instruments as in Small [2007] and Kang et al. [2016b]. For two possible values of the exposure $d', d$ and instruments $\mathbf{z}', \mathbf{z}$, we assume the following potential outcome model

$$Y_i^{(d', \mathbf{z}')} - Y_i^{(d, \mathbf{z})} = (\mathbf{z}' - \mathbf{z})^\mathsf{T} \boldsymbol{\kappa}^* + (d' - d)\beta^*, \quad \mathbf{E}(Y_i^{(0, \mathbf{0})} \mid \mathbf{Z}_{i.}, \mathbf{X}_{i.}) = \mathbf{Z}_{i.}^\mathsf{T} \boldsymbol{\eta}^* + \mathbf{X}_{i.}^\mathsf{T} \boldsymbol{\phi}^* \quad (1)$$

where $\boldsymbol{\kappa}^*, \beta^*, \boldsymbol{\eta}^*$, and $\boldsymbol{\phi}^*$ are unknown parameters. The parameter $\beta^*$ represents the causal parameter of interest, the causal effect (divided by $d' - d$) of changing the exposure from $d'$ to $d$ on the outcome. The parameter $\boldsymbol{\phi}^*$ represents the impact of covariates on the baseline potential outcome $Y_i^{(0, \mathbf{0})}$. The parameter $\boldsymbol{\kappa}^*$ represents the violation of the no direct effect assumption between the instrument and the outcome. The parameter $\boldsymbol{\eta}^*$ represents the presence of unmeasured confounding between the instrument and the outcome.

The model parameters $\boldsymbol{\kappa}^*$ and $\boldsymbol{\eta}^*$ parametrizes the assumptions (A2) and (A3) in Section 1.1 as follows. The exclusion restriction (A2), which is typically stated [Angrist et al., 1996] as $Y_i^{(d, z)} = Y_i^{(d, z')}$ for all $z, z' \in \mathbb{R}$ implies $\boldsymbol{\kappa}^* = 0$. Also, the assumption of no unmeasured confounding of the IV-outcome relationship (A3), which is typically stated [Angrist et al., 1996] as $Y_i^{(d, z)}$ and $D_i^{(z)}$ are independent of $Z_i$ for all $d, z \in \mathbb{R}$ implies $\boldsymbol{\eta}^* = 0$; we note that Angrist et al. [1996] considered the instrument to have a non-zero average causal effect on the exposure, hence the potential outcome notation for the exposure $D_i^{(z)}$. Our model parameters $\boldsymbol{\kappa}^*$ and $\boldsymbol{\eta}^*$ encode a particular case of the definitions of exclusion restriction (A2) and no unmeasured confounding instrument (A3) in Angrist et al. [1996] where we assume an additive, linear, and constant treatment effect $\beta^*$; see Holland [1988] and Appendix, Section 1.4 of Hernán and Robins [2006] for additional discussions about different formalizations of the IV assumptions (A2) and (A3).

Let $\boldsymbol{\pi}^* = \boldsymbol{\kappa}^* + \boldsymbol{\eta}^*$, $e_i = Y_i^{(0,0)} - \mathbf{E}(Y_i^{(0,0)} \mid \mathbf{Z}_{i.}, \mathbf{X}_{i.})$, and $\mathrm{Var}(e_i|\mathbf{Z}_{i.}, \mathbf{X}_{i.}) = \sigma^2$. When we combine equation (1) along with the definition of $e_i$, we have the following observed data model, which is also known as the under-identified single-equation linear model in econometrics (page 83 of Wooldridge [2010])

$$Y_i = \mathbf{Z}_{i.}^{\mathsf{T}}\boldsymbol{\pi}^* + D_i\beta^* + \mathbf{X}_{i.}^{\mathsf{T}}\boldsymbol{\phi}^* + e_i, \quad \mathbf{E}(e_i \mid \mathbf{Z}_{i.}, \mathbf{X}_{i.}) = 0 \qquad (2)$$

The observed model is not a usual regression model because $D_i$ might be correlated with $e_i$. In particular, the parameter $\beta^*$ measures the causal effect of changing $D$ on $Y$ rather than an association. Also, the parameter $\boldsymbol{\pi}^*$ in model (2) combines two assumptions, the exclusion restriction (A2) parametrized by $\boldsymbol{\kappa}^*$, and the no unmeasured confounding parametrized by $\boldsymbol{\eta}^*$. If both assumptions are satisfied, $\boldsymbol{\kappa}^* = \boldsymbol{\eta}^* = 0$ so that $\boldsymbol{\pi}^* = 0$ and the instruments are said to be valid [Murray, 2006]. Hence, $\boldsymbol{\pi}^*$ can be used to define valid IVs and we formalize the definition of valid IVs as follows.

DEFINITION 1. *Suppose we have $p_z$ candidate instruments along with the models (1)–(2). We say that instrument $j = 1, \ldots, p_z$ satisfies both (A2) and (A3), or is valid, if $\pi_j^* = 0$ and we denote $\mathcal{P}^*$ to be the set of valid instruments.*

We now discuss our definition of valid instruments in relation with those in the literature. When there is one instrument (i.e. $p_z = 1$), Definition 1 of a valid instrument matches the formal definition of a valid instrument in Holland [1988] and Angrist et al. [1996]. Specifically, as discussed before, by Angrist et al. [1996], the exclusion restriction and no unmeasured confounding imply $\phi^* = \psi^* = 0$ and consequently, $\pi^* = 0$, which is the definition of a valid IV in Definition 1. By using $\pi^*$ to define a valid instrument, we combine the two assumptions, (A2) and (A3); this is common in econometrics [Holland, 1988, Imbens and Angrist, 1994, Angrist et al., 1996, Wooldridge, 2010] where the combination of both assumptions is referred to as instrument exogeneity. When there are several instruments (i.e. $p_z > 1$), our models and definition of valid IVs can be viewed as a generalization of Holland [1988] and Angrist et al. [1996]. It is important to note that in this generalization, Definition 1 defines the validity of an instrument $j$ in the context of the set of instruments $\{1, \ldots, p_z\}$ being considered. Specifically, an instrument $j$ could be valid in the context of the set $\{1, \ldots, p_z\}$ (i.e. $\pi_j^* = 0$), but invalid if considered alone because $\mathbf{Z}_{.j}$ may be associated with or causally affect another IV $\mathbf{Z}_{.j'}$, $j \neq j'$ where $\pi_{j'}^* \neq 0$. Finally, in Mendelian Randomization literature, the parameter $\boldsymbol{\pi}^*$ is used to define a valid instrument in exactly the same way as our definition of a valid instrument [Bowden et al., 2016, Burgess et al., 2016].

In addition to the model for the outcome, we assume a linear association/observational model between the exposure $D_i$, the instruments $\mathbf{Z}_{i.}$, and the covariates $\mathbf{X}_{i.}$.

$$D_i = \mathbf{Z}_{i.}^{\mathsf{T}}\boldsymbol{\gamma}^* + \mathbf{X}_{i.}^{\mathsf{T}}\boldsymbol{\psi}^* + \epsilon_{i2}, \quad \mathbf{E}(\epsilon_{i2}|\mathbf{Z}_{i.}, \mathbf{X}_{i.}) = 0 \qquad (3)$$

Each element $\gamma_j^*$, $j = 1, \ldots, L$ is the partial correlation between the $j$th instrument and $D_i$. The parameter $\boldsymbol{\psi}^*$ represents the association between the covariates and $D_i$. Also, unlike the models (1)-(2), we do not need a causal model between $D_i$, $\mathbf{Z}_{i\cdot}$, and $\mathbf{X}_{i\cdot}$; this is because the constant effect assumption we make in model (1) eliminates the need to assume a causal instrument (see Angrist et al. [1996] for details).

Based on model (3), we formally define assumption (A1), the instruments' relevance to the exposure; this is sometimes referred to as existence of non-redundant instruments in econometrics [Cheng and Liao, 2015].

DEFINITION 2. *Suppose we have $p_z$ candidate instruments along with the model (3). We say that instrument $j = 1, \ldots, p_z$ satisfies (A1), or is a non-redundant IV, if $\gamma_j^* \neq 0$ and we denote $\mathcal{S}^*$ to be the set of these instruments.*

We now discuss our definition of instruments' relevance with those in the literature. In econometrics, satisfying (A1) has been defined in a global sense where (A1) is satisfied if $\boldsymbol{\gamma}^* \neq 0$ (see Section 5.2.1 of Wooldridge [2010]). However, this global definition can be misleading in the presence of multiple candidate instruments. For example, it is possible that $\gamma_1^* \neq 0$ while $\gamma_j^* = 0$ for all $j \neq 1$ so that only the first instrument has an effect on the exposure while the rest do not. Using the global definition would imply that all $p_z$ instruments satisfy (A1) while Definition 2 makes it explicit and, perhaps less ambiguous, that it is only the first instrument $j = 1$ that satisfies (A1). Both the global definition and Definition 2 are the same if $\gamma_j^* \neq 0$ for all $j$, that is where we only include relevant instruments; this is typically the case in practice, especially in Mendelian Randomization. Finally, like Definition 1, when there is only one candidate instrument (i.e. $p_z = 1$), Definition 2 is a special case of the more general definition of (A1) in Angrist et al. [1996] where we assume a model.

Combining Definitions 1 and 2, we can formally define the usual three core conditions for instruments, i.e. (A1)-(A3).

DEFINITION 3. *Suppose we have $p_z$ candidate instruments along with the models (1)–(3). We say that $Z_{ij}$, $j = 1, \ldots, p_z$, is an instrument if (A1) − (A3) are satisfied, i.e. if $\pi_j^* = 0$ and $\gamma_j^* \neq 0$. Let $\mathcal{V}^* = \mathcal{S}^* \cap \mathcal{P}^*$ be the set of instruments.*

We discuss our three core conditions that define an instrument in Definition 3 with those in the literature. Again, when there is only one instrument (i.e. $p_z = 1$), Definition 3 is identical to the three core conditions that define an instrument in Holland [1988]. Definition 3 is also a special case of the three core conditions in Angrist et al. [1996] where, as mentioned before, we assume an additive, linear, constant effects model. When there are multiple instruments (i.e. $p_z > 1$), our models (1)–(3) and Definition 3 can be viewed as a generalization of the definition of instruments in Holland [1988].

Finally, the models presented above are ubiquitous IV models in econometrics [Wooldridge, 2010], genetic epidemiology and Mendelian Randomization [Didelez and Sheehan, 2007]. However, we generalize these models in two important ways: (i) the model in equation (2) allows for possibly invalid instruments and (ii) the number of covariates $p_x$ (and even the number of instruments $p_z$) may be larger than the sample size $n$.

For the rest of the paper, we define the sparsity level of $\boldsymbol{\pi}^*, \boldsymbol{\phi}^*, \boldsymbol{\gamma}^*$ and $\boldsymbol{\psi}^*$ as, $s_{z2} = \|\boldsymbol{\pi}^*\|_0, s_{x2} = \|\boldsymbol{\phi}^*\|_0, s_{z1} = \|\boldsymbol{\gamma}^*\|_0$ and $s_{x1} = \|\boldsymbol{\psi}^*\|_0$. Let $s = \max\{s_{z2}, s_{x2}, s_{z1}, s_{x1}\}$.

### 2.2. Identification of Model Parameters

Identification of the model parameters with invalid instruments has been discussed in prior works [Han, 2008, Bowden et al., 2015, Kolesár et al., 2015, Kang et al., 2016b]. This section reviews and extends these works to guide the discussion of our inferential method for the treatment effect $\beta^*$.

A popular identifying assumption is that majority of the instruments are valid

(IN1-M)  **Majority rule/50% rule**: $|\mathcal{V}^*| > \frac{1}{2}|\mathcal{S}^*|$

Han [2008] and Kang et al. [2016b] discussed a special case of (IN1-M) where all the instruments are relevant, i.e. $|\mathcal{S}^*| = p_z$. However, as stressed in Kang et al. [2016b], the 50% rule is only a sufficient condition, not a necessary condition, to identify the model parameters. A necessary and sufficient condition for identification, which we state in Theorem 1, is that the valid instruments form a plurality defined by ratios of $\boldsymbol{\pi}^*$ and $\boldsymbol{\gamma}^*$.

THEOREM 1. *Suppose models* (2) *and* (3) *hold,* $\boldsymbol{\Sigma}^*$ *exists and is invertible. The parameters* $\beta^*$ *and* $\boldsymbol{\pi}^*$ *are identified if and only if the following hold*

(IN1-P)  **Plurality rule**: $|\mathcal{V}^*| > \max_{c \neq 0} \left| \left\{ j \in \mathcal{S}^* : \frac{\pi_j^*}{\gamma_j^*} = c \right\} \right|$

As an illustration of Theorem 1, suppose there are $p_z = 5$ instruments that are all relevant, $\gamma_j^* \neq 0, j = 1, \ldots, 5$, and the first two are valid, $\pi_1^* = \pi_2^* = 0$. Under this setup, condition (IN1-M) is not sufficient for identification because $|\mathcal{V}^*| = 2$, which is less than $|\mathcal{S}^*|/2 = 5/2$. Instead, using Theorem 1, identification is possible if the last three invalid instruments have different values of $\pi_3^*/\gamma_3^*, \pi_4^*/\gamma_4^*$ and $\pi_5^*/\gamma_5^*$. For example, if all five instruments have identical $\gamma_j^*$ values, then identification is possible if the last three invalid instruments have different $\pi_j^*$ values because $\max_{c \neq 0} |\{j \in \mathcal{S}^* : \pi_j^*/\gamma_j^* = c\}| = 1$ while $|\mathcal{V}^*| = 2$. In other words, the valid instruments form a plurality defined by $\pi_j^*/\gamma_j^*$ compared to the invalid instruments. Conversely, by the necessary aspect of Theorem 1, if $\pi_3^*/\gamma_3^* = \pi_4^*/\gamma_4^*$, then $\max_{c \neq 0} |\{j \in \mathcal{S}^* : \pi_j^*/\gamma_j^* = c\}| = 2$. Subsequently, the valid instruments cannot form a plurality defined by $\pi_j^*/\gamma_j^*$, condition (IN1-P) fails, and identification is not possible. Note that the same claim cannot be said about the condition (IN1-M) because (IN1-M) is only a sufficient condition.

More generally, condition (IN1-M) implies condition (IN1-P) . This is because if more than 50% of the relevant instruments are valid, the maximum on the right-hand side of (IN1-P) is less than $1/2 \cdot |\mathcal{S}^*|$ because there are less than $1/2 \cdot |\mathcal{S}^*|$ where $\pi_j^* \neq 0$, irrespective of the values of $\gamma_j^*$. In short, condition (IN1-M) is more stringent than condition (IN1-P) .

We also compare two identifying assumptions that are common in the IV literature. First, if there are no relevant instruments, i.e. $|\mathcal{S}^*| = 0$, it is well-known in the literature that $\beta^*$ is not identified. This conclusion parallels condition (IN1-P) where $|\mathcal{S}^*| = 0$ implies that condition (IN1-P) fails and hence, identification is not possible. Second, if there are relevant instruments, i.e. $|\mathcal{S}^*| > 0$, and all the instruments are valid, i.e. $\pi_j^* = 0$ for all $j \in \mathcal{S}^*$, traditional IV literature tells us that $\beta^*$ is identified. Again, we can derive the same conclusion from condition (IN1-P) because $|\mathcal{V}^*| = |\mathcal{S}^*|$, which is the upper bound on the size of the set $\left\{ j \in \mathcal{S}^* \ : \ \pi_j^*/\gamma_j^* = c \right\}$ for any $c \neq 0$ and thus, condition (IN1-P) is satisfied.

Finally, Theorem 1 is different than the the necessary and sufficient condition stated in Kang et al. [2016b] as follows. First, Kang et al. [2016b] only considered the case where all the instruments were relevant, i.e. $\gamma_j^* \neq 0$, whereas Theorem 1 generalizes this to the case where some (or all) instruments may be irrelevant. Second, Kang et al. [2016b] considered identification where the parameter space of $\boldsymbol{\pi}^*$ was restricted by an upper bound $U$ on the number of invalid instruments, i.e. $\|\boldsymbol{\pi}^*\|_0 < U$. This upper bound $U$ can also be thought of as a sensitivity parameter where $U = 1$ represents the best case scenario where all the instruments are valid and deviations from $U = 1$ represent violations of instrument validity [Kang et al., 2016a]. In contrast, Theorem 1 does not consider any upper bound $U$. Nevertheless, the two conditions are equivalent if all instruments are relevant and $U = p_z - |\mathcal{V}^*| + 1$; in this case, the proof in Kang et al. [2016b] can be used to prove Theorem 1.

The identification result in Theorem 1 provides a blueprint for building a confidence interval of $\beta^*$. Specifically, Theorem 1 implies that we need an estimate of IVs that satisfy (A1), i.e. the set $\mathcal{S}^*$, and an estimate of IVs that satisfy (A2) and (A3), i.e. $\mathcal{P}^*$. Additionally, these estimates must satisfy condition (IN1-P) to identify and eventually construct a confidence interval of $\beta^*$. Our method, TSHT with voting, does exactly this. In particular, the first stage of TSHT estimates $\mathcal{S}^*$ and the second stage of TSHT generates many candidate estimates of $\mathcal{P}^*$. The voting step ensures asymptotically that we provide a good estimator of $\mathcal{V}^*$ under condition (IN1-P) .

## 3. Confidence Interval Estimation via Two-Stage Hard Thresholding (TSHT) with Voting

### 3.1. An Illustration of TSHT in Low Dimensional Settings

We first illustrate TSHT under the low dimensional setting where $n \gg p_z + p_x$. The low dimensional setting is common in many applications of instrumental variables, such as economics, social sciences, and medical sciences, including Mendelian Randomization.

We start by rewriting the models of $Y$ and $D$ in equations (2) and (3) into reduced-forms, i.e. models of $Y$ and $D$ are functions of $\mathbf{Z}_{i.}$ and $\mathbf{X}_{i.}$ only.

$$Y_i = \mathbf{Z}_{i.}^{\mathsf{T}}\boldsymbol{\Gamma}^* + \mathbf{X}_{i.}^{\mathsf{T}}\Psi^* + \epsilon_{i1}, \quad \mathbf{E}(\epsilon_{i1} \mid \mathbf{Z}_{i.}, \mathbf{X}_{i.}) = 0 \tag{3}$$

$$D_i = \mathbf{Z}_{i.}^{\mathsf{T}}\boldsymbol{\gamma}^* + \mathbf{X}_{i.}^{\mathsf{T}}\boldsymbol{\psi}^* + \epsilon_{i2}, \quad \mathbf{E}(\epsilon_{i2} \mid \mathbf{Z}_{i.}, \mathbf{X}_{i.}) = 0 \tag{4}$$

where $\boldsymbol{\Gamma}^* = \beta^*\boldsymbol{\gamma}^* + \boldsymbol{\pi}^*$, $\Psi^* = \beta^*\boldsymbol{\psi}^* + \boldsymbol{\phi}^*$, and $\epsilon_{i1} = \beta^*\epsilon_{i2} + e_i$ is the reduced-form error terms. Specifically, the term $\boldsymbol{\Gamma}^*$ represents the intent-to-treat (ITT) effect of the instruments on the outcome and the term $\boldsymbol{\gamma}^*$ represents the effect of the instruments on the treatment. The terms $\epsilon_{i1}$ and $\epsilon_{i2}$ are reduced-form errors with covariance matrix $\mathbf{\Theta}^*$ where $\Theta_{11}^* = \mathrm{Var}(\epsilon_{i1}|\mathbf{Z}_{i.}, \mathbf{X}_{i.})$, $\Theta_{22}^* = \mathrm{Var}(\epsilon_{i2}|\mathbf{Z}_{i.}, \mathbf{X}_{i.})$, and $\Theta_{12}^* = \mathrm{Cov}(\epsilon_{i1}, \epsilon_{i2}|\mathbf{Z}_{i.}, \mathbf{X}_{i.})$.

Each reduced-form model is the usual regression model with regressors $\mathbf{Z}_{i.}$ and $\mathbf{X}_{i.}$ and outcomes $D_i$ and $Y_i$. There are consistent and asymptotically normal estimators of the regression model parameters in low dimensional settings, for instance estimators based on ordinary least squares (OLS) stated below.

$$(\widehat{\boldsymbol{\gamma}}, \widehat{\boldsymbol{\psi}})^{\intercal} = (\mathbf{W}^{\intercal}\mathbf{W})^{-1}\mathbf{W}^{\intercal}\mathbf{D}, \quad (\widehat{\boldsymbol{\Gamma}}, \widehat{\boldsymbol{\Psi}})^{\intercal} = (\mathbf{W}^{\intercal}\mathbf{W})^{-1}\mathbf{W}^{\intercal}\mathbf{Y}$$

$$\widehat{\Theta}_{11} = \frac{1}{n}\left\|\mathbf{Y} - \mathbf{Z}\widehat{\boldsymbol{\Gamma}} - \mathbf{X}\widehat{\boldsymbol{\Psi}}\right\|_2^2, \quad \widehat{\Theta}_{22} = \frac{1}{n}\left\|\mathbf{D} - \mathbf{Z}\widehat{\boldsymbol{\gamma}} - \mathbf{X}\widehat{\boldsymbol{\psi}}\right\|_2^2,$$

$$\widehat{\Theta}_{12} = \frac{1}{n}\left(\mathbf{Y} - \mathbf{Z}\widehat{\boldsymbol{\Gamma}} - \mathbf{X}\widehat{\boldsymbol{\Psi}}\right)^{\intercal}\left(\mathbf{D} - \mathbf{Z}\widehat{\boldsymbol{\gamma}} - \mathbf{X}\widehat{\boldsymbol{\psi}}\right)$$

Let $\widehat{\mathbf{U}}$ denote an estimate of $(\mathbf{\Sigma}^*)^{-1}$, the precision matrix of $\mathbf{W}$, i.e. $\widehat{\mathbf{U}} = (\mathbf{W}^{\intercal}\mathbf{W}/n)^{-1}$. Then, $\Theta_{11}\widehat{\mathbf{U}}/n$ and $\Theta_{22}\widehat{\mathbf{U}}/n$ are the covariance matrices of the OLS estimators $(\widehat{\boldsymbol{\Gamma}}, \widehat{\boldsymbol{\Psi}})$ and $(\widehat{\boldsymbol{\gamma}}, \widehat{\boldsymbol{\psi}})$, respectively.

The estimators above are the only necessary inputs for TSHT: (i) **c**onsistent and **a**symptotically **n**ormal (CAN) estimators of the reduced-form coefficients in equations (3) and (4), (ii) a consistent estimator of the error variance matrix $\mathbf{\Theta}^*$, and (iii) the instrument-covariate matrix $\mathbf{W}$, primarily to estimate its precision matrix, $(\mathbf{\Sigma}^*)^{-1}$. While our discussion was restricted to OLS estimators, any estimator that satisfies the input requirements will work for TSHT . For example, in Section 3.7, we discuss input estimators for TSHT when the data is high dimensional. Finally, we emphasize that no additional choices/inputs are needed for TSHT, such as tuning parameters, beyond the inputs stated above.

### 3.2. First Hard Thresholding: Select Strong IVs Satisfying (A1), $\mathcal{S}^*$

The first thresholding step estimates the set of instruments that satisfy (A1), or the set $\mathcal{S}^* = \{j : \gamma_j^* \neq 0\}$ defined in Definition 2. To do this, we use one of the inputs for TSHT, the estimator for $\boldsymbol{\gamma}^*$, $\widehat{\boldsymbol{\gamma}}$. Specifically, if the $j$th component of $\widehat{\boldsymbol{\gamma}}$ exceeds some threshold, then $j$ mostly likely belongs to the set $\mathcal{S}^*$. Estimating $\mathcal{S}^*$ based on this principle is called hard thresholding [Donoho and Johnstone, 1994, Donoho, 1995] and we denote an estimator of $\mathcal{S}^*$ as $\widehat{\mathcal{S}}$.

$$\widehat{\mathcal{S}} = \left\{ j : |\widehat{\gamma}_j| \geq \frac{\sqrt{\widehat{\Theta}_{22}}\|\mathbf{W}\widehat{\mathbf{U}}_{.j}\|_2}{\sqrt{n}}\sqrt{\frac{2.01\log\max(p_z, n)}{n}} \right\} \tag{5}$$

The threshold to declare whether the estimate $\widehat{\gamma}_j$ is away from zero consists of two terms. The first term $\sqrt{\widehat{\Theta}_{22}}\|\mathbf{W}\widehat{\mathbf{U}}_{.j}\|_2/n$ represents the standard error of $\widehat{\gamma}_j$. The second term $\sqrt{2.01\log\max(p_z, n)}$ represents a multiplicity correction for checking whether normally distributed estimators, like $\widehat{\gamma}_j$, are away from zero. In particular, the $\sqrt{2.01\log(\cdot)}$ part comes from the tail bound of a normal distribution. The $\max(p_z, n)$ part comes from checking multiple $\widehat{\gamma}_j$'s distance from zero. Without the multiplier term $\max(p_z, n)$ in equation (5) and if we have many instruments, some estimates $\widehat{\gamma}_j$ may exceed the threshold by chance and be part of the set $\widehat{\mathcal{S}}$, even though their true $\gamma_j^*$s may actually be zero. In practice, $\max(p_z, n)$ is often replaced by $p_z$ or $n$ to improve the finite-sample performance of hard thresholding procedures and we explore this numerically in Section 5. But, so long as this term grows with $n$, like $\max(p_z, n)$ or $p_z$ that grow with $n$ in high dimensional asymptotics, the asymptotic properties of our procedure described in Section 4 hold.

To summarize, the two terms, $\sqrt{\widehat{\Theta}_{22}} \|\mathbf{W} \widehat{\mathbf{U}}_{.j}\|_2 / n$ and $\sqrt{2.01 \log \max(p_z, n)}$, make up the threshold that account for the variability of the estimate $\widehat{\gamma}_j$ as well as the repeated testing of whether an IV satisfies (A1). If $\widehat{\gamma}_j$ is away from zero by the amount specified by these two term, there is a high probability that the underlying estimand, $\gamma_j^*$, is away from zero and hence, satisfies (A1). Also, the estimator of $\mathcal{S}^*$ does not require selection of tuning parameters, which is in contrast to other variable selection procedures like the Lasso [Tibshirani, 1996] which typically uses cross-validation to select the tuning parameters [Hastie et al., 2016]; all the components of our threshold in equation (5) are pre-determined from the inputs provided in Section 3.1.

If external information suggests that instruments are associated with exposure and these associations are strong, then the first thresholding step may not be necessary and we can simply set $\widehat{\mathcal{S}} = \{1, \ldots, p_z\}$. However, when some of these associations may be weak, we recommend running the first-thresholding step to improve the finite-sample performance of TSHT since the first thresholding should eliminate weak instruments and make TSHT more robust.

### 3.3. Second Hard Thresholding: Select Valid IVs Satisfying (A2) and (A3), $\mathcal{P}^*$

The second thresholding step estimates the set of instruments that satisfy (A2) and (A3), or the set $\mathcal{P}^* = \{j : \pi_j^* \neq 0\}$ defined in Definition 1. Unfortunately, unlike the first thresholding step, none of the inputs for TSHT in Section 3.1 directly estimates $\boldsymbol{\pi}^*$, which we can use to estimate $\mathcal{P}^*$ via hard thresholding. Instead, we propose many estimates of $\mathcal{P}^*$ and combine the information from each estimate via voting. We illustrate the estimation of $\mathcal{P}^*$ in this section and the voting in the next section.

To estimate $\mathcal{P}^*$, we need an estimate of $\boldsymbol{\pi}^*$ that define $\mathcal{P}^*$, which we propose as follows. For each individually strong IV $j \in \widehat{\mathcal{S}}$, we propose a plug-in estimate of $\boldsymbol{\pi}^*$, denoted as $\widehat{\boldsymbol{\pi}}^{[j]}$, based on the relationship between the model parameters $\boldsymbol{\Gamma}^* = \beta^* \boldsymbol{\gamma}^* + \boldsymbol{\pi}^*$ in equation (4)

$$\widehat{\boldsymbol{\pi}}^{[j]} = \widehat{\boldsymbol{\Gamma}} - \frac{\widehat{\Gamma}_j}{\widehat{\gamma}_j} \widehat{\boldsymbol{\gamma}} \tag{6}$$

The terms $\widehat{\boldsymbol{\Gamma}}$ and $\widehat{\boldsymbol{\gamma}}$ in equation (6) are directly from the inputs to TSHT. The term $\widehat{\Gamma}_j / \widehat{\gamma}_j$ in equation (6) is a Wald-type/ratio estimate of $\beta^*$ based on instrument $j$. Note that the estimate of $\beta^*$ comes only from individually strong instruments in the set $\widehat{\mathcal{S}}$. We also propose an estimate of the variance $\sigma^2$ based on this $j$th IV as $\widehat{\sigma^2}^{[j]} = \widehat{\Theta}_{11} + (\widehat{\beta}^{[j]})^2 \widehat{\Theta}_{22} - 2\widehat{\beta}^{[j]}\widehat{\Theta}_{12}$. In total, we should have $|\widehat{\mathcal{S}}|$ estimates of $\boldsymbol{\pi}^*$ and $\sigma^2$.

For each estimate of $\boldsymbol{\pi}^*$, we can estimate the set $\mathcal{P}^*$ in a similar fashion to the first hard thresholding step; the only difference is that we are selecting instruments $k$ with $\pi_k^* = 0$ whereas in the first thresholding step, we are selecting instruments with $\gamma_k^* \neq 0$. Specifically, for each estimate $\boldsymbol{\pi}^{[j]}$, we threshold each component of the vector $\boldsymbol{\pi}^{[j]}$ below some threshold and we denote the set consisting of these components as $\widehat{\mathcal{P}}^{[j]}$.

$$\widehat{\mathcal{P}}^{[j]} =: \left\{ k : |\widehat{\pi}_k^{[j]}| \leq \sqrt{\widehat{\sigma^2}^{[j]}} \frac{\|\mathbf{W}(\widehat{\mathbf{U}}_{\cdot k} - \frac{\widehat{\gamma}_k}{\widehat{\gamma}_j}\widehat{\mathbf{U}}_{\cdot j})\|_2}{\sqrt{n}} \sqrt{\frac{2.01^2 \log \max(p_z, n)}{n}} \right\} \quad (7)$$

Like the first threshold in equation (5), the threshold in equation (7) comprises of two terms. The first term $\sqrt{\widehat{\sigma^2}^{[j]}}\|\mathbf{W}(\widehat{\mathbf{U}}_{\cdot k} - \frac{\widehat{\gamma}_k}{\widehat{\gamma}_j}\widehat{\mathbf{U}}_{\cdot j})\|_2/n$ represents the standard error of $\widehat{\pi}_k^{[j]}$. The second term $\sqrt{2.01^2 \log \max(p_z, n)}$ represents the multiplicity correction. The constant $2.01^2$ is due to the fact that we are performing (at most) $p_z^2$ hypothesis testing for all candidate-component combinations. Combined, the two terms account for the variability of the estimate $\widehat{\boldsymbol{\pi}}^{[j]}$ as well as the multiplicity of testing whether each IV satisfies (A2) and (A3) and in the end, we have $|\widehat{\mathcal{S}}|$ estimates of $\mathcal{P}^*$, $\widehat{\mathcal{P}}^{[j]}, j \in \widehat{\mathcal{S}}$. Also, similar to the first thresholding step, all the components in equation (7) are pre-determined.

To recap, the first thresholding step estimates a set of IVs that satisfy (A1), $\widehat{\mathcal{S}}$. The second thresholding step estimates sets of IVs that satisfy (A2) and (A3), $\widehat{\mathcal{P}}^{[j]}$ for each $j \in \widehat{\mathcal{S}}$. Combining the two thresholding steps gives estimates of IVs that satisfy all (A1)-(A3), or the set $\mathcal{V}^*$ in Definition 3. Specifically, each intersection $\widehat{\mathcal{V}}^{[j]} = \widehat{\mathcal{S}} \cap \widehat{\mathcal{P}}^{[j]}$ is an estimate of $\mathcal{V}^*$ and we have $|\widehat{\mathcal{S}}|$ estimates of $\mathcal{V}^*$ (i.e. $\widehat{\mathcal{V}}^{[j]}, j \in \widehat{\mathcal{S}}$). The remaining task is to combine the information from these estimates in a manner that is consistent with the identification conditions in Section 2.2 to produce a single, consistent estimate of the set $\mathcal{V}^*$.

### 3.4. Majority and Plurality Voting

To explain how we combine several estimates $\widehat{\mathcal{V}}^{[j]}, j \in \widehat{\mathcal{S}}$ to produce a consistent estimate of $\mathcal{V}^*$, it's helpful to consider a voting analogy where each $j \in \widehat{\mathcal{S}}$ is an expert and $\widehat{\mathcal{V}}^{[j]}$ is expert $j$'s ballot that contains expert's $j$'s opinion about which instruments that he/she believes satisfy (A1)-(A3). Because $\widehat{\mathcal{V}}^{[j]} \subseteq \widehat{\mathcal{S}}$ for any $j$, all experts have to pick instruments from the set $\widehat{\mathcal{S}}$ when they cast their ballots. For example, $k \in \widehat{\mathcal{V}}^{[j]}$ indicates that expert $j$ voted on instrument $k$ as satisfying (A1)-(A3). Another expert $j' \neq j$ may have not voted for instrument $k$ as satisfying (A1)-(A3), i.e. $k \notin \widehat{\mathcal{V}}^{[j']}$

Following the voting analogy, we can tally the number of experts who cast their votes for a particular candidate IV as satisfying (A1)-(A3). Specifically, let $\mathbf{1}(k \in \widehat{\mathcal{V}}^{[j]})$ denote the indicator function of the $k$th instrument belonging to $\widehat{\mathcal{V}}^{[j]}$ and $VM_k = \sum_{j \in \widehat{\mathcal{S}}} \mathbf{1}(k \in \widehat{\mathcal{V}}^{[j]})$ denote the tally of votes from experts for each candidate IV $k \in \widehat{\mathcal{S}}$. For example, $VM_k = 3$ indicates that three out of $|\widehat{\mathcal{S}}|$ total experts have voted instrument $k$ as satisfying (A1)-(A3).

Now, suppose a candidate IV $k \in \widehat{\mathcal{S}}$ received votes from a majority of experts, i.e. more than 50% of experts, as satisfying (A1)-(A3), i.e. $VM_k > 1/2 \cdot |\widehat{\mathcal{S}}|$. Let $\widehat{\mathcal{V}}_M$ consist of such instruments and we refer to this type of voting as majority voting,

$$\widehat{\mathcal{V}}_M = \left\{ k \in \widehat{\mathcal{S}} \mid VM_k > \frac{1}{2}|\widehat{\mathcal{S}}| \right\}. \tag{8}$$

For example, suppose there are $|\widehat{\mathcal{S}}| = 5$ experts voting on five candidate IVs labeled from 1 to 5. Each candidate IV received the following total number of votes from the 5 experts: $VM_1 = 4, VM_2 = 4, VM_3 = 4, VM_4 = 1$, and $VM_5 = 1$. Then, $\widehat{\mathcal{V}}_M = \{1, 2, 3\}$ and the first three instruments has won the majority vote. Also, if the identification condition (IN1-M) held so that majority of instruments satisfy (A1)-(A3), then the set $\widehat{\mathcal{V}}_M$ would be similar to the set $\mathcal{V}^*$. This is because majority of experts would produce similar estimates of $\beta^*$, $\pi^*$, and $\mathcal{V}^*$. Specifically, if instrument $k$ truly satisfy (A1)-(A3), then a majority of experts' estimates $\widehat{\mathcal{V}}^{[j]}$ would contain $k$ and the tally of votes for instrument $k$ across experts, i.e. $VM_k$, would exceed 50% of the total number of experts. Consequently, if the identification condition (IN1-M) held, the set $\widehat{\mathcal{V}}_M$ would be a good proxy of $\mathcal{V}^*$.

Now, suppose that the condition (IN1-M) does not hold, then it is highly likely that no instrument won support from a majority of the voters and $\widehat{\mathcal{V}}_M$ is empty. In this case, suppose a candidate IV $k$ received a plurality of votes to satisfy (A1)-(A3), i.e. $VM_k = \max_l VM_l$. Let $\widehat{\mathcal{V}}_P$ denote instruments that received a plurality of votes and we refer to this type of voting as plurality voting,

$$\widehat{\mathcal{V}}_P = \left\{ k \in \widehat{\mathcal{S}} \mid VM_k = \max_l VM_l \right\}. \tag{9}$$

For example, going back to the previous example with 5 experts, suppose the tallies for the 5 candidate instruments are $VM_1 = 2, VM_2 = 1, VM_3 = 2, VM_4 = 1$, and $VM_5 = 1$. Then while none of the five instruments received a majority vote, instruments 1 and 3 received plurality of votes to satisfy (A1)-(A3) and thus, $\widehat{\mathcal{V}}_P = \{1, 3\}$. While weaker than a majority vote, plurality of support from experts shows strong evidence for instruments 1 and 3 satisfying (A1)-(A3). More importantly, if the more general identification condition (IN1-P) held, there would be more experts using valid instruments and producing similar estimates of $\beta^*$, $\pi^*$, and $\mathcal{V}^*$, i.e. the left-hand side of condition (IN1-P), than those using invalid instruments and producing similar estimates of these quantities, i.e. the right-hand side of condition (IN1-P). Thus, if an instrument $k$ truly satisfied (A1)-(A3), those experts using valid instruments to inform their ballots $\widehat{\mathcal{V}}^{[j]}$ would include $k$ in their ballots and the tally for instrument $k$ across all experts, i.e. $VM_k$, would be largest because there are more experts using the valid instruments. Hence, $\widehat{\mathcal{V}}_P$ would be a good proxy of $\mathcal{V}^*$ if condition (IN1-P) held.

To summarize, under the majority rule condition (IN1-M) , $\widehat{\mathcal{V}}_M$ will be sufficient to provide a good estimate of $\mathcal{V}^*$. However, if the majority rule fails but the necessary and sufficient condition (IN1-P) holds, $\widehat{\mathcal{V}}_P$ will be a good estimate of $\mathcal{V}^*$. Then, a single, robust estimate of $\mathcal{V}^*$ under any of the two conditions is the union of the two sets $\widehat{\mathcal{V}} = \widehat{\mathcal{V}}_M \cup \widehat{\mathcal{V}}_P$. Technically speaking, because condition (IN1-P) is both sufficient and necessary, the union can only consist of the set $\widehat{\mathcal{V}}_P$. However, we find that in simulation studies and in practice, taking the union of the two sets provide robustness in finite-samples.

### 3.5. Point Estimate, Standard Error, and Confidence Interval

Once we have an estimate of the set of instruments that satisfy (A1)-(A3), i.e. $\widehat{\mathcal{V}}$, estimation and inference of $\beta^*$ is straightforward in the low dimensional setting. In particular, we can use two stage least squares with $\widehat{\mathcal{V}}$ as the set of IVs that satisfy (A1)-(A3) and obtain a point estimate for $\beta^*$, which we denote as $\widehat{\beta}_L$

$$\widehat{\beta}_L = \frac{\widehat{\gamma}_{\widehat{\mathcal{V}}}^{\mathsf{T}} \widehat{\mathbf{A}} \widehat{\mathbf{\Gamma}}_{\widehat{\mathcal{V}}}}{\widehat{\gamma}_{\widehat{\mathcal{V}}}^{\mathsf{T}} \widehat{\mathbf{A}} \widehat{\gamma}_{\widehat{\mathcal{V}}}}, \quad \widehat{\mathbf{A}} = \widehat{\mathbf{\Sigma}}_{\widehat{\mathcal{V}},\widehat{\mathcal{V}}} - \widehat{\mathbf{\Sigma}}_{\widehat{\mathcal{V}},\widehat{\mathcal{V}}^c} \widehat{\mathbf{\Sigma}}_{\widehat{\mathcal{V}}^c,\widehat{\mathcal{V}}^c}^{-1} \widehat{\mathbf{\Sigma}}_{\widehat{\mathcal{V}}^c,\widehat{\mathcal{V}}} \tag{10}$$

The $\widehat{\mathbf{A}}$ is a weighing matrix for the estimates $\widehat{\gamma}$ and $\widehat{\mathbf{\Gamma}}$, which, among other things, comprises of $\widehat{\mathbf{\Sigma}} = \mathbf{W}^{\mathsf{T}}\mathbf{W}/n$, the inverse of the estimate of the precision matrix of $\mathbf{W}$ that we used in the inputs for TSHT. The estimated variance of $\widehat{\beta}_L$ is

$$\widehat{\mathrm{Var}}_L = \frac{\widehat{\gamma}_{\widehat{\mathcal{V}}}^{\mathsf{T}} \widehat{\mathbf{A}} \left(\widehat{\mathbf{\Sigma}}^{-1}\right)_{\widehat{\mathcal{V}},\widehat{\mathcal{V}}} \widehat{\mathbf{A}} \widehat{\gamma}_{\widehat{\mathcal{V}}}}{\left(\widehat{\gamma}_{\widehat{\mathcal{V}}}^{\mathsf{T}} \widehat{\mathbf{A}} \widehat{\gamma}_{\widehat{\mathcal{V}}}\right)^2} \left(\widehat{\Theta}_{11} + \widehat{\beta}_L^2 \widehat{\Theta}_{22} - 2\widehat{\beta}_L \widehat{\Theta}_{12}\right) \tag{11}$$

which simplifies to

$$\widehat{\mathrm{Var}}_L = \frac{\widehat{\Theta}_{11} + \widehat{\beta}_L^2 \widehat{\Theta}_{22} - 2\widehat{\beta}_L \widehat{\Theta}_{12}}{\widehat{\gamma}_{\widehat{\mathcal{V}}}^{\mathsf{T}} \widehat{\mathbf{A}} \widehat{\gamma}_{\widehat{\mathcal{V}}}}.$$

Finally, for any $\alpha$ with $0 < \alpha < 1$, the $1 - \alpha$ confidence interval for $\beta^*$ is

$$\left(\widehat{\beta}_L - z_{1-\alpha/2}\sqrt{\widehat{\mathrm{Var}}_L/n}, \quad \widehat{\beta}_L + z_{1-\alpha/2}\sqrt{\widehat{\mathrm{Var}}_L/n}\right). \tag{12}$$

where $z_{1-\alpha/2}$ is the $1 - \alpha/2$ quantile of the standard normal distribution.

In Section 4.1, we show that the $\widehat{\beta}_L$ achieves optimal performance in the sense that $\widehat{\beta}_L$ converges to an asymptotic normal distribution that is identical to the asymptotic normal distribution of the TSLS estimator for $\beta^*$ that knows which IVs are valid a priori, i.e. the set $\mathcal{V}^*$.

### 3.6. Comparison to Other Methods

We make some remarks about our method and the methods proposed in the literature on invalid IVs. The work by Windmeijer et al. [2016] is the methodologically most similar to our method in that it also estimates $\mathcal{V}^*$ and uses the estimate of $\mathcal{V}^*$ to obtain oracle-optimal point estimate and confidence interval of $\beta^*$ like we do, i.e. via two-stage least squares as outlined in Section 3.5. To estimate $\mathcal{V}^*$, Windmeijer et al. [2016] utilizes the adaptive Lasso with a median estimator of Han [2008] and Bowden et al. [2016] as the initial estimator; the tuning parameter in the adaptive Lasso is chosen via data by cross-validation. In contrast, TSHT with voting utilizes hard thresholding steps to estimate $\mathcal{V}^*$ where our 'tuning' parameters, i.e. the thresholds, are pre-determined and theoretically motivated. Based on the numerical results in Section 5.2, we suspect that in low dimensional settings, their method and TSHT with voting are asymptotically equivalent when condition (IN1-M) holds.

Another inferential method in the invalid IV literature is bootstrapping the median estimator [Bowden et al., 2016, Burgess et al., 2016]. The key idea is to directly go after the target estimand, $\beta^*$, with the median estimator mentioned above and bootstrap the estimate with sufficient statistics. Their works are under the two-sample designs with summary data where the errors in the reduced-form models are independent of each other. In contrast, TSHT with voting and Windmeijer et al. [2016]'s method focus on the one-sample design with individual level data and correlated error terms. Also, both methods do not utilize the bootstrap to generate our inferential quantities.

We argue that TSHT with voting is a major improvement from the methods of Windmeijer et al. [2016], Bowden et al. [2016] and Burgess et al. [2016] for the following three reasons. First, all three methods rely on the 50% rule because their estimators rely on the median estimator, which is only consistent whenever the 50% rule holds. In contrast, TSHT with voting does not rely on an initial consistent estimator and our inferential guarantees are possible under the more general (IN1-P) condition. Second, the median methods of Bowden et al. [2016] and Burgess et al. [2016] may not be oracle-optimal in the sense that it may not be as efficient as the oracle estimator that knows, a priori, which instruments are invalid. Windmeijer et al. [2016]'s method is oracle-optimal in low dimensional settings, but only when the 50% rule holds. In contrast, TSHT with voting is oracle-optimal in low dimensional settings under the more general (IN1-P) condition; see Section 4.1. Third, there are no theoretical guarantees that the bootstrap approach to inference for the median method will always generate a confidence interval that will cover the true parameter with probability $1 - \alpha$, although it does perform well in large numerical studies under two-sample designs [Burgess et al., 2016]. Similarly, the theoretical properties of Windmeijer et al. [2016]'s method are under the assumption that the tuning parameter is not chosen via cross-validation, despite the fact that Windmeijer et al. [2016] utilizes cross-validation when they use their method in simulations and in real data example. In contrast, TSHT with voting uses pre-determined thresholding values both in theory and in numerical studies and have theoretical guarantees on inference; see Section 4.

Finally, the work by Kang et al. [2016b] which is the precursor of this paper, also proposes a joint estimator of $\beta^*$ and $\boldsymbol{\pi}^*$ called sisVIVE. Kang et al. [2016b]'s estimator is based on the Lasso that minimizes the sum of squared errors from the model (2) with respect to an $\ell_1$ penalty on $\boldsymbol{\pi}^*$. The tuning parameter of the Lasso is chosen via cross-validation. A nice feature of sisVIVE is that it is a one-step method to estimate $\beta^*$. In contrast, TSHT requires two thresholding steps plus voting to estimate $\beta^*$. Unfortunately, sisVIVE requires more stringent conditions for consistency than the identification condition (IN1-M). Also, like the method of Windmeijer et al. [2016], the theory behind consistency is developed under the assumption that the tuning parameter is not chosen via cross-validation. More importantly, sisVIVE did not resolve the issue of confidence interval construction.

### 3.7. High Dimensional Setting

TSHT with voting can also accommodate settings where we have high dimensional covariates and/or instruments. The modifications we have to make are the estimation of the reduced-form model parameters in equations (3) (4), the weighing matrix $A$ in (10), and the formula for the standard error; the rest of the procedure is identical.

Specifically, instead of using OLS estimators in Section 3.1, we have to resort to estimators that can handle the csae when $n \ll p$ and are CAN so that the input requirements for TSHT are met. There are many estimators in high dimensions that meet this criterion, such as the debiased Lasso or its variants laid out in Zhang and Zhang [2014], Javanmard and Montanari [2014], van de Geer et al. [2014] and Cai and Guo [2017]. For completeness, we present one estimator in high dimensional regression that is CAN, the debiased square-root Lasso estimator [Belloni et al., 2011, Javanmard and Montanari, 2014]; see the references cited for additional details on CAN estimators in high dimensions. First, the square-root Lasso estimator [Belloni et al., 2011] estimates high dimensional reduced-form model parameters in equations (3) and (4) based on the following optimization problems.

$$\{\widetilde{\boldsymbol{\Gamma}}, \widetilde{\boldsymbol{\Psi}}\} = \operatorname*{argmin}_{\boldsymbol{\Gamma} \in \mathbb{R}^{p_z}, \boldsymbol{\Psi} \in \mathbb{R}^{p_x}} \frac{\|\mathbf{Y} - \mathbf{Z}\boldsymbol{\Gamma} - \mathbf{X}\boldsymbol{\Psi}\|_2}{\sqrt{n}} + \frac{\sqrt{2.01 \log p}}{n} \left( \sum_{j=1}^{p_z} \|\mathbf{Z}_{.j}\|_2 |\Gamma_j| + \sum_{j=1}^{p_x} \|\mathbf{X}_{.j}\|_2 |\Psi_j| \right)$$

$$\{\widetilde{\boldsymbol{\gamma}}, \widetilde{\boldsymbol{\psi}}\} = \operatorname*{argmin}_{\boldsymbol{\Gamma} \in \mathbb{R}^{p_z}, \boldsymbol{\Psi} \in \mathbb{R}^{p_x}} \frac{\|\mathbf{D} - \mathbf{Z}\boldsymbol{\gamma} - \mathbf{X}\boldsymbol{\psi}\|_2}{\sqrt{n}} + \frac{\sqrt{2.01 \log p}}{n} \left( \sum_{j=1}^{p_z} \|\mathbf{Z}_{.j}\|_2 |\gamma_j| + \sum_{j=1}^{p_x} \|\mathbf{X}_{.j}\|_2 |\psi_j| \right)$$

Also, the corresponding estimates of the variances $\Theta_{11}^*, \Theta_{22}^*$, and $\Theta_{12}^*$ from the square-root Lasso are

$$\widehat{\Theta}_{11} = \frac{1}{n} \left\| \mathbf{Y} - \mathbf{Z}\widetilde{\boldsymbol{\Gamma}} - \mathbf{X}\widetilde{\boldsymbol{\Psi}} \right\|_2^2, \quad \widehat{\Theta}_{22} = \frac{1}{n} \left\| \mathbf{D} - \mathbf{Z}\widetilde{\boldsymbol{\gamma}} - \mathbf{X}\widetilde{\boldsymbol{\psi}} \right\|_2^2, \widehat{\Theta}_{12} = \frac{1}{n} \left( \mathbf{Y} - \mathbf{Z}\widetilde{\boldsymbol{\Gamma}} - \mathbf{X}\widetilde{\boldsymbol{\Psi}} \right)^{\mathsf{T}} \left( \mathbf{D} - \mathbf{Z}\widetilde{\boldsymbol{\gamma}} - \mathbf{X}\widetilde{\boldsymbol{\psi}} \right).$$

Unfortunately, the square-root Lasso estimator is biased due to the penalty term and Javanmard and Montanari [2014] proposes a way to debias the square-root Lasso estimator and turn it into CAN estimators. Specifically, Javanmard and Montanari [2014] proposes $p_z$ optimization problems where the solution to each $p_z$ optimization problem, denoted as $\widehat{\mathbf{U}}_{.j} \in \mathbb{R}^p$, $j = 1, \ldots, p_z$, is

$$\widehat{\mathbf{U}}_{.j} = \operatorname*{argmin}_{\mathbf{u} \in \mathbb{R}^p} \frac{1}{n} \|\mathbf{W}\mathbf{u}\|_2^2 \quad \text{s.t.} \quad \|\widehat{\boldsymbol{\Sigma}}\mathbf{u} - \mathbf{I}_{.j}\|_\infty \leq 12 M_1^2 \sqrt{\frac{\log p}{n}} \qquad (13)$$

with $\widehat{\boldsymbol{\Sigma}} = \frac{1}{n}\mathbf{W}^{\mathsf{T}}\mathbf{W}$. Here, $\mathbf{I}_{\cdot j}$ denotes the $j$-th column of the identity matrix $\mathbf{I}$ and $M_1$ denotes the largest eigenvalue of $\boldsymbol{\Sigma}^*$. Let $\widehat{\mathbf{U}}$ denote the matrix concatenation of the $p_z$ solutions to the optimization problem. Then, the debiased estimates of $\widetilde{\boldsymbol{\Gamma}}$ and $\widetilde{\boldsymbol{\gamma}}$, denoted as $\widehat{\boldsymbol{\Gamma}}$ and $\widehat{\boldsymbol{\gamma}}$, are

$$\widehat{\boldsymbol{\Gamma}} = \widetilde{\boldsymbol{\Gamma}} + \frac{1}{n}\widehat{\mathbf{U}}\mathbf{W}^{\mathsf{T}}\left(\mathbf{Y} - \mathbf{Z}\widetilde{\boldsymbol{\Gamma}} - \mathbf{X}\widetilde{\boldsymbol{\Psi}}\right), \quad \widehat{\boldsymbol{\gamma}} = \widetilde{\boldsymbol{\gamma}} + \frac{1}{n}\widehat{\mathbf{U}}\mathbf{W}^{\mathsf{T}}\left(\mathbf{D} - \mathbf{Z}\widetilde{\boldsymbol{\gamma}} - \mathbf{X}\widetilde{\psi}\right). \quad (14)$$

We now have obtained all the ingredients for TSHT for the high dimensional setting: (i) the CAN estimators of $\boldsymbol{\Gamma}^*$ and $\boldsymbol{\gamma}^*$, $\widehat{\boldsymbol{\Gamma}}$ and $\widehat{\boldsymbol{\gamma}}$, respectively, based on the debiased square-root Lasso; (ii) consistent estimators of the error variances $\Theta_{11}^*, \Theta_{22}^*$, and $\Theta_{12}^*$, $\widehat{\Theta}_{11}, \widehat{\Theta}_{22}$ and $\widehat{\Theta}_{12}$ respectively from the square-root Lasso; and (iii) an estimate of the precision matrix of $\mathbf{W}$, $\widehat{\mathbf{U}}$ from the de-biasing procedure. Running TSHT with these inputs will estimate the set of valid instruments $\widehat{\mathcal{V}}$ in high dimensional settings.

For point estimation of $\beta^*$ in high dimensions, we simply replace $\widehat{\mathbf{A}}$ in equation (10) with the identity matrix

$$\widehat{\beta} = \frac{\widehat{\boldsymbol{\gamma}}_{\widehat{\mathcal{V}}}^{\mathsf{T}}\widehat{\boldsymbol{\Gamma}}_{\widehat{\mathcal{V}}}}{\widehat{\boldsymbol{\gamma}}_{\widehat{\mathcal{V}}}^{\mathsf{T}}\widehat{\boldsymbol{\gamma}}_{\widehat{\mathcal{V}}}}. \quad (15)$$

The variance estimate of the point estimate in equation (15) uses a high dimensional estimate of the precision matrix in equation (11), i.e.

$$\widehat{\mathrm{Var}} = \frac{\widehat{\boldsymbol{\gamma}}_{\widehat{\mathcal{V}}}^{\mathsf{T}}(\widehat{\mathbf{U}}_{\cdot,\widehat{\mathcal{V}}})^{\mathsf{T}}\widehat{\mathbf{U}}_{\cdot,\widehat{\mathcal{V}}}\frac{\mathbf{W}^{\mathsf{T}}\mathbf{W}}{n}\widehat{\mathbf{U}}_{\cdot,\widehat{\mathcal{V}}}\widehat{\boldsymbol{\gamma}}_{\widehat{\mathcal{V}}}}{\left(\widehat{\boldsymbol{\gamma}}_{\widehat{\mathcal{V}}}^{\mathsf{T}}\widehat{\boldsymbol{\gamma}}_{\widehat{\mathcal{V}}}\right)^2}\left(\widehat{\Theta}_{11} + \widehat{\beta}^2\widehat{\Theta}_{22} - 2\widehat{\beta}\widehat{\Theta}_{12}\right) \quad (16)$$

Given the point estimate and the variance, the confidence interval for $\beta^*$ follows the usual form

$$\left(\widehat{\beta} - z_{1-\alpha/2}\sqrt{\widehat{\mathrm{Var}}/n}, \quad \widehat{\beta} + z_{1-\alpha/2}\sqrt{\widehat{\mathrm{Var}}/n}\right), \quad (17)$$

## 4. Theoretical Results

In this section, we state the asymptotic properties of TSHT with voting. In Section 4.1, we consider the low dimensional setting where $p_{\mathrm{x}}$ and $p_{\mathrm{z}}$ are fixed. In Section 4.2, we consider the general case when $p_{\mathrm{x}}$ and/or $p_{\mathrm{z}}$ are allowed to grow and exceed sample size $n$.

### 4.1. Invalid IVs in Low Dimensional Setting

First, we prove that the estimated set $\widehat{\mathcal{V}}$ is an asymptotically consistent estimator of the true set $\mathcal{V}^*$ in the low dimensional setting where $p_{\mathrm{x}}$ and $p_{\mathrm{z}}$ are fixed.

LEMMA 1. *Suppose the assumption* (IN1-P) *holds. Then, we have*

$$\lim_{n\to\infty}\mathbf{P}\left(\widehat{\mathcal{V}} = \mathcal{V}^*\right) = 1 \quad (18)$$

Lemma 1 confirms our intuition in Section 3.4 that the voting process correctly generates the set of instruments that are relevant and valid. In fact, a useful feature of our method is that it provably and correctly selects the IVs that satisfy (A1)-(A3), something that is not possible with prior methods that only target $\beta^*$, e.g., the median method.

Next, Theorem 2 states that the confidence interval outlined in Section 3.5 has the desired coverage and optimal length in the low dimensional settings with fixed $p_x$ and $p_z$.

THEOREM 2. *Suppose the assumption* (IN1-P) *holds. Then, as $n \to \infty$, we have*

$$\sqrt{n}\left(\widehat{\beta}_L - \beta^*\right) \xrightarrow{d} N\left(0, \frac{\sigma^2}{\gamma_{\mathcal{V}^*}^{*\intercal}\left(\mathbf{\Sigma}_{\mathcal{V}^*\mathcal{V}^*}^* - \mathbf{\Sigma}_{\mathcal{V}^*(\mathcal{V}^*)^c}^* \mathbf{\Sigma}_{(\mathcal{V}^*)^c(\mathcal{V}^*)^c}^{*-1} \mathbf{\Sigma}_{(\mathcal{V}^*)^c\mathcal{V}^*}^*\right)\gamma_{\mathcal{V}^*}^*}\right).$$
(19)

*Consequently, the confidence interval given in equation* (12) *has asymptotic coverage probability of $1 - \alpha$, i.e.*

$$\mathbf{P}\left\{\beta^* \in \left(\widehat{\beta}_L - z_{1-\alpha/2}\sqrt{\widehat{\mathrm{Var}}_L/n}, \quad \widehat{\beta}_L + z_{1-\alpha/2}\sqrt{\widehat{\mathrm{Var}}_L/n}\right)\right\} \to 1 - \alpha.$$
(20)

We note that the proposed estimator $\widehat{\beta}_L$ has the same asymptotic variance as the oracle TSLS estimator with prior knowledge of $\mathcal{V}^*$, which is efficient under the homoskedastic variance assumption; see Theorem 5.2 in Wooldridge [2010] for details. Consequently, our confidence interval in equation (12) asymptotically performs like the oracle TSLS confidence interval and is of optimal length. But, unlike TSLS, we achieve this oracle performance without prior knowledge of $\mathcal{V}^*$. We remind readers that the previous estimators proposed by Bowden et al. [2015, 2016], Burgess et al. [2016], Kang et al. [2016a] do not achieve oracle performance and TSLS-like efficiency while the estimator proposed by Windmeijer et al. [2016] does achieve this, but only when condition (IN1-M) holds.

### 4.2. Invalid IVs in High Dimensional Setting

We now consider the asymptotic properties of TSHT with voting under the general case when $p_z$ and/or $p_x$ are allowed to grow, potentially exceeding sample size $n$. As noted in Section 3.2, to be in alignment with the traditional high dimensional literature where $p_z$ and/or $p_x$ are always larger than $n$ and growing faster than $n$, we simplify TSHT by replacing the thresholds in equations (5) and (7) from $\log \max\{p_z, n\}$ to $\log p_z$.

We first introduce the regularity assumptions used in high-dimensional statistics [Bickel et al., 2009, Bühlmann and van de Geer, 2011, Cai and Guo, 2017].

(R1) (Coherence): The matrix $\mathbf{\Sigma}^*$ satisfies $1/M_1 \leq \lambda_{\min}(\mathbf{\Sigma}^*) \leq \lambda_{\max}(\mathbf{\Sigma}^*) \leq M_1$ for some constant $M_1 \geq 1$ and has bounded sub-Gaussian norm.

(R2) (normality): The error terms in (3) and (4) follow a bivariate normal distribution.

(R3) (Global IV Strength): The IVs are globally strong with $\sqrt{(\gamma_{\mathcal{V}^*}^*)^\intercal \Sigma_{\mathcal{V}^*,\mathcal{V}^*} \gamma_{\mathcal{V}^*}^*} \gg s_{z1} \log p/\sqrt{n}$, where $\mathcal{V}^*$ is the set of valid IVs defined in Definition 3.

Assumption (R1) places a condition on the spectrum of the design matrix $\mathbf{W}$ and the tail distribution of $\mathbf{W}_{i\cdot}$, which is related to the restricted eigenvalue condition in Bickel et al. [2009]. For simplicity, we also assume that the sub-Gaussian norm of $\mathbf{W}_{i\cdot}$ is upper bounded by $M_1$, that is, $\sup_{\boldsymbol{v} \in S^{p-1}} \sup_{q \geq 1} (\mathbf{E}|\boldsymbol{v}^\intercal \mathbf{W}_{i\cdot}|^q/q)^{\frac{1}{q}} \leq M_1$ where $S^{p-1}$ is the unit sphere in $\mathbb{R}^p$; see Vershynin [2012] for details on sub-Gaussian random variables and bounds. Assumption (R2) states that the errors $(e_{i1}, \epsilon_{i2})$ are bivariate normal. We make the normality assumption out of simplicity, similar to the normal error assumption made in the work on inference in weak IV literature (e.g. Section 2 of Moreira [2003] and Section 2.2.1 Andrews et al. [2007]) and the work on inference in high dimensional linear models (e.g. Theorem 2.5 in Javanmard and Montanari [2014] and Theorem 2.1 in van de Geer et al. [2014]). Finally, assumption (R3) states that the global strength of the instruments, measured by the weighted $\ell_2$ norm of $\gamma_{\mathcal{V}^*}^*$, is bounded away from zero. This type of global strength assumption is commonly made in the IV literature under the guise of a concentration parameter, which is a measure of instrument strength and is the weighted $\ell_2$ norm of $\gamma_{\mathcal{V}^*}^*$; this type of assumption is sometimes referred to as "traditional/strong" asymptotics [Stock et al., 2002, Wooldridge, 2010]. Recent works by Belloni et al. [2012] and Chernozhukov et al. [2015], which considered the setting where all IVs were valid after conditioning on high dimensional covariates, also make this type of assumption, specifically condition SM in Belloni et al. [2012] and condition RF in the supplementary materials of Chernozhukov et al. [2015]. Essentially, both of these works require $\|\gamma^*\|_2$ to be bounded away from zero by a constant and are actually stronger than our assumption (R3). In practice, assumption (R3) is satisfied so long as there is at least one IV that has a constant non-zero effect on the treatment, or a non-zero effect that doesn't diminish with sample size. However, if the IVs are arbitrary weak in the sense of Staiger and Stock [1997], then assumption (R3), let alone the assumptions of Belloni et al. [2012] in high dimensional valid IV literature, does not hold, and we leave this as a future topic of research to deal with arbitrary weak IVs in invalid IV settings.

Section A in the supplementary materials shows that if the IVs are valid, then Assumptions (R1)-(R3) are sufficient to show that a valid confidence interval can be constructed in high dimensional setting. However, when IVs are invalid, we need to make two additional assumptions that are not part of high dimensional statistics or instrumental variables literature and may be of theoretical interest in future work.

(IN2) (Individual IV Strength) For IVs in $\mathcal{S}^*$, $\delta_{\min} = \min_{j \in \mathcal{S}^*} |\gamma_j^*| \gg \sqrt{\log p_z/n}$.

(IN3) (Separated Levels of Violation) For the pair $j, k \in \mathcal{S}^*$ with $\frac{\pi_j^*}{\gamma_j^*} \neq \frac{\pi_k^*}{\gamma_k^*}$,

$$\left| \frac{\pi_j^*}{\gamma_j^*} - \frac{\pi_k^*}{\gamma_k^*} \right| \geq \frac{12 \left( 1 + \max_{j \in \mathcal{S}^*} \left| \frac{\Gamma_j^*}{\gamma_j^*} \right| \right)}{\delta_{\min}} \sqrt{\frac{M_1 \log p_z}{\lambda_{\min}(\boldsymbol{\Theta}^*) n}}. \tag{21}$$

Assumption (IN2) requires individual IV strength to be bounded away from zero so that all IVs in selected $\widehat{\mathcal{S}}$ are strong. Without this condition, an individually weak IV with small $\widehat{\gamma}_j$ may be included in the first thresholding step and subsequently cause trouble to the second thresholding step in equation (7) that uses $\widehat{\gamma}_j$ in the denominator to construct a candidate estimate of $\boldsymbol{\pi}^*$ and $\mathcal{P}^*$. In the literature, assumption (IN2) is similar to the "beta-min" assumption in high dimensional linear regression without IVs, with the exception that this condition is not imposed on our inferential quantity of interest, $\beta^*$. Also, (IN2) is different from Assumption (R3) in that (R3) only requires the global IV strength to be bounded away from zero. Next, Assumption (IN3) requires that the difference between different levels of ratios $\pi_j^*/\gamma_j^*$ is sufficiently large. Without this assumption, it would be difficult to distinguish subsets of instruments with different $\pi_j^*/\gamma_j^*$ values from the data and to identify the set of valid IVs based on the plurality rule (IN1-P). For example, consider instruments $k$ and $l$ with $\pi_k^*/\gamma_k^* \neq \pi_j^*/\gamma_j^*$. If equation (21) is satisfied, then $k \notin \widehat{\mathcal{P}}^{[j]}$ with high probability because $\pi_k^*/\gamma_k^*$ and $\pi_j^*/\gamma_j^*$ are far apart from each other. In contrast, if equation (21) does not hold, then $\widehat{\mathcal{P}}^{[j]}$ might contain instrument $k$ by chance because $\pi_k^*/\gamma_k^*$ and $\pi_j^*/\gamma_j^*$ are close to each other.

Lemma 2 shows that with assumptions (R1)-(R3) and (IN1-P), (IN2) and (IN3), TSHT with voting produces a consistent estimator of the set of valid instruments in the high dimensional setting.

LEMMA 2. *Suppose $\sqrt{s_{z1}} s \log p / \sqrt{n} \to 0$ and assumptions (R1)−(R3), (IN1−P) and (IN2)−(IN3) are satisfied. With probability larger than $1 - c \left( p^{-c} + \exp(-cn) \right)$ for some $c > 0$,*

$$\widehat{\mathcal{V}} = \mathcal{V}^* \tag{22}$$

Next, the following theorem shows that $\widehat{\beta}$ is a consistent and asymptotic normal estimator of $\beta^*$.

THEOREM 3. *Under the same assumption as Lemma 2, we have*

$$\sqrt{n} \left( \widehat{\beta} - \beta^* \right) = T^{\beta^*} + \Delta^{\beta^*} \tag{23}$$

*where $T^{\beta^*} \mid \mathbf{W} \sim N(0, \mathrm{Var})$ and $\mathrm{Var} = \sigma^2 \boldsymbol{\gamma}_{\mathcal{V}^*}^{\mathsf{T}} (\widehat{\mathbf{U}}_{\cdot, \mathcal{V}^*})^{\mathsf{T}} \mathbf{W}^{\mathsf{T}} \mathbf{W} \widehat{\mathbf{U}}_{\cdot, \mathcal{V}^*} \boldsymbol{\gamma}_{\mathcal{V}^*} / n \left( \boldsymbol{\gamma}_{\mathcal{V}^*}^{\mathsf{T}} \boldsymbol{\gamma}_{\mathcal{V}^*} \right)^2$. As $\sqrt{s_{z1}} s \log p / \sqrt{n} \to 0$, $\Delta^{\beta^*} / \sqrt{\mathrm{Var}} \xrightarrow{p} 0$ and the confidence interval given in equation (17) has asymptotic coverage probability of $1 - \alpha$, i.e.*

$$\mathbf{P} \left\{ \beta^* \in \left( \widehat{\beta} - z_{1-\alpha/2} \sqrt{\widehat{\mathrm{Var}}/n}, \quad \widehat{\beta} + z_{1-\alpha/2} \sqrt{\widehat{\mathrm{Var}}/n} \right) \right\} \to 1 - \alpha. \tag{24}$$

## 5. Simulation

### 5.1. Setup: Low Dimensional Setting

In addition to the theoretical analysis of our method in Section 4, we also conduct a simulation study to investigate the numerical performance of our method. The design of the simulation study follows closely that of Windmeijer et al. [2016] where we use the models (2) and (3) in Section 2.1. Specifically, (i) there are $p_z = 7$ or 10 instruments, (ii) there are no covariates, (iii) the instruments are generated from a multivariate normal with mean zero and identity covariance, (iv) the treatment effect is fixed to be $\beta^* = 1$, and (v) the errors have variance 1 and covariance 0.25. Similar to Windmeijer et al. [2016], we vary (i) the sample size $n$, (ii) the strength of the IV by manipulating $\boldsymbol{\gamma}^* = (1, \ldots, 1) \cdot C_\gamma$ with different values of $C_\gamma$, and the (iii) the degree of violations of (A2) and (A3) by manipulating $\boldsymbol{\pi}^*$. With respect to the last variation, if $p_z = 7$, we set $\boldsymbol{\pi}^* = (1, 1, 0.5, 0.5, 0, 0, 0) \cdot C_\pi$ where $C_\pi$ is a constant that we vary to change the magnitude of $\boldsymbol{\pi}^*$. If $p_z = 10$, we set $\boldsymbol{\pi}^* = (1, 1, 1, 1, 0, \ldots, 0) \cdot C_\pi$. The first setting mimics the case where (IN1-M) holds, similar to Windmeijer et al. [2016] while the second setting mimics the case where the 50% rule (IN1-M) fails but (IN1-P) holds .

Under this data generating mechanism, we compare our procedure TSHT with voting to (i) the naive TSLS that assumes all the instruments satisfy (A1)-(A3), (ii) the oracle TSLS that knows, a priori, which instruments satisfy (A1)-(A3), (iii) the method of Windmeijer et al. [2016] that uses the adaptive Lasso tuned via cross-validation and the median estimator, and (iv) the median estimator of Bowden et al. [2016] and Burgess et al. [2016] with bootstrapped confidence intervals by using the R package [Yavorska and Burgess, 2017] under default settings. For (i), we implement TSLS so that it mimics most practitioners' use of TSLS by simply assuming all the instruments $\mathbf{Z}$ are valid. For (ii), we have the oracle TSLS where an oracle provides us with the true set of valid IVs, which will not occur in practice. Note that because TSLS is not robust against weak instruments, we purposely set our $C_\gamma$ to correspond to strong IV regimes. Finally, for (iii) and (iv), see Section 3.6 for discussions of methods (iii) and (iv). Our simulations are repeated 500 times and we measure the median absolute error, the empirical coverage proportion, and average length of the confidence interval computed across simulations.

### 5.2. Result: Low Dimensional Setting

We first present the setting where the 50% rule (IN1-M) holds. Specifically, following Windmeijer et al. [2016], Table 1 shows the cases where we have 10 IVs with $s_{z2} = 3$, $n = (500, 1000, 2000, 5000, 10000)$, $C_\gamma = (0.2, 0.6, 1)$, and $C_\pi = 0.2$. For reference, with $n = 500, 2000, 5000$ and $C_\gamma = 0.2$, the expected concentration parameter is $7nC_\gamma^2$, or, for each $n$, $140, 560$ and $1400$, respectively. Because condition (IN1-M) holds, TSHT, Windmeijer et al. [2016]'s method, and the median method should do well. Indeed, between TSHT and Windmeijer et al. [2016]'s method, there is little difference in terms of median absolute error, coverage, and length of the confidence interval. Both of the methods struggle with low sample size at $n = 500$, but once $n \geq 2000$, the two methods perform as well as the oracle.

| $n$ | $C_\gamma$ | TSHT | | | Adaptive Lasso | | | Median | | | Naive TSLS | | | Oracle TSLS | | |
|---|---|---|---|---|---|---|---|---|---|---|---|---|---|---|---|---|
| 500 | 0.2 | 0.09 | 0.72 | 0.32 | 0.09 | 0.72 | 0.32 | 0.13 | 0.17 | 0.09 | 0.30 | 0.03 | 0.28 | 0.08 | 0.96 | 0.44 |
| 500 | 0.6 | 0.06 | 0.84 | 0.11 | 0.03 | 0.81 | 0.11 | 0.05 | 0.31 | 0.06 | 0.10 | 0.01 | 0.09 | 0.03 | 0.96 | 0.15 |
| 500 | 1.0 | 0.02 | 0.83 | 0.07 | 0.02 | 0.78 | 0.07 | 0.03 | 0.71 | 0.07 | 0.06 | 0.02 | 0.06 | 0.01 | 0.95 | 0.09 |
| 1000 | 0.2 | 0.04 | 0.93 | 0.24 | 0.04 | 0.91 | 0.23 | 0.09 | 0.09 | 0.05 | 0.30 | 0.00 | 0.20 | 0.05 | 0.94 | 0.31 |
| 1000 | 0.6 | 0.01 | 0.95 | 0.08 | 0.01 | 0.93 | 0.08 | 0.03 | 0.27 | 0.03 | 0.10 | 0.00 | 0.07 | 0.02 | 0.96 | 0.10 |
| 1000 | 1.0 | 0.01 | 0.94 | 0.05 | 0.01 | 0.94 | 0.05 | 0.02 | 0.49 | 0.04 | 0.06 | 0.00 | 0.04 | 0.01 | 0.96 | 0.06 |
| 2000 | 0.2 | 0.03 | 0.93 | 0.17 | 0.03 | 0.93 | 0.17 | 0.07 | 0.08 | 0.02 | 0.30 | 0.00 | 0.14 | 0.04 | 0.94 | 0.22 |
| 2000 | 0.6 | 0.01 | 0.96 | 0.06 | 0.01 | 0.96 | 0.06 | 0.02 | 0.16 | 0.02 | 0.10 | 0.00 | 0.05 | 0.01 | 0.96 | 0.07 |
| 2000 | 1.0 | 0.01 | 0.95 | 0.03 | 0.01 | 0.95 | 0.03 | 0.01 | 0.33 | 0.02 | 0.06 | 0.00 | 0.03 | 0.01 | 0.97 | 0.04 |
| 5000 | 0.2 | 0.02 | 0.96 | 0.11 | 0.02 | 0.95 | 0.10 | 0.04 | 0.06 | 0.01 | 0.30 | 0.00 | 0.09 | 0.03 | 0.95 | 0.14 |
| 5000 | 0.6 | 0.01 | 0.96 | 0.04 | 0.01 | 0.96 | 0.03 | 0.01 | 0.12 | 0.01 | 0.10 | 0.00 | 0.03 | 0.01 | 0.95 | 0.05 |
| 5000 | 1.0 | 0.00 | 0.94 | 0.02 | 0.00 | 0.95 | 0.02 | 0.01 | 0.25 | 0.01 | 0.06 | 0.00 | 0.02 | 0.00 | 0.95 | 0.03 |
| 10000 | 0.2 | 0.01 | 0.97 | 0.08 | 0.01 | 0.97 | 0.07 | 0.03 | 0.04 | 0.00 | 0.30 | 0.00 | 0.06 | 0.02 | 0.95 | 0.10 |
| 10000 | 0.6 | 0.00 | 0.96 | 0.03 | 0.00 | 0.97 | 0.02 | 0.01 | 0.06 | 0.00 | 0.10 | 0.00 | 0.02 | 0.01 | 0.94 | 0.03 |
| 10000 | 1.0 | 0.00 | 0.94 | 0.02 | 0.00 | 0.95 | 0.01 | 0.01 | 0.16 | 0.00 | 0.06 | 0.00 | 0.01 | 0.00 | 0.95 | 0.02 |

Table 1: Comparison of methods when the 50% rule holds. Adaptive Lasso stands for the method of Windmeijer et al. [2016] using the median estimator and tuning with cross validation. For each setting and method, say TSHT under $n = 500$ and $C_\gamma = 0.2$, the corresponding row of numbers $(0.09, 0.72, 0.32)$ represent the median absolute error, the empirical coverage and average length of the confidence interval.

The median method does well with respect to median absolute error, but not as well as TSHT or Windmeijer et al. [2016]'s method, and is not near oracle-level performance. Also, we notice that the confidence interval based on bootstrapping the median estimator is not a wise strategy in the one-sample setting. However, this is expected since the median estimator and the R package that implements it assume a two-sample setting with independent samples. Finally, the naive TSLS consistently has the worst performance across all simulation settings. For example, the naive TSLS performs worse than TSHT even at $n = 500$.

Next, we present the setting where the 50% rule is violated, i.e. with $p_z = 7$ IVs where only three satisfy assumptions (A1)-(A3). The other parameters of the simulation remain the same, i.e. $n = (500, 1000, 2000, 5000, 10000)$, $C_\gamma = (0.2, 0.6, 1)$, and $C_\pi = 0.2$. We drop the median method from our comparison because of its poor performance in Table 1.

As expected, in Table 2, the adaptive Lasso approach of Windmeijer et al. [2016] performs as worse as the naive TSLS since the adaptive Lasso depends on (IN1-M) for consistency. In contrast, TSHT, which relies on a more general identifying condition (IN1-P), has low error along with much better coverage than the adaptive Lasso. Also, TSHT requires more samples to achieve the desired level of coverage when the data is generated under the identifying condition (IN1-P) than condition (IN1-M).

Overall, the simulation study shows that TSHT with voting performs no worse than the competing approaches in the literature. When condition (IN1-M) holds, TSHT performs as well as the method proposed by Windmeijer et al. [2016]. But, when condition (IN1-M) fails to hold, but condition (IN1-P) holds, TSHT performs much better with respect to absolute error, coverage and length of confidence

| $n$ | $C_\gamma$ | TSHT | | | Adaptive Lasso | | | Naive TSLS | | | Oracle TSLS | | |
|---|---|---|---|---|---|---|---|---|---|---|---|---|---|
| 500 | 0.2 | 0.37 | 0.17 | 0.38 | 0.35 | 0.18 | 0.39 | 0.41 | 0.00 | 0.33 | 0.09 | 0.97 | 0.51 |
| 500 | 0.6 | 0.11 | 0.24 | 0.13 | 0.13 | 0.17 | 0.13 | 0.14 | 0.00 | 0.11 | 0.03 | 0.93 | 0.17 |
| 500 | 1.0 | 0.07 | 0.21 | 0.08 | 0.08 | 0.18 | 0.08 | 0.09 | 0.00 | 0.07 | 0.02 | 0.94 | 0.10 |
| 1000 | 0.2 | 0.37 | 0.17 | 0.36 | 0.37 | 0.10 | 0.33 | 0.42 | 0.00 | 0.24 | 0.07 | 0.96 | 0.36 |
| 1000 | 0.6 | 0.09 | 0.32 | 0.13 | 0.12 | 0.12 | 0.11 | 0.14 | 0.00 | 0.08 | 0.02 | 0.96 | 0.12 |
| 1000 | 1.0 | 0.06 | 0.24 | 0.07 | 0.07 | 0.11 | 0.07 | 0.09 | 0.00 | 0.05 | 0.01 | 0.93 | 0.07 |
| 2000 | 0.2 | 0.19 | 0.45 | 0.32 | 0.44 | 0.01 | 0.27 | 0.42 | 0.00 | 0.17 | 0.05 | 0.95 | 0.25 |
| 2000 | 0.6 | 0.04 | 0.62 | 0.10 | 0.15 | 0.02 | 0.09 | 0.14 | 0.00 | 0.06 | 0.01 | 0.94 | 0.08 |
| 2000 | 1.0 | 0.03 | 0.55 | 0.06 | 0.09 | 0.02 | 0.05 | 0.09 | 0.00 | 0.03 | 0.01 | 0.94 | 0.05 |
| 5000 | 0.2 | 0.04 | 0.90 | 0.19 | 0.49 | 0.00 | 0.19 | 0.42 | 0.00 | 0.11 | 0.03 | 0.95 | 0.16 |
| 5000 | 0.6 | 0.01 | 0.91 | 0.06 | 0.17 | 0.00 | 0.06 | 0.14 | 0.00 | 0.03 | 0.01 | 0.94 | 0.05 |
| 5000 | 1.0 | 0.01 | 0.91 | 0.04 | 0.10 | 0.00 | 0.04 | 0.09 | 0.00 | 0.02 | 0.01 | 0.94 | 0.03 |
| 10000 | 0.2 | 0.02 | 0.92 | 0.13 | 0.50 | 0.00 | 0.14 | 0.43 | 0.00 | 0.07 | 0.02 | 0.96 | 0.11 |
| 10000 | 0.6 | 0.01 | 0.92 | 0.04 | 0.17 | 0.00 | 0.04 | 0.14 | 0.00 | 0.02 | 0.01 | 0.95 | 0.04 |
| 10000 | 1.0 | 0.00 | 0.94 | 0.03 | 0.10 | 0.00 | 0.03 | 0.09 | 0.00 | 0.01 | 0.00 | 0.94 | 0.02 |

Table 2: Comparison of methods when the 50% rule is violated but plurality rule holds. Adaptive Lasso stands for the method of Windmeijer et al. [2016] using the median estimator and tuning with cross validation. For each setting and method, say TSHT under $n = 500$ and $C_\gamma = 0.2$, the corresponding row of numbers $(0.37, 0.17, 0.38)$ represent the median absolute error, the empirical coverage and average length of the confidence interval.

intervals.

### 5.3. High Dimensional Setting

In this section, we present simulations in high dimensions. We use the same data generating models as before, except we have $p_z = 100$ instruments with the first $s_{z1} = 7$ being relevant and the first $s_{z2} = 5$ being valid. We also have $p_x = 150$ covariates with $s_{x2} = s_{x1} = 10$. We refer to this case as the high dimensional instruments and covariates setting. We also consider $p_z = 9$ and $p_x = 150$, which we refer to as the low dimensional instruments and high dimensional covariates setting. The only difference between these two settings is the dimension of instrumental variables. However, from a theoretical standpoint, both settings are considered high dimensional.

Both the instruments and covariates $\mathbf{W}_{i\cdot}$ are generated from a multivariate normal with mean zero and covariance $\Sigma_{ij}^* = 0.5^{|i-j|}$ for $1 \leq i, j \leq p_x + p_z$. The other parameters for the models are: $\beta^* = 1$, $\phi^* = (0.6, 0.7, 0.8, \cdots, 1.5, 0, 0, \cdots, 0) \in \mathbb{R}^{150}$, $\psi^* = (1.1, 1.2, 1.3, \cdots, 2.0, 0, 0, \cdots, 0) \in \mathbb{R}^{150}$, $\mathrm{Var}(\epsilon_{i1}) = \mathrm{Var}(\epsilon_{i2}) = 1.5$, and $\mathrm{Cov}(\epsilon_{i1}, \epsilon_{i2}) = 0.75$.

We vary (i) the sample size $n$, (ii) the strength of IV via $\boldsymbol{\gamma}^*$, and the (iii) the degree of violations of (A2) and (A3) via $\boldsymbol{\pi}^*$. For the sample size, we let $n = (200, 300, 1000, 2500)$. For the IV strength, we set $\boldsymbol{\gamma}^* = C_\gamma \cdot (1, 1, 1, 1, 1, 1, 1, 0, 0, \ldots, 0)$ with $C_\gamma = 0.5$. For violations of (A2) and (A3), we set $\boldsymbol{\pi}^* = (0, 0, 0, 0, 0, 1, 1, 0, 0, \ldots, 0) \cdot C_\pi$ where $C_\pi$ is a constant that we vary to change the magnitude of $\boldsymbol{\pi}^*$.

We compare TSHT to the oracle TSLS where the oracle uses only the relevant and valid instruments, i.e. knows the seven relevant instruments, of which the first five are valid. We do not include the naive TSLS because it is not feasible in high dimensions. We also do not include other methods because they were not designed with high dimensionality in mind. The high dimensional instruments and covariate setting is present in Table 3 while the low dimensional instruments and high dimensional covariates setting is present in Table 4. To mimic the low dimensional results, Table 3 presents the result for $C_\pi = (0.25, 0.5, 1)$. In both settings, for $n = 200$, our method TSHT does not achieve the desired level of coverage, although coverage improves dramatically once the violation of (A2) and (A3) becomes bigger, i.e. when $C_\pi = 1$. When $n \geq 300$ and if the violations of (A2) and (A3) are substantial, TSHT achieves the desired level of coverage with absolute error and length of the confidence interval that is comparable to the oracle.

| $n$ | $C_\pi$ | | TSHT | | | Oracle | |
|---|---|---|---|---|---|---|---|
| 200 | 0.25 | 0.162 | 0.162 | 0.202 | 0.038 | 0.956 | 0.219 |
| 200 | 0.50 | 0.129 | 0.448 | 0.232 | 0.036 | 0.962 | 0.218 |
| 200 | 1.00 | 0.056 | 0.876 | 0.259 | 0.036 | 0.956 | 0.221 |
| 300 | 0.25 | 0.155 | 0.080 | 0.164 | 0.033 | 0.952 | 0.179 |
| 300 | 0.50 | 0.093 | 0.516 | 0.197 | 0.029 | 0.952 | 0.177 |
| 300 | 1.00 | 0.041 | 0.906 | 0.209 | 0.029 | 0.946 | 0.176 |
| 1000 | 0.25 | 0.136 | 0.062 | 0.094 | 0.016 | 0.936 | 0.096 |
| 1000 | 0.50 | 0.020 | 0.942 | 0.119 | 0.016 | 0.936 | 0.095 |
| 1000 | 1.00 | 0.020 | 0.958 | 0.120 | 0.016 | 0.964 | 0.096 |
| 2500 | 0.25 | 0.015 | 0.802 | 0.068 | 0.011 | 0.946 | 0.060 |
| 2500 | 0.50 | 0.012 | 0.956 | 0.069 | 0.011 | 0.948 | 0.060 |
| 2500 | 1.00 | 0.011 | 0.954 | 0.069 | 0.010 | 0.942 | 0.060 |

Table 3: Performance of TSHT in high dimensional instruments and covariates with $p_x = 150$ and $p_z = 100$. For each setting and method, say TSHT under $n = 200$ and $C_\pi = 0.25$, the row of numbers $(0.162, 0.162, 0.202)$ represents the median absolute error, the empirical coverage and average length of the confidence interval.

## 6. Application: Causal Effect of Years of Education on Annual Earnings

To demonstrate our method in a real setting, we analyze the causal effect of years of education on yearly earnings, which has been studied extensively in economics using IV methods [Angrist and Krueger, 1991, Card, 1993, 1999]. The data comes from the Wisconsin Longitudinal Study (WLS), a longitudinal study that has kept track of American high school graduates from Wisconsin since 1957, and we examine the relationship between graduates' earnings and education from the 1974 survey Hauser [2005], roughly 20 years after they graduated from high school. Our analysis includes $N = 3772$ individuals, 1784 males and 1988 females. For our outcome, we use imputed log total yearly earnings prepared by the WLS (see WLS documentation and Hauser [2005] for details) and for the treatment, we use the total years of education, all from the 1974 survey. The median total earnings is $\$9,200$ with a 25% quartile of $\$1,000$ and a 75% quartile of $\$15,320$ in 1974 dollars. The mean years of total education is 13.7 years with a standard deviation of 2.3 years.

| $n$ | $C_\pi$ | TSHT | | | Oracle | | |
|---|---|---|---|---|---|---|---|
| 200 | 0.25 | 0.169 | 0.196 | 0.214 | 0.037 | 0.928 | 0.221 |
| 200 | 0.5 | 0.167 | 0.362 | 0.240 | 0.039 | 0.926 | 0.221 |
| 200 | 1 | 0.057 | 0.852 | 0.276 | 0.041 | 0.942 | 0.222 |
| 300 | 0.25 | 0.155 | 0.094 | 0.170 | 0.031 | 0.938 | 0.178 |
| 300 | 0.5 | 0.123 | 0.426 | 0.198 | 0.031 | 0.956 | 0.177 |
| 300 | 1 | 0.043 | 0.916 | 0.222 | 0.030 | 0.960 | 0.177 |
| 1000 | 0.25 | 0.133 | 0.076 | 0.090 | 0.015 | 0.944 | 0.095 |
| 1000 | 0.5 | 0.019 | 0.962 | 0.113 | 0.016 | 0.954 | 0.096 |
| 1000 | 1 | 0.020 | 0.958 | 0.113 | 0.016 | 0.950 | 0.095 |
| 2500 | 0.25 | 0.012 | 0.860 | 0.067 | 0.009 | 0.948 | 0.060 |
| 2500 | 0.5 | 0.012 | 0.952 | 0.068 | 0.010 | 0.950 | 0.060 |
| 2500 | 1 | 0.012 | 0.958 | 0.068 | 0.011 | 0.944 | 0.060 |

Table 4: Performance of TSHT in low dimension instruments and high dimension covariates with $p_x = 150$ and $p_z = 9$ . For each setting and method, say TSHT under $n = 200$ and $C_\pi = 0.25$, the row of numbers $(0.169, 0.196, 0.214)$ represents the median absolute error, the empirical coverage and average length of the confidence interval.

We incorporate many covariates, including sex, graduate's hometown population, educational attainment of graduates' parents, graduates' family income, relative income in graduates' hometown, graduates' high school denomination, high school class size, all measured in 1957 when the participants were high school seniors. We also include 81 genetic covariates, specifically single nucleotide polymorphisms (SNPs), that were part of the WLS to further control for potential variations between graduates; see Section B in Supplementary Materials for details on the non-genetic and genetic covariates. In summary, our data analysis includes 7 non-genetic covariates and 81 genetic covariates.

We used five instruments in our analysis, all derived from past studies of education on earnings [Card, 1993, Blundell et al., 2005, Gary-Bobo et al., 2006]. They are (i) total number of sisters, (ii) total number of brothers, (iii) individual's birth order in the family, all from Gary-Bobo et al. [2006], (iv) proximity to college from Card [1993], and (v) teacher's interest in individual's college education from Blundell et al. [2005], all measured in 1957. Although all these IVs have been suggested to be valid with varying explanations as to why they satisfy (A2) and (A3) after controlling for the aforementioned covariates, in practice, we are always uncertain due to the lack of complete socioeconomic knowledge about the effect of these IVs. Our method should provide some protection against this uncertainty compared to traditional methods where they simply assume that all five IVs are valid. Also, the first-stage F-test produces an F-statistic of 90.3 with a p-value less than $10^{-16}$, which indicates a very strong set of instruments. For more details on the instruments, see Section B of the Supplementary Materials.

Table 5 summarizes the results of our data analysis. OLS refers to running a regression of the treatment and the covariates on the outcome and looking at the slope coefficient of the treatment variable. TSLS refers to running two-stage least squares as described in Section 5 under the operating assumption that all the five instruments are valid; this is the usual and most popular analysis in the IV literature. Finally, we run our procedure TSHT with voting.

| Method | Point Estimate | 95% Confidence Interval |
|--------|----------------|-------------------------|
| OLS    | 0.097          | (0.051, 0.143)          |
| TSLS   | 0.169          | (0.029, 0.301)          |
| TSHT   | 0.062          | (0.046, 0.077)          |

Table 5: Estimates of the Effect of Years of Education on Log Earnings. OLS is ordinary least squares, TSLS is two-stage least squares, and TSHT is our method.

The OLS estimate suggests a positive association between education and earnings, with a statistically significant result at $\alpha = 0.05$ level. This agrees with previous literature which suggests a statistically significant positive association between years of education and log earnings Card [1999]. However, OLS does not completely control for confounding even after controlling for covariates. TSLS provides an alternative method of controlling for confounding by using instruments so long as all the instruments satisfy the three core assumptions and the inclusion of covariates helps make these assumptions more plausible. Unfortunately, we notice that the TSLS estimate in Table 5 is inconsistent with previous studies' estimates among individuals from the U.S. between the 1950s to the 1970s, which range from 0.06 to 0.13 (see Table 4 in Card [1999]). Our method, which addresses the concern for invalid instruments with TSLS, provides an estimate of 0.062, which is more consistent with previous studies' estimates of the effect of years of education on earnings.

The data analysis suggests that our method can be a useful tool in IV analysis when there is concern for invalid instruments, even after attempting to mitigate this problem via covariates. Our method provides much more accurate estimates of the returns on education than TSLS, which naively assumes all the instruments are valid.

## 7. Conclusion and Discussion

We present a method to estimate the effect of the treatment on the outcome using instrumental variables where we do not make the assumption that all the instruments are valid. Our approach is based on the novel TSHT procedure with majority and plurality voting, which is shown to succeed in selecting valid IVs in the presence of possibly invalid IVs even when the 50% rule is violated and produces robust confidence intervals. In simulation and in real data settings, our approach provides a more robust analysis than the traditional IV approaches or recent methods in the invalid IV literature, most notably TSLS, by providing some protection against possibly invalid instruments. In fact, our numerical studies suggest that one should always use TSHT and TSHT will generally reach oracle performance with $n \geq 2000$. Overall, based on the theoretical analysis and the numerical studies, we believe TSHT with voting can be a valuable tool for researchers in Mendelian Randomization and instrumental variables to treatment effects whenever there are concerns for invalid IVs, which is often the case in practice.

Finally, as discussed in Section 4.2, our theoretical analysis for the case of invalid IVs in high dimensions require Assumptions (IN2) and (IN3). We believe assumption (IN3) is most likely necessary for the invalid IV problem in high dimensions. We believe this, in part, because of the model selection literature by Leeb and Pötscher [2005] who pointed out that "in general no model selector can be uniformly consistent for the most parsimonious true model" and hence the post-model-selection inference is generally non-uniform. Consequently, the set of competing models has to be "well separated" such that we can consistently select a correct model. Assumption (IN3) serves as this "well separated" condition in our invalid IV problem. While some recent work in high dimensional inference [Zhang and Zhang, 2014, Javanmard and Montanari, 2014, van de Geer et al., 2014, Chernozhukov et al., 2015, Cai and Guo, 2017] does not make this "well separated" assumption, our invalid IV problem is of a different nature than the prior work because a single invalid IV declared as valid can ruin inference while said prior works assume covariates are exogenous and moments are known perfectly. It is certainly possible that advanced methods can weaken (IN3) and we leave this as a direction for further research.

## Acknowledgments

The research of Hyunseung Kang was supported in part by NSF Grant DMS-1502437. The research of T. Tony Cai was supported in part by NSF Grants DMS-1208982 and DMS-1403708, and NIH Grant R01 CA127334. The research of Dylan S. Small was supported in part by NSF Grant SES-1260782.

## A. Supplementary Materials: Valid IVs After Controlling for High Dimensional Covariates

In this section, we present the method and theory for the case where the instruments are assumed to be valid (i.e. no direct effect and no unmeasured confounding) after conditioning on high dimensional covariates. This setup was considered in Gautier and Tsybakov [2011], Belloni et al. [2012], Fan and Liao [2014], Chernozhukov et al. [2015]. We can simply the procedure present in Section 3.7 by taking $\widehat{\mathcal{V}} = \widehat{\mathcal{S}}$. The resulting estimator for $\beta^*$ is

$$\widehat{\beta}_H = \frac{\sum_{j \in \widehat{\mathcal{S}}} \widehat{\gamma}_j \widehat{\Gamma}_j}{\sum_{j \in \widehat{\mathcal{S}}} \widehat{\gamma}_j^2}. \tag{25}$$

The corresponding confidence interval for $\beta^*$ would be

$$\left( \widehat{\beta}_H - z_{1-\alpha/2} \sqrt{\widehat{\mathrm{Var}}_H/n}, \quad \widehat{\beta}_H + z_{1-\alpha/2} \sqrt{\widehat{\mathrm{Var}}_H/n} \right). \tag{26}$$

Here, $\widehat{\mathrm{Var}}_H$ is $\widehat{\mathrm{Var}}$ in (16) except we replace $\widehat{\mathcal{V}} = \widehat{\mathcal{S}}$ and $\widehat{\beta}$ with $\widehat{\beta}_H$.

REMARK 1. *When there is no baseline covariates, a procedure different from (25) can be proposed to estimate $\beta$. Note that $\beta = \langle \gamma, \Gamma \rangle / \langle \gamma, \Gamma \rangle$. Hence, estimation of $\beta$ can be decomposed into two steps, estimation of inner product $\langle \gamma, \Gamma \rangle$ and estimation of the quadratic norm $\|\gamma\|_2^2$. We can estimate these two quantities by the method proposed in Guo et al. [2016]. See also Section 8 of Cai and Guo [2016b] for the inference of quadratic functional in high-dimensional linear regression.*

In the following, we state the theoretical results for valid IVs after controlling for high dimensional covariates. Under the assumptions (R1)-(R3), Theorem 4 shows if the instruments are valid after conditioning on many covariates, then the estimator $\widehat{\beta}_H$ in our procedure is consistent and asymptotically normal.

THEOREM 4. *Suppose we have valid IVs, that is $\boldsymbol{\pi}^* = 0$ in (2), and the assumptions (R1) − (R3) hold. The following property holds for the estimator $\widehat{\beta}_H$,*

$$\sqrt{n} \left( \widehat{\beta}_H - \beta^* \right) = T^{\beta^*} + \Delta^{\beta^*}, \tag{27}$$

*where $T^{\beta^*} \mid \mathbf{W} \sim N\left(0, \mathrm{Var}_H\right)$, $\mathrm{Var}_H = \sigma^2 \left\| \sum_{j \in \mathcal{S}^*} \boldsymbol{\gamma}_j^* \mathbf{W} \widehat{\mathbf{U}}_{\cdot j} / \sqrt{n} \right\|_2^2 / \|\boldsymbol{\gamma}^*\|_2^4$ and $\Delta^{\beta^*} / \sqrt{\mathrm{Var}_H} \xrightarrow{p} 0$ as $\sqrt{s_{z1}} s \log p / \sqrt{n} \to 0$.*

Theorem 4 states that if the IVs satisfy the exclusion restriction and no unmeasured confounding after conditioning on many covariates, $\widehat{\beta}_H$ defined in (25) is a consistent estimator and the dominating part of the scaled difference $\sqrt{n}(\widehat{\beta}_H - \beta)$ is normal. Based on the asymptotic normality established in (27), the following theorem justifies the coverage property of the confidence interval proposed in (26) under the assumption that all instruments have no direct effect and are unconfounded after conditioning on many covariates.

THEOREM 5. *Suppose we have valid IVs, that is $\boldsymbol{\pi}^* = 0$ in (2) and the assumptions $(R1)-(R3)$ hold. Assuming $\sqrt{s_{z1}}s\log p/\sqrt{n} \to 0$, the confidence interval given in (26) has asymptotically coverage probability $1 - \alpha$, i.e.,*

$$\mathbf{P}\left\{\beta^* \in \left(\widehat{\beta}_H - z_{1-\alpha/2}\sqrt{\widehat{\mathrm{Var}}_H/n}, \quad \widehat{\beta}_H + z_{1-\alpha/2}\sqrt{\widehat{\mathrm{Var}}_H/n}\right)\right\} \to 1 - \alpha, \quad (28)$$

Theorem 5 is similar to a result given in Chernozhukov et al. [2015], who studied IV estimators in high dimensional regime where all the instruments are valid after conditioning. However, there are some notable differences between our results and those in Chernozhukov et al. [2015] in terms of sparsity and instrument-covariate modeling assumptions that are required to achieve $1 - \alpha$ coverage. A simulation study is carried out in Section 5 to compare our procedure to that of the oracle.

## B.  Supplementary Materials: Details on Data Analysis

### B.1.  Genetic Data

This is the list of genetic covariates, i.e. 81 SNPs, used in the data analysis: rs1018381, rs1042838, rs11178997, rs1137100, rs11564774, rs11902591, rs12152850, rs12313279, rs12602084, rs12664989, rs12913832, rs1421085, rs1424954, rs1435252, rs144848, rs1501299, rs1535_chr11, rs17070145, rs174575, rs17529477, rs17561, rs17571, rs1799913, rs1799945, rs1799966, rs1799998, rs1800046, rs1800497, rs1800795, rs1800955, rs1805420, rs1937_chr10, rs2059693, rs2061174, rs2071219, rs2237436, rs2241766, rs2242592, rs2254298, rs2306604, rs25533, rs2760118, rs2779562, rs2963238, rs35118453, rs363050, rs3749386, rs3788862, rs3797297, rs3802657, rs3853248, rs3990403, rs4073366, rs4245147, rs429358, rs4502885, rs45537037, rs4680, rs4986852, rs5031016, rs592389, rs6152_chrX, rs6166, rs6169, rs6265, rs6277, rs6312, rs6314 rs6318 rs669, rs707555, rs7412, rs760761, rs7761133, rs7997012, rs8039957, rs8076005, rs8191992, rs821616, rs878567, rs908867.

Similar to another analysis with the WLS genetic data, we remove individuals with more than 10% missing genotype data Roetker et al. [2012]. We also remove SNPs with no variation in our data set or SNPs with more than 20% missing. The other missing values were imputed using the most frequent genotype for that SNP. Following a Mendelian randomization study by Timpson et al. [2005], we code the SNPs using an additive model.

### B.2.  Covariates in the Data Analysis

We include the following nine covariates for our dat analysis: sex, graduate's hometown population, educational attainment of graduates' parents, graduates' family income, relative income in graduates' hometown, graduates' high school denomination, high school class size, all measured in 1957 when the participants were high school seniors. The details of these covariates are in Table 6

| Covariate | Statistics |
|---|---|
| Sex | 1784 males (47.3%), 1988 females (52.7%) |
| Hometown Population | |
|     Under 1,000 | 585 (15.5%) |
|     1,000 - 2,499 | 584 (15.5%) |
|     2,500 - 9,999 | 771 (20.4 %) |
|     10,000 - 24,999 | 375 (9.9%) |
|     25,000 - 49,999 | 602 (16.0 %) |
|     50,000 - 99.999 | 226 (7.1%) |
|     100,000 - 150,000 | 112 (3.0%) |
|     Over 150,000 | 477 (12.6%) |
| Father's Total Years of Education | 10.4 years (SD: 3.2) |
| Mother's Total Years of Education | 10.8 years (SD: 3.0) |
| Family Income Relative to Community | |
|     Considerably below average | 15 (0.4%) |
|     Somewhat below average | 246 (6.5 %) |
|     Average | 2670 (70.8%) |
|     Somewhat above average | 783 (20.8%) |
|     Considerably above average | 58 (1.5%) |
| School Type | 479 catholic/private (12.7%) , 3293 public (87.3%) |
| Graduating Class Size | 166.8 students (SD: 131.5) |

Table 6: Covariates used in the data analysis. All the covariates were measured in 1957. SD stands for standard deviation.

### B.3.  Instruments in the Data Analysis

They are five instruments in the analysis: (i) total number of sisters, (ii) total number of brothers, (iii) individual's birth order in the family, (iv) proximity to college, and (v) teacher's interest in individuals' college education, all measured in 1957.

A small proportion of individuals had missing teacher's interest in individual's college education. To accommodate this, we created a binary covariate to indicate missingness and imputed missing observations as the average value for that covariate. In total, there were 119 individuals with missing values for teacher's interest in education out of 3772 individuals.

## C.  Supplementary Materials: Proofs of Theorems

In this section, we provide detailed proofs of Lemma 1, Theorems 1, 2, 4, 5 and 3. Proof of extra lemmas are presented in next section. Before presenting the proof, we will introduce the notations used throughout the proof.

### C.1.  Notations

For any vector $\mathbf{v} \in \mathbb{R}^p$, let $\mathbf{v}_j$ denote the $j$th element of $\mathbf{v}$. Let $\|\mathbf{v}\|_1$, $\|\mathbf{v}\|_2$, and $\|\mathbf{v}\|_\infty$ be the usual $1, 2$ and $\infty$-norms, respectively. Let $\|\mathbf{v}\|_0$ denote the 0-norm, i.e. the number of non-zero elements in $\mathbf{v}$. The support of $\mathbf{v}$, denoted as $\operatorname{supp}(\mathbf{v}) \subseteq \{1, \ldots, p\}$, is defined as the set containing the non-zero elements of the vector $\mathbf{v}$, i.e.

| Instrument | Statistics |
| --- | --- |
| Total number of sisters | Median: 1 (25Q: 1, 75Q: 2) |
| Total number of brothers | Median: 1 (25Q: 1, 75Q: 2) |
| Birth order | Median: 1 (25Q: 1, 75Q: 3) |
| Distance from High School (HS) to College | |
|   HS more than 15 miles from any college | 1,487 (39.4%) |
|   HS 15 miles or less from extension center | 168 (4.4 %) |
|   HS 15 miles or less from state college | 196 (5.2%) |
|   HS less than 15 miles from private college | 159 (4.2 %) |
|   HS in city with extension center | 446 (1.2 %) |
|   HS in city with stage college | 283 (7.5 %) |
|   HS in city with private college | 205 (5.4 %) |
|   HS 15 miles or less from state university | 239 (6.3 %) |
|   HS in city with state university (Milwaukee or Madison) | 589 (15.6%) |
| Teacher's Interest in Individual's College Education | |
|   Discouraged to attend college | 47 (1.2%) |
|   Had no effect or no response | 1825 (48.4%) |
|   Encouraged to attend college | 1900 (50.4%) |

Table 7: Instruments used in the data analysis. All the instruments were measured in 1957. 25Q and 75Q stand for 25% and 75% quartile, respectively. There were 119 individuals who had no response to teacher's interest in individual's college education.

$j \in \text{supp}(\mathbf{v})$ if and only if $\mathbf{v}_j \neq 0$. Also, for a vector $\mathbf{v} \in \mathbb{R}^p$ and set $J \subseteq \{1, \ldots, p\}$, we denote $\mathbf{v}_J \in \mathbb{R}^p$ to be the vector where all the elements except whose indices are in $J$ are zero. For a set $J$, $|J|$ denotes its cardinality.

For any $n$ by $p$ matrix $\mathbf{M} \in \mathbb{R}^{n \times p}$, we denote the $(i, j)$ element of matrix $\mathbf{M}$ as $\mathbf{M}_{ij}$, the $i$th row as $\mathbf{M}_{i.}$, and the $j$th column as $\mathbf{M}_{.j}$. Let $\mathbf{M}^{\intercal}$ be the transpose of $\mathbf{M}$. Finally, $\|\mathbf{M}\|_{\infty}$ represents the element-wise matrix sup norm of matrix $\mathbf{M}$.

For a sequence of random variables $X_n$, we use $X_n \xrightarrow{p} X$ and $X_n \xrightarrow{d} X$ to represent that $X_n$ converges to $X$ in probability and in distribution, respectively. For any two sequences $a_n$ and $b_n$, we will write $a_n \gg b_n$ if $\limsup \frac{b_n}{a_n} = 0$ and write $a_n \ll b_n$ if $b_n \gg a_n$. We use $c$ and $C$ to denote generic positive constants that may vary from place to place.

Throughout the whole proof section, we will use $\beta, \boldsymbol{\gamma}, \boldsymbol{\Gamma}, \boldsymbol{\psi}, \boldsymbol{\Psi}, \boldsymbol{\pi}, \boldsymbol{\Theta}, \boldsymbol{\Sigma}, T^{\beta}, \Delta^{\beta}$ to stand for $\beta^*, \boldsymbol{\gamma}^*, \boldsymbol{\Gamma}^*, \boldsymbol{\psi}^*, \boldsymbol{\Psi}^*, \boldsymbol{\pi}^*, \boldsymbol{\Theta}^*, \boldsymbol{\Sigma}^*, T^{\beta^*}, \Delta^{\beta^*}$ respectively and define

$$\widehat{\mathbf{v}}^{[j]} = \mathbf{W}^{\intercal} \widehat{\mathbf{U}}_{.j} \quad \text{for} \quad 1 \leq j \leq p_z.$$

We also introduce the notation $\boldsymbol{\Omega} = \boldsymbol{\Sigma}^{-1}$, $\sigma_1 = \sqrt{\boldsymbol{\Theta}_{11}}$, $\sigma_2 = \sqrt{\boldsymbol{\Theta}_{22}}$ and $\boldsymbol{\Pi}_{i.} = (\epsilon_{i1}, \epsilon_{i2})$. Let $M_2 = \max\{1/\lambda_{\min}(\boldsymbol{\Theta}), \lambda_{\max}(\boldsymbol{\Theta})\}$ and hence $1/M_2 \leq \lambda_{\min}(\boldsymbol{\Theta}) \leq \lambda_{\max}(\boldsymbol{\Theta}) \leq M_2$. We normalize the columns of $\mathbf{W}$ as $\mathbf{H}_{.j} = \sqrt{n}\mathbf{W}_{.j}/\|\mathbf{W}_{.j}\|_2$ for $j \in [p]$. Let $\text{Diag} = \text{diag}(\|\mathbf{W}_{.j}\|_2/\sqrt{n})_{1 \leq j \leq p}$ denote the $p \times p$ diagonal matrix with $(j, j)$ entry to be $\|\mathbf{W}_{.j}\|_2/\sqrt{n}$. We set $\lambda_0 = \sqrt{2.05 \log p/n} = (1 + \gamma_0)\sqrt{2\delta_0 \log p/n}$, where $\delta_0 = \sqrt{1.025} > 1$ and $\gamma_0 = (1.025)^{\frac{1}{4}} - 1 > 0$. Take $\epsilon_0 = 2.01/\gamma_0 + 1$, $\nu_0 = 0.01$, $\tau_0 = 0.01$, $C_1 = 2.25$, $c_0 = 1/6$ and $C_0 = 3$. We also assume that $\log p/n \to 0$ and $\delta_0 \log p > 2$. Rather than use the constants directly in the following discussion, we use $\delta_0, \pi_0, \epsilon_0, \nu_0, C_1, C_0$ and $c_0$ to represent the above fixed constants in the following discussion. We review the following definition of restricted eigenvalue introduced in Bickel et al. [2009],

$$\kappa(X, k, \alpha_0) = \min_{\substack{J_0 \subset \{1, \cdots, p\}, \\ |J_0| \le k}} \min_{\substack{\delta \ne 0, \\ \|\delta_{J_0^c}\|_1 \le \alpha_0 \|\delta_{J_0}\|_1}} \frac{\|X\delta\|_2}{\sqrt{n}\|\delta_{J_0}\|_2}. \tag{29}$$

Define the oracle estimator of $\sigma_1$ and $\sigma_2$ as

$$\sigma_1^{ora} = \frac{1}{\sqrt{n}}\|Y - \mathbf{Z\Gamma} - \mathbf{X\Psi}\|_2 \text{ and } \sigma_2^{ora} = \frac{1}{\sqrt{n}}\|D - \mathbf{Z\gamma} - \mathbf{X\psi}\|_2,$$

and

$$\tau = \sqrt{1 + \epsilon_0} \frac{2\sqrt{s}\lambda_0}{\kappa(\mathbf{H}, 4s, 1 + 2\epsilon_0)}. \tag{30}$$

### C.2. Proof of Lemma 1

In the proof, we simply the notation by using $\beta, \boldsymbol{\gamma}, \boldsymbol{\Gamma}, \boldsymbol{\psi}, \boldsymbol{\Psi}, \boldsymbol{\pi}, \boldsymbol{\Theta}, \boldsymbol{\Sigma}$ to represent $\beta^*, \boldsymbol{\gamma}^*, \boldsymbol{\Gamma}^*, \boldsymbol{\psi}^*, \boldsymbol{\Psi}^*, \boldsymbol{\pi}^*, \boldsymbol{\Theta}^*, \boldsymbol{\Sigma}^*$ respectively. Define $\boldsymbol{\Pi} = (\epsilon_{\cdot 1}, \epsilon_{\cdot 2})$ and $\mathcal{I} = \{1, 2, \cdots, p_z\}$. We first note the following expression for $\widehat{\boldsymbol{\gamma}}_j$ and $\widehat{\boldsymbol{\Gamma}}_j$ for $j \in \mathcal{I}$,

$$\sqrt{n}\left(\widehat{\boldsymbol{\gamma}}_j - \boldsymbol{\gamma}_j\right) = \left(\widehat{\boldsymbol{\Sigma}}^{-1}\right)_{j,\cdot} \frac{1}{\sqrt{n}}\mathbf{W}^{\mathsf{T}}\boldsymbol{\Pi}_{\cdot 2} \quad \text{and} \quad \sqrt{n}\left(\widehat{\boldsymbol{\Gamma}}_j - \boldsymbol{\Gamma}_j\right) = \left(\widehat{\boldsymbol{\Sigma}}^{-1}\right)_{j,\cdot} \frac{1}{\sqrt{n}}\mathbf{W}^{\mathsf{T}}\boldsymbol{\Pi}_{\cdot 1} \tag{31}$$

and the following limiting theorem ( Theorem 3.1 in Wooldridge [2010]),

$$\widehat{\boldsymbol{\gamma}} \xrightarrow{p} \boldsymbol{\gamma} \quad \text{and} \quad \widehat{\boldsymbol{\Gamma}} \xrightarrow{p} \boldsymbol{\Gamma}, \tag{32}$$

$$\sqrt{n}\left(\widehat{\boldsymbol{\gamma}} - \boldsymbol{\gamma}\right) \xrightarrow{d} N\left(0, \Theta_{22}\left(\boldsymbol{\Sigma}^{-1}\right)_{\mathcal{I},\mathcal{I}}\right) \quad \text{and} \quad \sqrt{n}\left(\widehat{\boldsymbol{\Gamma}} - \boldsymbol{\Gamma}\right) \xrightarrow{d} N\left(0, \Theta_{11}\left(\boldsymbol{\Sigma}^{-1}\right)_{\mathcal{I},\mathcal{I}}\right). \tag{33}$$

Note that

$$\frac{\sqrt{\widehat{\boldsymbol{\Theta}}_{22}}\|\mathbf{W}(\widehat{\boldsymbol{\Sigma}}^{-1})_{\cdot j}\|_2}{\sqrt{n}} \xrightarrow{p} \sqrt{\Theta_{22}\left(\boldsymbol{\Sigma}^{-1}\right)_{jj}}. \tag{34}$$

We define the following events

$$\begin{aligned}
\mathcal{B}_1 &= \left\{\widehat{\mathcal{S}} = \mathcal{S}^*\right\} \\
\mathcal{B}_2 &= \left\{\min_{k \in \mathcal{V}^*} \|\mathrm{VM}_k\|_0 > \frac{q}{2} > \max_{k \in \mathcal{S}^* \backslash \mathcal{V}^*} \|\mathrm{VM}_k\|_0\right\} \\
\mathcal{B}_3 &= \left\{\min_{k \in \mathcal{V}^*} \|\mathrm{VM}_k\|_0 = \max_{k \in \mathcal{S}^*} \|\mathrm{VM}_k\|_0 > \max_{k \in \mathcal{S}^* \backslash \mathcal{V}^*} \|\mathrm{VM}_k\|_0\right\}
\end{aligned} \tag{35}$$

where $q = |\widehat{\mathcal{S}}|$. On the event $\mathcal{B} = \mathcal{B}_1 \cap (\mathcal{B}_2 \cup \mathcal{B}_3)$, we have $\widehat{\mathcal{V}} = \mathcal{V}^*$ and it is sufficient to show that $\lim_{n\to\infty} \mathbf{P}(\mathcal{B}) = 1$. For $j \in \mathcal{S}^*$, we have

$$|\widehat{\gamma}_j| - \frac{\sqrt{\widehat{\boldsymbol{\Theta}}_{22}}\|\mathbf{W}(\widehat{\boldsymbol{\Sigma}}^{-1})_{\cdot j}\|_2}{\sqrt{n}}\sqrt{\frac{2.01\log n}{n}} \xrightarrow{p} |\gamma_j| > 0, \tag{36}$$

where the convergence follows from (32) and (34). For $j \in (\mathcal{S}^*)^c$, we have

$$\sqrt{\frac{n}{2.01\log n}}|\widehat{\gamma}_j| - \frac{\sqrt{\widehat{\boldsymbol{\Theta}}_{22}}\|\mathbf{W}(\widehat{\boldsymbol{\Sigma}}^{-1})_{\cdot j}\|_2}{\sqrt{n}} \xrightarrow{p} -\sqrt{\boldsymbol{\Theta}_{22}\left(\boldsymbol{\Sigma}^{-1}\right)_{jj}} < 0, \tag{37}$$

where the convergence follows from (33) and (34). Combining (36) and (37), we establish that

$$\lim_{n\to\infty} \mathbf{P}\left(\mathcal{B}_1\right) = 1. \tag{38}$$

Without loss of generality, we assume $\widehat{\mathcal{S}} = \{1, 2, \cdots, q\}$ and analyze $\widehat{\mathcal{V}}^{[1]}$. Note that

$$\widehat{\beta}^{[1]} - \frac{\boldsymbol{\Gamma}_1}{\gamma_1} \xrightarrow{p} 0, \tag{39}$$

and hence we have

$$\sqrt{\widehat{\boldsymbol{\Theta}}_{11} + (\widehat{\beta}^{[1]})^2\widehat{\boldsymbol{\Theta}}_{22} - 2\widehat{\beta}^{[1]}\widehat{\boldsymbol{\Theta}}_{12}}\frac{\|\mathbf{W}((\widehat{\boldsymbol{\Sigma}}^{-1})_{k\cdot} - \frac{\widehat{\gamma}_k}{\widehat{\gamma}_1}(\widehat{\boldsymbol{\Sigma}}^{-1})_{1\cdot})\|_2}{\sqrt{n}}$$

$$\xrightarrow{p} \sqrt{\boldsymbol{\Theta}_{11} + \left(\frac{\boldsymbol{\Gamma}_1}{\gamma_1}\right)^2\boldsymbol{\Theta}_{22} - 2\frac{\boldsymbol{\Gamma}_1}{\gamma_1}\boldsymbol{\Theta}_{12}}\sqrt{(\boldsymbol{\Sigma}^{-1})_{kk} + \left(\frac{\gamma_k}{\gamma_1}\right)^2(\boldsymbol{\Sigma}^{-1})_{11} - 2\frac{\gamma_k}{\gamma_1}(\boldsymbol{\Sigma}^{-1})_{k1}}. \tag{40}$$

We also have the following expression

$$\begin{aligned}
\widehat{\boldsymbol{\pi}}_k^{[1]} - \left(\boldsymbol{\Gamma}_k - \frac{\boldsymbol{\Gamma}_1}{\gamma_1}\gamma_k\right) &= \left(\widehat{\boldsymbol{\Gamma}}_k - \frac{\widehat{\boldsymbol{\Gamma}}_1}{\widehat{\gamma}_1}\widehat{\gamma}_k\right) - \left(\boldsymbol{\Gamma}_k - \frac{\boldsymbol{\Gamma}_1}{\gamma_1}\gamma_k\right) \\
&= \left(\widehat{\boldsymbol{\Gamma}}_k - \boldsymbol{\Gamma}_k\right) - \frac{\boldsymbol{\Gamma}_1}{\gamma_1}\left(\widehat{\gamma}_k - \gamma_k\right) - \frac{\gamma_k}{\gamma_1^2}\left(\gamma_1\left(\widehat{\boldsymbol{\Gamma}}_1 - \boldsymbol{\Gamma}_1\right) - \boldsymbol{\Gamma}_1\left(\widehat{\gamma}_1 - \gamma_1\right)\right) \\
&\quad + \left(\frac{\widehat{\boldsymbol{\Gamma}}_1}{\widehat{\gamma}_1} - \frac{\boldsymbol{\Gamma}_1}{\gamma_1}\right)\left(\frac{\gamma_k}{\gamma_1}\left(\widehat{\gamma}_1 - \gamma_1\right) - \left(\widehat{\gamma}_k - \gamma_k\right)\right)
\end{aligned} \tag{41}$$

Note that

$$\sqrt{n}\left(\left(\widehat{\boldsymbol{\Gamma}}_k - \boldsymbol{\Gamma}_k\right) - \frac{\boldsymbol{\Gamma}_1}{\gamma_1}\left(\widehat{\gamma}_k - \gamma_k\right) - \frac{\gamma_k}{\gamma_1^2}\left(\gamma_1\left(\widehat{\boldsymbol{\Gamma}}_1 - \boldsymbol{\Gamma}_1\right) - \boldsymbol{\Gamma}_1\left(\widehat{\gamma}_1 - \gamma_1\right)\right)\right)$$

$$= \left(\left(\widehat{\boldsymbol{\Sigma}}^{-1}\right)_{\cdot k} - \frac{\gamma_j}{\gamma_1}\left(\widehat{\boldsymbol{\Sigma}}^{-1}\right)_{\cdot 1}\right)\frac{1}{\sqrt{n}}\mathbf{W}^{\mathsf{T}}\left(\boldsymbol{\Pi}_{\cdot 2} - \frac{\boldsymbol{\Gamma}_1}{\gamma_1}\boldsymbol{\Pi}_{\cdot 1}\right)$$

$$\xrightarrow{d} N\left(0, \left(\boldsymbol{\Theta}_{11} + \left(\frac{\boldsymbol{\Gamma}_1}{\gamma_1}\right)^2\boldsymbol{\Theta}_{22} - 2\frac{\boldsymbol{\Gamma}_1}{\gamma_1}\right)\left(\boldsymbol{\Theta}_{12}(\boldsymbol{\Sigma}^{-1})_{kk} + \left(\frac{\gamma_k}{\gamma_1}\right)^2(\boldsymbol{\Sigma}^{-1})_{11} - 2\frac{\gamma_k}{\gamma_1}(\boldsymbol{\Sigma}^{-1})_{k1}\right)\right), \tag{42}$$

where the convergence follows from Theorem 3.1 in Wooldridge [2010]. By (32) and (33), we have

$$\left(\frac{\widehat{\boldsymbol{\Gamma}}_1}{\widehat{\gamma}_1} - \frac{\boldsymbol{\Gamma}_1}{\gamma_1}\right)\left(\frac{\gamma_k}{\gamma_1}\left(\widehat{\gamma}_1 - \gamma_1\right) - \left(\widehat{\gamma}_k - \gamma_k\right)\right) \xrightarrow{p} 0.$$

Combined with (39) and (42), we have

$$\frac{\sqrt{n}}{\sqrt{\widehat{\boldsymbol{\Theta}}_{11} + (\widehat{\beta}^{[1]})^2\widehat{\boldsymbol{\Theta}}_{22} - 2\widehat{\beta}^{[1]}\widehat{\boldsymbol{\Theta}}_{12}}\frac{\|\mathbf{W}((\widehat{\boldsymbol{\Sigma}}^{-1})_{k\cdot} - \frac{\widehat{\gamma}_k}{\widehat{\gamma}_1}(\widehat{\boldsymbol{\Sigma}}^{-1})_{1\cdot})\|_2}{\sqrt{n}}}\left(\widehat{\boldsymbol{\pi}}_k^{[1]} - \left(\boldsymbol{\Gamma}_k - \frac{\boldsymbol{\Gamma}_1}{\gamma_1}\gamma_k\right)\right) \xrightarrow{p} N(0,1) \tag{43}$$

and hence

$$\widehat{\boldsymbol{\pi}}_k^{[1]} \xrightarrow{p} \boldsymbol{\Gamma}_k - \frac{\boldsymbol{\Gamma}_1}{\gamma_1}\gamma_k. \tag{44}$$

To analyze $\widehat{\mathcal{V}}^{[1]}$, we will first establish the following results,

$$\text{If } \frac{\pi_k}{\gamma_k} = \frac{\pi_1}{\gamma_1}, \text{ then } \lim_{n\to\infty} \mathbf{P}\left(k \in \widehat{\mathcal{V}}^{[1]}\right) = 1. \tag{45}$$

$$\text{If } \frac{\pi_k}{\gamma_k} \neq \frac{\pi_1}{\gamma_1}, \text{ then } \lim_{n\to\infty} \mathbf{P}\left(k \notin \widehat{\mathcal{V}}^{[1]}\right) = 1. \tag{46}$$

<u>Proof of (45)</u> In this case, $\boldsymbol{\Gamma}_k - \frac{\boldsymbol{\Gamma}_1}{\gamma_1}\gamma_k = 0$. By (43), we have

$$\frac{\left|\widehat{\boldsymbol{\pi}}_k^{[1]}\right|}{2.05\sqrt{\widehat{\boldsymbol{\Theta}}_{11} + (\widehat{\beta}^{[1]})^2\widehat{\boldsymbol{\Theta}}_{22} - 2\widehat{\beta}^{[1]}\widehat{\boldsymbol{\Theta}}_{12}}\frac{\|\mathbf{W}((\widehat{\boldsymbol{\Sigma}}^{-1})_{k\cdot} - \frac{\widehat{\gamma}_k}{\widehat{\gamma}_1}(\widehat{\boldsymbol{\Sigma}}^{-1})_{1\cdot})\|_2}{\sqrt{n}}\sqrt{\frac{\log n}{n}}} \xrightarrow{p} 0 \tag{47}$$

and

$$\lim_{n\to\infty} \mathbf{P}\left(k \in \widehat{\mathcal{V}}^{[1]}\right) = 1. \tag{48}$$

<u>Proof of (46)</u> In this case, $\boldsymbol{\Gamma}_k - \frac{\boldsymbol{\Gamma}_1}{\gamma_1}\gamma_k = \left(\frac{\pi_k}{\gamma_k} - \frac{\pi_1}{\gamma_1}\right)\gamma_k \neq 0$. Hence, we have

$$\left|\widehat{\boldsymbol{\pi}}_k^{[1]}\right| - 2.05\sqrt{\widehat{\boldsymbol{\Theta}}_{11} + (\widehat{\beta}^{[1]})^2\widehat{\boldsymbol{\Theta}}_{22} - 2\widehat{\beta}^{[1]}\widehat{\boldsymbol{\Theta}}_{12}}\frac{\|\mathbf{W}((\widehat{\boldsymbol{\Sigma}}^{-1})_{k\cdot} - \frac{\widehat{\gamma}_k}{\widehat{\gamma}_1}(\widehat{\boldsymbol{\Sigma}}^{-1})_{1\cdot})\|_2}{\sqrt{n}}\sqrt{\frac{\log\max\{n, p_z\}}{n}}$$

$$\xrightarrow{p} \left|\left(\frac{\pi_k}{\gamma_k} - \frac{\boldsymbol{\pi}_1}{\boldsymbol{\gamma}_1}\right)\gamma_k\right| > 0 \tag{49}$$

and

$$\lim_{n\to\infty} \mathbf{P}\left(k \notin \widehat{\mathcal{V}}^{[1]}\right) = 1. \tag{50}$$

The results (45) and (46) can be generalized with the same proofs,

$$\text{If } \frac{\pi_k}{\gamma_k} = \frac{\pi_j}{\gamma_j}, \text{ then } \lim_{n\to\infty} \mathbf{P}\left(k \in \widehat{\mathcal{V}}^{[j]}\right) = 1. \tag{51}$$

$$\text{If } \frac{\pi_k}{\gamma_k} \neq \frac{\pi_j}{\gamma_j}, \text{ then } \lim_{n\to\infty} \mathbf{P}\left(k \notin \widehat{\mathcal{V}}^{[j]}\right) = 1. \tag{52}$$

In the following, we will apply (51) and (52) and establish $\lim_{n\to\infty} \mathbf{P}\left(\mathcal{B}_1 \cap (\mathcal{B}_2 \cup \mathcal{B}_3)\right) = 1$. On the event $\mathcal{B}_1$, by (51) and (52), the voting matrix have the following expression

$$\lim_{n\to\infty} \mathbf{P}\left(k \in \widehat{\mathcal{V}}^{[j]} \text{ if and only if } \frac{\pi_k}{\gamma_k} = \frac{\pi_j}{\gamma_j}\right) = 1. \tag{53}$$

- If the assumption (IN1-M) holds, in this case, (IN1-P) automatically holds, then by (53), we have

$$\lim_{n\to\infty} \mathbf{P}\left(\mathcal{B}_2\right) = \lim_{n\to\infty} \mathbf{P}\left(\mathcal{B}_3\right) = 1. \tag{54}$$

- If the assumption (IN1-M) does not hold, but (IN1-P) holds, then then by (53), we have

$$\lim_{n\to\infty} \mathbf{P}\left(\mathcal{B}_2\right) = 0, \quad \lim_{n\to\infty} \mathbf{P}\left(\mathcal{B}_3\right) = 1. \tag{55}$$

Combining (38), (54) and (55), we establish the lemma.

### C.3.    Proof of Theorem 1

Without loss of generality, we prove the case where covariates $\mathbf{X}_{i\cdot}$ are absent. The models (2) and (3) imply that

$$E(\mathbf{Z}_{i\cdot}(Y_i - D_i\beta^* - \mathbf{Z}_{i\cdot}^{\mathsf{T}}\boldsymbol{\pi}^*)) = 0, \quad E(\mathbf{Z}_{i\cdot}(D_i - \mathbf{Z}_{i\cdot}^{\mathsf{T}}\boldsymbol{\gamma}^*)) = 0$$

Combining the two expressions, we obtain $E(\mathbf{Z}_{i\cdot}Y_i) = E(\mathbf{Z}_{i\cdot}\mathbf{Z}_{i\cdot}^{\mathsf{T}})(\beta^*\boldsymbol{\gamma}^* + \boldsymbol{\pi}^*)$. Let $\boldsymbol{\Gamma}^* = \beta^*\boldsymbol{\gamma}^* + \boldsymbol{\pi}^*$. Because we can identify both $\boldsymbol{\gamma}^*$ and $\boldsymbol{\Gamma}^*$ via $E(\mathbf{Z}_{i\cdot}\mathbf{Z}_{i\cdot}^{\mathsf{T}})^{-1}E(\mathbf{Z}_{i\cdot}D_i)$ and $E(\mathbf{Z}_{i\cdot}\mathbf{Z}_{i\cdot}^{\mathsf{T}})^{-1}E(\mathbf{Z}_{i\cdot}Y_i)$, respectively, proving identification for $\beta^*$ and $\boldsymbol{\pi}^*$ (or lack thereof) reduces to proving that there is a unique solution to the equation $\boldsymbol{\Gamma}^* = \beta^*\boldsymbol{\gamma}^* + \boldsymbol{\pi}^*$ given $\boldsymbol{\Gamma}^*$ and $\boldsymbol{\gamma}^*$. Note that for any $\boldsymbol{\Gamma}^*$ and $\boldsymbol{\gamma}^*$, a solution always exist, say by letting $\boldsymbol{\pi}^* = \boldsymbol{\Gamma}^*$ and $\beta^* = 0$.

First, we prove the sufficient part of the theorem. Suppose we have any two sets of solutions $\beta_1^*, \boldsymbol{\pi}_1^*$ with $\mathcal{V}_1^* = \{j \in \mathcal{S}^* \mid \pi_{j,1}^* = 0\}$ and $\beta_2^*, \boldsymbol{\pi}_2^*$ with $\mathcal{V}_2^* = \{j \in \mathcal{S}^* \mid \pi_{j,2}^* = 0\}$, both of which satisfy $\boldsymbol{\Gamma}^* = \beta_1^*\boldsymbol{\gamma}^* + \boldsymbol{\pi}_1^*$ and $\boldsymbol{\Gamma}^* = \beta_2^*\boldsymbol{\gamma}^* + \boldsymbol{\pi}_2^*$. We show that these two sets of the solutions must be equal to each other when (IN1-P) holds.

When (IN1-P) holds, the sets $\mathcal{V}_1^*$ and $\mathcal{V}_2^*$ must be non-empty because (IN1-P) places a lower bound on the size of the two sets to be strictly greater than 0. Now, consider the two cases.

(a) If $\mathcal{V}_1^* \cap \mathcal{V}_2^*$ is non-empty, then any element $j \in \mathcal{V}_1^* \cap \mathcal{V}_2^*$ leads to $\Gamma_j^*/\gamma_j^* = \beta_1^*$ and $\Gamma_j^*/\gamma_j^* = \beta_2^*$, which imply $\beta_1^* = \beta_2^*$ and $\boldsymbol{\pi}_1^* = \boldsymbol{\pi}_2^*$.

(b) If $\mathcal{V}_1^* \cap \mathcal{V}_2^*$ is empty, then for any element $j \in \mathcal{V}_1^*$, we have

$$\Gamma_j^*/\gamma_j^* = \beta_1^*, \quad \Gamma_j^*/\gamma_j^* = \beta_2^* + \pi_{j,2}^*/\gamma_j^*, \quad j \in \mathcal{V}_1^*$$

This simplifies to $\pi_{j,2}^*/\gamma_j^* = (\beta_1^* - \beta_2^*)$, implying that $j \in \mathcal{V}_1^*$ also belongs to the set $\left\{j \in \mathcal{S}^* \mid \pi_{j,2}^*/\gamma_j^* = c\right\}$ where $c = \beta_1^* - \beta_2^*$ and $c$ is non-zero, i.e.

$$\mathcal{V}_1 \subseteq \left\{j \in \mathcal{S}^* \mid \frac{\pi_{j,2}^*}{\gamma_j^*} = \beta_1^* - \beta_2^*, \beta_1^* - \beta_2^* \neq 0\right\}$$

Combining the inequality implied the set relation above and the inequality under (IN1-P) , we have

$$|\mathcal{V}_1^*| \leq \left| \left\{ j \in \mathcal{S}^* \mid \frac{\pi_{j,2}^*}{\gamma_j^*} = \beta_1^* - \beta_2^*, \beta_1^* - \beta_2^* \neq 0 \right\} \right| < |\mathcal{V}_2^*|$$

and thus $|\mathcal{V}_1^*| < |\mathcal{V}_2^*|$. Similar argument for $j \in \mathcal{V}_2^*$ would lead to

$$|\mathcal{V}_2^*| \leq \left| \left\{ j \in \mathcal{S}^* \mid \frac{\pi_{j,1}^*}{\gamma_j^*} = \beta_2^* - \beta_1^*, \beta_2^* - \beta_1^* \neq 0 \right\} \right| < |\mathcal{V}_1^*|$$

and thus $|\mathcal{V}_2^*| < |\mathcal{V}_1^*|$. This is clearly a contradiction and $\mathcal{V}_1^* \cap \mathcal{V}_2^*$ must be non-empty.

Second, we prove the necessary part of the theorem. Suppose $\beta_1^*, \boldsymbol{\pi}_1^*$ with $\mathcal{V}_1^* = \{j \in \mathcal{S}^* \mid \pi_{j,1}^* = 0\}$ satisfies $\boldsymbol{\Gamma}^* = \beta_1^* \boldsymbol{\gamma}^* + \boldsymbol{\pi}^*$. We show that we can obtain another $(\beta_2, \boldsymbol{\pi}_2) \neq (\beta_2^*, \boldsymbol{\pi}_2^*)$ that also satisfies $\boldsymbol{\Gamma}^* = \beta_2^* \boldsymbol{\gamma}^* + \boldsymbol{\pi}_2^*$ when (IN1-P) fails to hold.

When (IN1-P) fails, $|\mathcal{V}_1^*| \leq \max_{c \neq 0} |\{j \in \mathcal{S}^* : \pi_{j,1}^*/\gamma_j^* = c\}|$. Let the maximizing $c$ be denoted a $c_m$; if $c_m$ is undefined, say when the set is empty, pick any $c_m \neq 0$. Then $\beta_2^* = \beta_1^* + (c_m - k)$ and $\boldsymbol{\pi}_2^* = \boldsymbol{\pi}_1^* - (c_m - k)\boldsymbol{\gamma}^*$ for $k \neq c_m$ and $k \neq 0$ satisfy $\boldsymbol{\Gamma}^* = \beta_2^* \boldsymbol{\gamma}^* + \boldsymbol{\pi}_2^*$ because

$$\beta_2^* \boldsymbol{\gamma}^* + \boldsymbol{\pi}_2^* = (\beta_1^* + c_m - k)\boldsymbol{\gamma}^* + \boldsymbol{\pi}_1^* - (c_m - k)\boldsymbol{\gamma}^* = \beta_1^* \boldsymbol{\gamma}^* + \boldsymbol{\pi}_1^* = \boldsymbol{\Gamma}^*$$

We now show that the $\beta_2^*, \boldsymbol{\pi}_2^*$ satisfy $|\mathcal{V}_2^*| \leq \max_{c \neq 0} |\{j \in \mathcal{S}^* : \pi_{j,2}^*/\gamma_j^* = c\}|$. First, if $\mathcal{S}^*$ is empty, the condition is trivially satisfied because $|\mathcal{V}_2^*| = |\{j \in \mathcal{S}^* : \pi_{j,2}/\gamma_j^* = c\}| = 0$ for any $c \neq 0$. If $\mathcal{S}^*$ is non-empty, for any $j \in \mathcal{S}^*$, we have the relation

$$\frac{\pi_{2,j}^*}{\gamma_j^*} = \frac{\pi_{1,j}^*}{\gamma_j^*} - (c_m - k)$$

This relation implies that the set $\mathcal{V}_2^* = \{j \in \mathcal{S}^* \mid \pi_{j,2}^*/\gamma_j^* = 0\}$ is equal to the set $\{j \in \mathcal{S}^* \mid \pi_{j,1}^*/\gamma_j^* = c_m - k\}$. Also, the set $\{j \in \mathcal{S}^* \mid \pi_{j,2}^*/\gamma_j^* = k\}$ is equal to the set $\{j \in \mathcal{S}^* \mid \pi_{1,j}^*/\gamma_j^* = c_m\}$. Combined, we have the following inequality.

$$|\mathcal{V}_2^*| = \left| \left\{ j \in \mathcal{S}^* \mid \frac{\pi_{j,1}^*}{\gamma_j^*} = c_m - k \right\} \right| \leq \left| \left\{ j \in \mathcal{S}^* \mid \frac{\pi_{1,j}^*}{\gamma_j^*} = c_m \right\} \right| = \left| \left\{ j \in \mathcal{S}^* \mid \frac{\pi_{j,2}^*}{\gamma_j^*} = k \right\} \right|$$

where the middle inequality is from the definition of $c_m$. Since $k \neq 0$, the second solution also satisfies the constraint $|\mathcal{V}_2^*| \leq \max_{c \neq 0} |\{j \in \mathcal{S}^* : \pi_{j,2}/\gamma_j^* = c\}|$.

## C.4. Proof of Theorem 2

Define

$$\widehat{A}(\mathcal{V}) = \widehat{\boldsymbol{\Sigma}}_{\mathcal{V},\mathcal{V}} - \widehat{\boldsymbol{\Sigma}}_{\mathcal{V},\mathcal{V}^c} \widehat{\boldsymbol{\Sigma}}_{\mathcal{V}^c,\mathcal{V}^c}^{-1} \widehat{\boldsymbol{\Sigma}}_{\mathcal{V}^c,\mathcal{V}} \quad \text{and} \quad A(\mathcal{V}) = \boldsymbol{\Sigma}_{\mathcal{V},\mathcal{V}} - \boldsymbol{\Sigma}_{\mathcal{V},\mathcal{V}^c} \boldsymbol{\Sigma}_{\mathcal{V}^c,\mathcal{V}^c}^{-1} \boldsymbol{\Sigma}_{\mathcal{V}^c,\mathcal{V}}.$$

We introduce the following lemma to facilitate the proof.

Lemma 3. *Under the assumptions of Theorem 2, we have*

$$\sqrt{n}\left(\frac{\widehat{\gamma}_{\mathcal{V}^*}^{\mathsf{T}}\widehat{A}(\mathcal{V}^*)\widehat{\Gamma}_{\mathcal{V}^*}}{\widehat{\gamma}_{\mathcal{V}^*}^{\mathsf{T}}\widehat{A}(\mathcal{V}^*)\widehat{\gamma}_{\mathcal{V}^*}} - \beta\right) \xrightarrow{d} N\left(0, \frac{\boldsymbol{\Theta}_{11} + \beta^2\boldsymbol{\Theta}_{22} - 2\beta\boldsymbol{\Theta}_{12}}{\gamma_{\mathcal{V}^*}^{\mathsf{T}}A(\mathcal{V}^*)\gamma_{\mathcal{V}^*}}\right). \qquad (56)$$

The estimator defined in (10) can be expressed as $\widehat{\beta}_L = \frac{\widehat{\gamma}_{\widehat{\mathcal{V}}}^{\mathsf{T}}A(\widehat{\mathcal{V}})\widehat{\Gamma}_{\widehat{\mathcal{V}}}}{\widehat{\gamma}_{\widehat{\mathcal{V}}}^{\mathsf{T}}A(\widehat{\mathcal{V}})\widehat{\gamma}_{\widehat{\mathcal{V}}}}$, and hence the difference $\sqrt{n}\left(\widehat{\beta}_L - \beta\right)$ can be expressed as

$$\sqrt{n}\left(\widehat{\beta}_L - \beta\right) = \sqrt{n}\left(\frac{\widehat{\gamma}_{\mathcal{V}^*}^{\mathsf{T}}\widehat{A}(\mathcal{V}^*)\widehat{\Gamma}_{\mathcal{V}^*}}{\widehat{\gamma}_{\mathcal{V}^*}^{\mathsf{T}}\widehat{A}(\mathcal{V}^*)\widehat{\gamma}_{\mathcal{V}^*}} - \beta\right)\mathbf{1}_{\widehat{\mathcal{V}}=\mathcal{V}^*} + \sum_{\mathcal{V}\neq\mathcal{V}^*}\sqrt{n}\left(\frac{\widehat{\gamma}_{\mathcal{V}}^{\mathsf{T}}\widehat{A}(\mathcal{V})\widehat{\Gamma}_{\mathcal{V}}}{\widehat{\gamma}_{\mathcal{V}}^{\mathsf{T}}\widehat{A}(\mathcal{V})\widehat{\gamma}_{\mathcal{V}}} - \beta\right)\mathbf{1}_{\widehat{\mathcal{V}}=\mathcal{V}}$$
$$(57)$$

By Lemma 1, we have $\mathbf{1}_{\widehat{\mathcal{V}}=\mathcal{V}^*} \xrightarrow{p} 1$ and $\mathbf{1}_{\widehat{\mathcal{V}}=\mathcal{V}} \xrightarrow{p} 0$ if $\mathcal{V} \neq \mathcal{V}^*$. Combined with Lemma 3 and Slutsky's theorem, we establish

$$\sqrt{n}\left(\frac{\widehat{\gamma}_{\mathcal{V}^*}^{\mathsf{T}}\widehat{A}(\mathcal{V}^*)\widehat{\Gamma}_{\mathcal{V}^*}}{\widehat{\gamma}_{\mathcal{V}^*}^{\mathsf{T}}\widehat{A}(\mathcal{V}^*)\widehat{\gamma}_{\mathcal{V}^*}} - \beta\right)\mathbf{1}_{\widehat{\mathcal{V}}=\mathcal{V}^*} \xrightarrow{d} N\left(0, \frac{\boldsymbol{\Theta}_{11} + \beta^2\boldsymbol{\Theta}_{22} - 2\beta\boldsymbol{\Theta}_{12}}{\gamma_{\mathcal{V}^*}^{\mathsf{T}}A(\mathcal{V}^*)\gamma_{\mathcal{V}^*}}\right).$$

Note that for any $\epsilon_0 > 0$,

$$\mathbf{P}\left(\left|\sqrt{n}\left(\widehat{\beta}_L - \beta\right) - \sqrt{n}\left(\frac{\widehat{\gamma}_{\mathcal{V}^*}^{\mathsf{T}}\widehat{A}(\mathcal{V}^*)\widehat{\Gamma}_{\mathcal{V}^*}}{\widehat{\gamma}_{\mathcal{V}^*}^{\mathsf{T}}\widehat{A}(\mathcal{V}^*)\widehat{\gamma}_{\mathcal{V}^*}} - \beta\right)\mathbf{1}_{\widehat{\mathcal{V}}=\mathcal{V}^*}\right| \geq \epsilon_0\right) \leq \mathbf{P}\left(\widehat{\mathcal{V}} \neq \mathcal{V}^*\right) \quad (58)$$

and it follows from Lemma 1 that

$$\sqrt{n}\left(\widehat{\beta}_L - \beta\right) - \sqrt{n}\left(\frac{\widehat{\gamma}_{\mathcal{V}^*}^{\mathsf{T}}\widehat{A}(\mathcal{V}^*)\widehat{\Gamma}_{\mathcal{V}^*}}{\widehat{\gamma}_{\mathcal{V}^*}^{\mathsf{T}}\widehat{A}(\mathcal{V}^*)\widehat{\gamma}_{\mathcal{V}^*}} - \beta\right)\mathbf{1}_{\widehat{\mathcal{V}}=\mathcal{V}^*} \xrightarrow{p} 0. \qquad (59)$$

By Lemma 3.7 in Wooldridge [2010], we establish (19).

### C.5. Preliminary lemmas for high dimension case

We first define the following events for the random design $\mathbf{W}$ (the normalized $\mathbf{H}$) and the error $\mathbf{\Pi}$,

$$G_1 = \left\{ \frac{2}{5} \frac{1}{\sqrt{M_1}} < \frac{\|\mathbf{W}_{\cdot j}\|_2}{\sqrt{n}} < \frac{7}{5}\sqrt{M_1} \text{ for } 1 \le j \le p \right\},$$

$$G_2 = \left\{ \left| \frac{(\sigma_i^{ora})^2}{\sigma_i^2} - 1 \right| \le 2\sqrt{\frac{\log p}{n}} + 2\frac{\log p}{n} \text{ for } i = 1, 2 \right\},$$

$$G_3 = \left\{ \left| \frac{\boldsymbol{\gamma}^{\mathsf{T}} \widehat{\boldsymbol{\Sigma}} \boldsymbol{\gamma}}{\boldsymbol{\gamma}^{\mathsf{T}} \boldsymbol{\Sigma} \boldsymbol{\gamma}} - 1 \right| \le 12\sqrt{\frac{\log p}{n}} \text{ and } \left| \frac{\boldsymbol{\Omega}_{j\cdot}^{\mathsf{T}} \widehat{\boldsymbol{\Sigma}} \boldsymbol{\Omega}_{j\cdot}}{\boldsymbol{\Omega}_{jj}} - 1 \right| \le 12\sqrt{\frac{\log p}{n}}, 1 \le j \le p_z \right\},$$

$$G_4 = \left\{ \kappa(\mathbf{H}, 4s, 1 + 2\epsilon_0) \ge \frac{1}{2\sqrt{M_1}} \right\},$$

$$G_5 = \left\{ \frac{\|\mathbf{H}^{\mathsf{T}} \mathbf{\Pi}_{i\cdot}\|_\infty}{n} \le \sigma_i \sqrt{\frac{2\delta_0 \log p}{n}} \text{ for } i = 1, 2 \right\},$$

$$S_1 = \left\{ \frac{\|\mathbf{H}^{\mathsf{T}} \mathbf{\Pi}_{i\cdot}\|_\infty}{n} \le \sigma_i^{ora} \lambda_0 \frac{\epsilon_0 - 1}{\epsilon_0 + 1} (1 - \tau) \text{ for } i = 1, 2 \right\},$$

$$S_2 = \left\{ (1 - \nu_0)\hat{\sigma}_i \le \sigma_i \le (1 + \nu_0)\hat{\sigma}_i \text{ for } i = 1, 2 \right\},$$

$$\tag{60}$$

and

$$A_1 = \left\{ \|e_j^{\mathsf{T}} \boldsymbol{\Omega} \widehat{\boldsymbol{\Sigma}} - e_j^{\mathsf{T}}\|_\infty \le \lambda_n, j = 1, 2, \cdots, p_z \right\}, \text{ where } \lambda_n = 2eC_0 M_1^2 \sqrt{\frac{\log p}{n}},$$

$$A_2 = \left\{ |\hat{\boldsymbol{\gamma}}_j - \boldsymbol{\gamma}_j| \le \frac{\|\widehat{\mathbf{v}}^{[j]}\|_2 \sigma_2}{\sqrt{n}} \sqrt{2.05 \log p_z} \text{ for } 1 \le j \le p_z \right\},$$

$$A_3 = \left\{ \max_{1 \le j \le p_z} \left\| \frac{1}{n}(\hat{v}^{[j]})^{\mathsf{T}} \mathbf{\Pi}_{i\cdot} \right\|_\infty \le \left(1 + 12\sqrt{\frac{\log p}{n}}\right) M_1 \sqrt{\frac{2.05 \log p_z}{n}} \sigma_i, \text{ for } i = 1, 2 \right\},$$

$$A_4 = \left\{ \frac{2}{\sqrt{n}} \sum_{j \in \mathcal{S}^*} \boldsymbol{\gamma}_j (\widehat{\mathbf{v}}^{[j]})^{\mathsf{T}} \mathbf{\Pi}_{\cdot 2} \le \frac{2\sqrt{\log p}}{\sqrt{n}} \left\| \sum_{j \in \mathcal{S}^*} \boldsymbol{\gamma}_j \widehat{\mathbf{v}}^{[j]} \right\|_2 \sqrt{\boldsymbol{\Theta}_{22}} \right\},$$

$$A_5 = \left\{ \frac{1}{\sqrt{n}} \sum_{j \in \mathcal{S}^*} \boldsymbol{\gamma}_j (\widehat{\mathbf{v}}^{[j]})^{\mathsf{T}} (\mathbf{\Pi}_{\cdot 1} + \beta \mathbf{\Pi}_{\cdot 2}) \le \frac{\sqrt{\log p}}{\sqrt{n}} \left\| \sum_{j \in \mathcal{S}^*} \boldsymbol{\gamma}_j \widehat{\mathbf{v}}^{[j]} \right\|_2 \sqrt{\boldsymbol{\Theta}_{11} + \beta^2 \boldsymbol{\Theta}_{22} + 2\beta \boldsymbol{\Theta}_{12}} \right\},$$

$$\tag{61}$$

where $\widehat{\boldsymbol{\Sigma}} = \frac{1}{n}\mathbf{W}^{\mathsf{T}}\mathbf{W}$ and $\widehat{\mathbf{v}}^{[j]} = \mathbf{W}^{\mathsf{T}}\widehat{\mathbf{U}}_{\cdot j}$. Define

$$G = \cap_{i=1}^5 G_i, \quad S = \cap_{i=1}^2 S_i \quad \text{and} \quad A = \cap_{i=1}^5 A_i.$$

We introduce the following lemmas to control the probability of events $G$, $S$ and $A$. The detailed proofs of the following lemmas are presented in Section D.3 and D.4.

LEMMA 4. *If $s \leq cn/\log p$, then*

$$\mathbf{P}(G) \geq 1 - \frac{6}{p} - 2p^{1-C_1} - \frac{1}{2\sqrt{\pi\delta_0\log p}}p^{1-\delta_0} - 2\exp\left(-\frac{c'n}{M_1^3}\right), \quad (62)$$

*and*

$$\mathbf{P}(G \cap S) \geq \mathbf{P}(G) - 2\exp\left(-\left(\frac{g_0 + 1 - \sqrt{2g_0 + 1}}{2}\right)n\right) - c''\frac{1}{\sqrt{\log p}}p^{1-\delta_0}, \quad (63)$$

*where $g_0 = \nu_0/(2 + 3\nu_0)$ and $c, c', c_*$ and $c''$ are universal positive constants, not depending on $n$ and $p$. We also have*

$$\mathbf{P}(A_1) \geq 1 - 2p_z p^{1-c_0 C_0^2}, \quad \text{and} \quad \mathbf{P}(A_4 \cap A_5) \geq 1 - p^{-c}, \quad (64)$$

$$\min\{\mathbf{P}(A_2), \mathbf{P}(A_3)\} \geq \mathbf{P}((A_1 \cap G_1 \cap G_3)) - \frac{1}{2\sqrt{\pi\log p_z}}p_z^{-0.02}. \quad (65)$$

LEMMA 5. *On the event $A_1 \cap G_1 \cap G_3$, we have*

$$\frac{(1-\lambda_n)^2}{2M_1} \leq \frac{\|\widehat{\mathbf{v}}^{[j]}\|_2^2}{n} \leq \left(1 + 12\sqrt{\frac{\log p}{n}}\right)M_1, \quad \text{for } 1 \leq j \leq p_z. \quad (66)$$

*If $s_{z1}\sqrt{\log p/n} \to 0$, on the event $G_3$, we have*

$$\frac{1}{n}\left\|\sum_{j\in\mathcal{S}^*}\gamma_j\widehat{\mathbf{v}}^{[j]}\right\|_2^2 \geq \frac{M_1\|\boldsymbol{\gamma}\|_2^2(1-s_{z1}\lambda_n)^2}{1 - 12\sqrt{\frac{\log p}{n}}} \quad \text{and} \quad \frac{1}{n}\left\|\sum_{j\in\mathcal{V}^*}\gamma_j\widehat{\mathbf{v}}^{[j]}\right\|_2^2 \geq \frac{M_1\|\boldsymbol{\gamma}_{\mathcal{V}^*}\|_2^2(1-s_{z1}\lambda_n)^2}{1 - 12\sqrt{\frac{\log p}{n}}}. \quad (67)$$

*Furthermore, we have*

$$\frac{M_1(1-s_{z1}\lambda_n)^2}{\|\boldsymbol{\gamma}\|_2^2\left(1 - 12\sqrt{\frac{\log p}{n}}\right)}\frac{1}{M_2} \leq \mathrm{Var}_H \leq \frac{4s_{z1}M_1^2M_2(1+\beta^2)}{\|\boldsymbol{\gamma}\|_2^2}, \quad (68)$$

*and*

$$\frac{M_1(1-s_{z1}\lambda_n)^2}{\|\boldsymbol{\gamma}_{\mathcal{V}^*}\|_2^2\left(1 - 12\sqrt{\frac{\log p}{n}}\right)}\frac{1}{M_2} \leq \mathrm{Var} \leq \frac{4s_{z1}M_1^2M_2(1+\beta^2)}{\|\boldsymbol{\gamma}_{\mathcal{V}^*}\|_2^2}. \quad (69)$$

### C.6.  Proof of Theorem 4

The proof of Theorem 4 is based on Lemma 6 and the following expression for the estimator $\widehat{\beta}_H$, $\widehat{\beta}_H = \widehat{\boldsymbol{\gamma}^{\mathsf{T}}\boldsymbol{\Gamma}}/\widehat{\|\boldsymbol{\gamma}\|}_2^2$, where $\widehat{\|\boldsymbol{\gamma}\|}_2^2 = \sum_{j\in\widehat{\mathcal{S}}}\widehat{\gamma}_j^2$ and $\widehat{\boldsymbol{\gamma}^{\mathsf{T}}\boldsymbol{\Gamma}} = \sum_{j\in\widehat{\mathcal{S}}}\widehat{\gamma}_j\widehat{\boldsymbol{\Gamma}}_j$.

LEMMA 6. *Suppose that $\sqrt{s_{z1}}s\log p/\sqrt{n} \to 0$, $\pi^* = 0$ and the assumptions* (R1)− (R3) *hold. Then we have the following decompositions,*

$$\sqrt{n}\left(\widehat{\|\boldsymbol{\gamma}\|_2^2} - \|\boldsymbol{\gamma}\|_2^2\right) = \frac{2}{\sqrt{n}}\sum_{j\in\mathcal{S}^*}\boldsymbol{\gamma}_j(\widehat{\mathbf{v}}^{[j]})^{\mathsf{T}}\boldsymbol{\Pi}_{\cdot 2} + R^{\gamma}, \tag{70}$$

*and*

$$\sqrt{n}\left(\widehat{\boldsymbol{\gamma}^{\mathsf{T}}\boldsymbol{\Gamma}} - \boldsymbol{\gamma}^{\mathsf{T}}\boldsymbol{\Gamma}\right) = \frac{1}{\sqrt{n}}\sum_{j\in\mathcal{S}^*}\boldsymbol{\gamma}_j(\widehat{\mathbf{v}}^{[j]})^{\mathsf{T}}\left(\boldsymbol{\Pi}_{\cdot 1} + \beta\boldsymbol{\Pi}_{\cdot 2}\right) + R^{\mathrm{inter}}, \tag{71}$$

*where*

$$\frac{2}{\sqrt{n}}\sum_{j\in\mathcal{S}^*}\boldsymbol{\gamma}_j(\widehat{\mathbf{v}}^{[j]})^{\mathsf{T}}\boldsymbol{\Pi}_{\cdot 2} \sim N\left(0, \frac{4}{n}\left\|\sum_{j\in\mathcal{S}^*}\boldsymbol{\gamma}_j\widehat{\mathbf{v}}^{[j]}\right\|_2^2\boldsymbol{\Theta}_{22}\right), \tag{72}$$

$$\frac{1}{\sqrt{n}}\sum_{j\in\mathcal{S}^*}\boldsymbol{\gamma}_j(\widehat{\mathbf{v}}^{[j]})^{\mathsf{T}}\left(\boldsymbol{\Pi}_{\cdot 1} + \beta\boldsymbol{\Pi}_{\cdot 2}\right) \sim N\left(0, \frac{1}{n}\left\|\sum_{j\in\mathcal{S}^*}\boldsymbol{\gamma}_j\widehat{\mathbf{v}}^{[j]}\right\|_2^2\left(\boldsymbol{\Theta}_{11} + \beta^2\boldsymbol{\Theta}_{22} + 2\beta\boldsymbol{\Theta}_{12}\right)\right), \tag{73}$$

*and on the event $A \cap S \cap G$, we have*

$$\max\left\{|R^{\gamma}|, |R^{\mathrm{inter}}|\right\} \le C\left(|\beta|+1\right)\|\boldsymbol{\gamma}\|_2\sqrt{s_{z1}}s\frac{\log p}{\sqrt{n}} + Cs_{z1}\frac{\log p_z}{\sqrt{n}}. \tag{74}$$

*Then on the event $A \cap S \cap G$, we have*

$$\max\left\{\left|\widehat{\|\boldsymbol{\gamma}\|_2^2} - \|\boldsymbol{\gamma}\|_2^2\right|, \left|\widehat{\boldsymbol{\gamma}^{\mathsf{T}}\boldsymbol{\Gamma}} - \boldsymbol{\gamma}^{\mathsf{T}}\boldsymbol{\Gamma}\right|\right\} \le C\|\boldsymbol{\gamma}\|_2 s_{z1}\sqrt{\frac{\log p}{n}} + Cs_{z1}\frac{\log p_z}{n} \le C\|\boldsymbol{\gamma}\|_2 s_{z1}\sqrt{\frac{\log p}{n}}. \tag{75}$$

In the following, we will prove (27) in the main paper. Note that

$$\widehat{\beta} - \beta = -\frac{\beta}{\|\boldsymbol{\gamma}\|_2^2}\left(\widehat{\|\boldsymbol{\gamma}\|_2^2} - \|\boldsymbol{\gamma}\|_2^2\right) + \frac{1}{\|\boldsymbol{\gamma}\|_2^2}\left(\widehat{\boldsymbol{\gamma}^{\mathsf{T}}\boldsymbol{\Gamma}} - \boldsymbol{\gamma}^{\mathsf{T}}\boldsymbol{\Gamma}\right) + \frac{\|\boldsymbol{\gamma}\|_2^2 - \widehat{\|\boldsymbol{\gamma}\|_2^2}}{\|\boldsymbol{\gamma}\|_2^2}\left(\frac{\widehat{\boldsymbol{\gamma}^{\mathsf{T}}\boldsymbol{\Gamma}}}{\widehat{\|\boldsymbol{\gamma}\|_2^2}} - \frac{\boldsymbol{\gamma}^{\mathsf{T}}\boldsymbol{\Gamma}}{\|\boldsymbol{\gamma}\|_2^2}\right). \tag{76}$$

By Lemma 6, we have the following decomposition,

$$\sqrt{n}\left(\widehat{\beta} - \beta\right) = T^{\beta} + \Delta^{\beta}, \tag{77}$$

*where*

$$\begin{aligned}
T^{\beta} &= -\frac{\beta}{\|\boldsymbol{\gamma}\|_2^2}\frac{2}{\sqrt{n}}\sum_{j\in\mathcal{S}^*}\boldsymbol{\gamma}_j(\widehat{\mathbf{v}}^{[j]})^{\mathsf{T}}\boldsymbol{\Pi}_{\cdot 2} + \frac{1}{\|\boldsymbol{\gamma}\|_2^2}\frac{1}{\sqrt{n}}\sum_{j\in\mathcal{S}^*}\boldsymbol{\gamma}_j(\widehat{\mathbf{v}}^{[j]})^{\mathsf{T}}\left(\boldsymbol{\Pi}_{\cdot 1} + \beta\boldsymbol{\Pi}_{\cdot 2}\right) \\
&= \frac{1}{\|\boldsymbol{\gamma}\|_2^2}\frac{1}{\sqrt{n}}\sum_{j\in\mathcal{S}^*}\boldsymbol{\gamma}_j(\widehat{\mathbf{v}}^{[j]})^{\mathsf{T}}\left(\boldsymbol{\Pi}_{\cdot 1} - \beta\boldsymbol{\Pi}_{\cdot 2}\right),
\end{aligned}$$

and $\Delta^\beta = \mathrm{Res}_1 + \mathrm{Res}_2$ with

$$\mathrm{Res}_1 = \frac{1}{\|\boldsymbol{\gamma}\|_2^2}\left(-\beta R^{\boldsymbol{\gamma}} + R^{\mathrm{inter}}\right) \text{ and } \mathrm{Res}_2 = \sqrt{n}\frac{\|\boldsymbol{\gamma}\|_2^2 - \widehat{\|\boldsymbol{\gamma}\|_2^2}}{\|\boldsymbol{\gamma}\|_2^2}\left(\frac{\widehat{\boldsymbol{\gamma^\intercal\Gamma}}}{\widehat{\|\boldsymbol{\gamma}\|_2^2}} - \frac{\boldsymbol{\gamma^\intercal\Gamma}}{\|\boldsymbol{\gamma}\|_2^2}\right).$$

By the distribution of $\boldsymbol{\Pi}$, we establish that

$$T^\beta \mid \mathbf{W} \sim N\left(0, \frac{1}{n\|\boldsymbol{\gamma}\|_2^4}\left\|\sum_{j\in\mathcal{S}^*}\boldsymbol{\gamma}_j\widehat{\mathbf{v}}^{[j]}\right\|_2^2\left(\boldsymbol{\Theta}_{11} + \beta^2\boldsymbol{\Theta}_{22} - 2\beta\boldsymbol{\Theta}_{12}\right)\right). \qquad (78)$$

By Lemma 6, on the event $G \cap S \cap A$, we have

$$\frac{1}{\sqrt{\mathrm{Var}_H}}\left|\mathrm{Res}_1\right| \le C\frac{1}{\|\boldsymbol{\gamma}\|_2}\left(|\beta|\left|R^{\boldsymbol{\gamma}}\right| + \left|R^{\mathrm{inter}}\right|\right) \le C\left(|\beta| + 1\right)\sqrt{s_{z1}}s\frac{\log p}{\sqrt{n}} + C\frac{1}{\|\boldsymbol{\gamma}\|_2}\frac{s_{z1}\log p}{\sqrt{n}}. \qquad (79)$$

Note that on the event $G \cap S \cap A$,

$$\frac{1}{\sqrt{\mathrm{Var}_H}}\mathrm{Res}_2 \le C\sqrt{n}\frac{\|\boldsymbol{\gamma}\|_2^2 - \widehat{\|\boldsymbol{\gamma}\|_2^2}}{\|\boldsymbol{\gamma}\|_2}\times\frac{\left(\widehat{\boldsymbol{\gamma^\intercal\Gamma}} - \boldsymbol{\gamma^\intercal\Gamma}\right) + \beta\left(\|\boldsymbol{\gamma}\|_2^2 - \widehat{\|\boldsymbol{\gamma}\|_2^2}\right)}{\|\boldsymbol{\gamma}\|_2^2 + \left(\widehat{\|\boldsymbol{\gamma}\|_2^2} - \|\boldsymbol{\gamma}\|_2^2\right)} \le C\frac{s_{z1}^3(\log p)^{\frac{3}{2}}}{n},$$
$$(80)$$

where the last inequality follows from (75). Combined with (79), by $\sqrt{s_{z1}}s\log p/\sqrt{n} \to 0$, we can establish that on the event $G \cap S \cap A$,

$$\left|\Delta^\beta/\sqrt{\mathrm{Var}_H}\right| \le C\sqrt{s_{z1}}s\frac{\log p}{\sqrt{n}} + C\frac{1}{\|\boldsymbol{\gamma}\|_2}\frac{s_{z1}\log p}{\sqrt{n}}. \qquad (81)$$

Since $\sqrt{s_{z1}}s\log p/\sqrt{n} \to 0$, we establish $\Delta^\beta/\sqrt{\mathrm{Var}_H} \xrightarrow{p} 0$. Combined with (78), we establish (27).

## C.7.   Proof of Theorem 5
We first introduce the following lemma to establish the coverage property.

LEMMA 7. *Suppose that $\boldsymbol{\pi}^* = 0$ and the assumptions* (R1) $-$ (R3) *hold. As* $\sqrt{s_{z1}}s\log p/\sqrt{n} \to 0$, *then we have*

$$\frac{\widehat{\mathrm{Var}}_H}{\mathrm{Var}_H} \xrightarrow{p} 1. \qquad (82)$$

By (78), we have $\frac{T^\beta}{\sqrt{\mathrm{Var}_H}} \sim N(0,1)$. Combined with (81) and Lemma 7, we have

$$\sqrt{n}\frac{\widehat{\beta}_H - \beta}{\sqrt{\widehat{\mathrm{Var}}_H}} = \frac{T^\beta + \Delta^\beta}{\sqrt{\mathrm{Var}_H}}\times\frac{\sqrt{\mathrm{Var}_H}}{\sqrt{\widehat{\mathrm{Var}}_H}} \xrightarrow{d} N(0,1). \qquad (83)$$

and hence the coverage property (28) follows.

### C.8.  Proof of Theorem 3

The proof of the theorem follows from Lemma 2, which characterizes the behavior of the selection process described in Section 3.2 to 3.4 the main paper when it is applied to the high dimensional setting. We have the following decomposition,

$$\widehat{\beta} - \beta = \left( \frac{\sum_{j \in \mathcal{V}^*} \widehat{\gamma}_j \widehat{\Gamma}_j}{\sum_{j \in \mathcal{V}^*} \widehat{\gamma}_j^2} - \frac{\sum_{j \in \mathcal{V}^*} \gamma_j \Gamma_j}{\sum_{j \in \mathcal{V}^*} \gamma_j^2} \right).$$

Based on the this expression, (23) in the main paper will follow the same argument with (27), which is presented in Section C.6. We introduce the following lemma to establish the coverage property.

LEMMA 8. *Suppose the assumptions* (R1)$-$(R5) *and* (IN1$-$P), (IN2) *and* (IN3) *are satisfied. As* $\sqrt{s_{z1}} s \log p / \sqrt{n} \to 0$, *we have*

$$\frac{\widehat{\mathrm{Var}}}{\mathrm{Var}} \xrightarrow{p} 1. \tag{84}$$

Similarly to the proof of Theorem 5 in Section C.7, we establish

$$\sqrt{n} \frac{\widehat{\beta} - \beta}{\sqrt{\widehat{\widehat{\mathrm{Var}}}}} = \frac{T^\beta + \Delta^\beta}{\sqrt{\mathrm{Var}}} \times \frac{\sqrt{\mathrm{Var}}}{\sqrt{\widehat{\widehat{\mathrm{Var}}}}} \xrightarrow{d} N(0,1), \tag{85}$$

and hence establish the coverage property (24) in the main paper.

## D.  Supplementary Materials: Proof of Extra Lemmas

In this section, we prove extra lemmas used in the proof of main theorems.

### D.1.  Proof of Lemma 3

Note that

$$\sqrt{n} \left( \frac{\widehat{\gamma}_{\mathcal{V}^*}^{\mathsf{T}} \widehat{A}(\mathcal{V}^*) \widehat{\Gamma}_{\mathcal{V}^*}}{\widehat{\gamma}_{\mathcal{V}^*}^{\mathsf{T}} \widehat{A}(\mathcal{V}^*) \widehat{\gamma}_{\mathcal{V}^*}} - \beta \right) = \frac{\widehat{\gamma}_{\mathcal{V}^*}^{\mathsf{T}} \widehat{A}(\mathcal{V}^*) \left( \widehat{\Sigma}^{-1} \right)_{\mathcal{V}^*}}{\widehat{\gamma}_{\mathcal{V}^*}^{\mathsf{T}} \widehat{A}(\mathcal{V}^*) \widehat{\gamma}_{\mathcal{V}^*}} \frac{1}{\sqrt{n}} \mathbf{W}^{\mathsf{T}} \left( \mathbf{\Pi}_{\cdot 1} - \beta \mathbf{\Pi}_{\cdot 2} \right). \tag{86}$$

Since

$$\frac{\widehat{\gamma}_{\mathcal{V}^*}^{\mathsf{T}} \widehat{A}(\mathcal{V}^*) \left( \widehat{\Sigma}^{-1} \right)_{\mathcal{V}^*}}{\widehat{\gamma}_{\mathcal{V}^*}^{\mathsf{T}} \widehat{A}(\mathcal{V}^*) \widehat{\gamma}_{\mathcal{V}^*}} \xrightarrow{p} \frac{\gamma_{\mathcal{V}^*}^{\mathsf{T}} A(\mathcal{V}^*) \left( \Sigma^{-1} \right)_{\mathcal{V}^*}}{\gamma_{\mathcal{V}^*}^{\mathsf{T}} A(\mathcal{V}^*) \gamma_{\mathcal{V}^*}} \tag{87}$$

and

$$\frac{1}{\sqrt{n}} \mathbf{W}^{\mathsf{T}} \left( \mathbf{\Pi}_{\cdot 1} - \beta \mathbf{\Pi}_{\cdot 2} \right) \xrightarrow{d} N \left( 0, \left( \mathbf{\Theta}_{11} + \beta^2 \mathbf{\Theta}_{22} - 2\beta \mathbf{\Theta}_{12} \right) \mathbf{\Sigma} \right), \tag{88}$$

we have

$$
\begin{aligned}
&\frac{\widehat{\gamma}_{\mathcal{V}^*}^{\mathsf{T}} \widehat{A}(\mathcal{V}^*) \left(\widehat{\boldsymbol{\Sigma}}^{-1}\right)_{\mathcal{V}^*} }{\widehat{\gamma}_{\mathcal{V}^*}^{\mathsf{T}} \widehat{A}(\mathcal{V}^*)\widehat{\gamma}_{\mathcal{V}^*}} \frac{1}{\sqrt{n}} \mathbf{W}^{\mathsf{T}} \left(\boldsymbol{\Pi}_{\cdot 2} - \beta \boldsymbol{\Pi}_{\cdot 1}\right) \\
&\xrightarrow{d} N\left(0, \frac{\gamma_{\mathcal{V}^*}^{\mathsf{T}} A(\mathcal{V}^*) \left(\boldsymbol{\Sigma}^{-1}\right)_{\mathcal{V}^* \mathcal{V}^*} A(\mathcal{V}^*)\gamma_{\mathcal{V}^*}}{\left(\gamma_{\mathcal{V}^*}^{\mathsf{T}} A(\mathcal{V}^*)\gamma_{\mathcal{V}^*}\right)^2} \left(\boldsymbol{\Theta}_{11} + \beta^2 \boldsymbol{\Theta}_{22} - 2\beta \boldsymbol{\Theta}_{12}\right)\right)
\end{aligned}
\tag{89}
$$

Since $A(\mathcal{V}^*) \left(\boldsymbol{\Sigma}^{-1}\right)_{\mathcal{V}^* \mathcal{V}^*} A(\mathcal{V}^*) = A(\mathcal{V}^*)$, we establish (56).

### D.2.    Lemmas for scaled Lasso and de-biasing Lasso

We introduce the following lemmas for scaled Lasso and de-biasing Lasso used in the later proofs. Lemma 9 establishes the convergence rate of the scaled Lasso method, which is based on the analysis in Sun and Zhang [2012].

LEMMA 9.  *On the event $G \cap S$, if $s \leq cn/\log p$, then*

$$
\|\widetilde{\boldsymbol{\Gamma}} - \boldsymbol{\Gamma}\|_1 + \|\widetilde{\boldsymbol{\Psi}} - \boldsymbol{\Psi}\|_1 \leq Cs\sqrt{\frac{\log p}{n}}\sigma_1, \quad \|\widetilde{\boldsymbol{\gamma}} - \boldsymbol{\gamma}\|_1 + \|\widetilde{\boldsymbol{\psi}} - \boldsymbol{\psi}\|_1 \leq Cs\sqrt{\frac{\log p}{n}}\sigma_2, \tag{90}
$$

$$
\frac{1}{\sqrt{n}}\|\mathbf{Z}(\widetilde{\boldsymbol{\Gamma}} - \boldsymbol{\Gamma}) + \mathbf{X}(\widetilde{\boldsymbol{\Psi}} - \boldsymbol{\Psi})\|_2 \leq C\sqrt{\frac{s\log p}{n}}\sigma_1, \tag{91}
$$

*and*

$$
\frac{1}{\sqrt{n}}\|\mathbf{Z}(\widetilde{\boldsymbol{\gamma}} - \boldsymbol{\gamma}) + \mathbf{X}(\widetilde{\boldsymbol{\psi}} - \boldsymbol{\psi})\|_2 \leq C\sqrt{\frac{s\log p}{n}}\sigma_2. \tag{92}
$$

The following lemma is the key result for the de-biasing Lasso estimator, established in Zhang and Zhang [2014], Javanmard and Montanari [2014], van de Geer et al. [2014].

LEMMA 10.  *We have the following expressions for the proposed de-biased estimator,*

$$
\widehat{\boldsymbol{\Gamma}} - \boldsymbol{\Gamma} = D^{\boldsymbol{\Gamma}} + \Delta^{\boldsymbol{\Gamma}}, \tag{93}
$$

*where*

$$
D_j^{\boldsymbol{\Gamma}} = \frac{1}{n}(\widehat{\mathbf{v}}^{[j]})^{\mathsf{T}} \boldsymbol{\Pi}_{\cdot 1} \quad and \quad \Delta_j^{\boldsymbol{\Gamma}} = \left(\frac{1}{n}(\widehat{\mathbf{U}}_{\cdot j})^{\mathsf{T}}\widehat{\boldsymbol{\Sigma}} - e_j^{\mathsf{T}}\right)\begin{pmatrix}\widehat{\boldsymbol{\Gamma}} - \boldsymbol{\Gamma} \\ \widehat{\boldsymbol{\Psi}} - \boldsymbol{\Psi}\end{pmatrix}, \; i = 1, \cdots, p_z. \tag{94}
$$

*We also have*

$$
\widehat{\boldsymbol{\gamma}} - \boldsymbol{\gamma} = D^{\boldsymbol{\gamma}} + \Delta^{\boldsymbol{\gamma}}, \tag{95}
$$

*where*

$$
D_j^{\boldsymbol{\gamma}} = \frac{1}{n}(\widehat{\mathbf{v}}^{[j]})^{\mathsf{T}} \boldsymbol{\Pi}_{\cdot 2} \quad and \quad \Delta_j^{\boldsymbol{\gamma}} = \left(\frac{1}{n}(\widehat{\mathbf{U}}_{\cdot j})^{\mathsf{T}}\widehat{\boldsymbol{\Sigma}} - e_j^{\mathsf{T}}\right)\begin{pmatrix}\widehat{\boldsymbol{\gamma}} - \boldsymbol{\gamma} \\ \widehat{\boldsymbol{\psi}} - \boldsymbol{\psi}\end{pmatrix}, \; i = 1, \cdots, p_z. \tag{96}
$$

*On the event $S \cap G \cap A$, we have*

$$
\max\left\{\|\Delta^{\boldsymbol{\gamma}}\|_\infty, \|\Delta^{\boldsymbol{\Gamma}}\|_\infty\right\} \leq Cs\frac{\log p}{n}\max\left\{\sigma_1, \sigma_2\right\}. \tag{97}
$$

### D.3. Proof of Lemma 4

The proof of Lemma 4 is a generalization of Lemma 4 in Cai and Guo [2017]. In the following, we extend the Gaussian design in Cai and Guo [2017] to sub-gaussian design considered in this paper. Since the error of the regression is still assumed to be Gaussian, it is sufficient to establish the probability bound of $G_1, G_3, G_4$ and $A_1$ for the sub-gaussian design matrix and control the events $A_2$ and $A_3$. The probability bound of the event $A_1$ for the sub-gaussian design is established in Lemma 4 of Cai and Guo [2016a]. By Corollary 5.17 in Vershynin [2012] and the union bound, we have

$$\mathbf{P}\left(\max_{1 \le j \le p}\left|\frac{1}{n}\left(\|\mathbf{W}_{\cdot j}\|_2^2 - \mathbf{E}\|\mathbf{W}_{\cdot j}\|_2^2\right)\right| \ge \epsilon\right) \le 2p\exp\left(-\frac{1}{6}\min\left\{\frac{\epsilon^2}{K^2}, \frac{\epsilon}{K}\right\}n\right),$$

where $K = 4M_1$. Taking $\epsilon = 12M_1\sqrt{\log p/n}$, we have

$$\mathbf{P}\left(\max_{1 \le j \le p}\left|\frac{1}{n}\left(\|\mathbf{W}_{\cdot j}\|_2^2 - \mathbf{E}\|\mathbf{W}_{\cdot j}\|_2^2\right)\right| \ge 12M_1\sqrt{\frac{\log p}{n}}\right) \le 2p^{-\frac{1}{2}} \quad \text{and} \quad \mathbf{P}(G_1) \ge 1 - 2p^{-\frac{1}{2}}.$$
$$(98)$$

Similarly, we have $\mathbf{P}\left(\left|\frac{1}{n}\left(\|\mathbf{W}u\|_2^2 - \mathbf{E}\|\mathbf{W}u\|_2^2\right)\right| \ge 12M_1\|u\|^2\sqrt{\log p/n}\right) \le 2p^{-\frac{3}{2}}$, and

$$\mathbf{P}\left(\left|\frac{1}{n}\left(\frac{\|\mathbf{W}u\|_2^2}{\mathbf{E}\|\mathbf{W}u\|_2^2} - 1\right)\right| \ge 12M_1\frac{\|u\|^2}{\mathbf{E}\|\mathbf{W}u\|_2^2}\sqrt{\log p/n}\right) \le 2p^{-\frac{3}{2}},$$

and hence

$$\mathbf{P}(G_3) \ge 1 - (p_z + 2)p^{-\frac{3}{2}}.$$

By Theorem 1.6 in Zhou [2009], if $n \ge 1/\theta^2 \times c'M_1^3 \max\left\{12(2+\gamma_0)^2 M_1 s \log(5ep/4s), 9\log p\right\}$, then with probability at least $1 - 2\exp\left(-c\theta^2 n/M_1^3\right)$, for all $\delta$ such that there exist $|J_0| \le 4s$ and $\|\delta_{J_0^c}\|_1 \le \gamma_0\|\delta_{J_0}\|_1$, we have $\|Z\delta\|_2/(\sqrt{n}\|\Sigma^{\frac{1}{2}}\delta\|_2) \ge 1 - \theta$. By taking $\theta = \frac{1}{2}$, if $n \ge 4c'M_1^3 \max\left\{12(2+\gamma_0)^2 M_1 s \log(5ep/4s), 9\log p\right\}$, then $\mathbf{P}(G_4) \ge 1 - 2\exp\left(-cn/M_1^3\right)$. In the following, we control the events $A_2$ and $A_3$,

$$\mathbf{P}(A_2^c) \le \mathbf{P}\left(\max_{1 \le i \le q}\frac{\left|D_j^\gamma\right|}{\sqrt{\operatorname{Var}\left(D_j^\gamma\right)}} \ge \sqrt{2.02\log p_z}\right) + \mathbf{P}\left(\max_{1 \le j \le p_z}\frac{\left|\Delta_j^\gamma\right|}{\sqrt{\operatorname{Var}\left(D_j^\gamma\right)}} \ge 0.01\sqrt{\log p_z}\right)$$
$$\le \frac{1}{2\sqrt{\pi\log p_z}}p_z^{-0.02} + \mathbf{P}\left((S \cap G \cap A_1)^c\right),$$

where the first inequality follows from (95) and the second inequality follows from (96) and (97). The control of $\mathbf{P}(A_4 \cap A_5)$ follows from (72) and (73). Note that

$$\mathbf{P}(A_3^c) \le \mathbf{P}\left((A_1 \cap G_1 \cap G_3)^c\right) + \mathbf{P}\left(A_3^c \cap A_1 \cap G_1 \cap G_3\right)$$
$$\le \mathbf{P}\left((A_1 \cap G_1 \cap G_3)^c\right) + 2\mathbf{P}\left(\max_{1 \le j \le p_z}\frac{1}{\|\widehat{\mathbf{v}}^{[j]}\|_2\sigma_1}\left|(\widehat{\mathbf{v}}^{[j]})^\intercal\mathbf{\Pi}_{\cdot 1}\right| \ge \sqrt{2.05\log p_z}\right)$$
$$\le \mathbf{P}\left((A_1 \cap G_1 \cap G_3)^c\right) + \frac{1}{\sqrt{\pi\log p_z}}p_z^{-0.05},$$

where the second inequality follows from (66) and the last inequality follows from the fact that $1/(\|\widehat{\mathbf{v}}^{[j]}\|_2\sigma_1) \times (\widehat{\mathbf{v}}^{[j]})^\intercal\mathbf{\Pi}_{\cdot 1}$ conditioning on $\mathbf{W}$ is normally distributed.

### D.4.  Proof of Lemma 5

In the following, we only establish the results for $\widehat{\mathbf{v}}^{[1]}$ and the same argument extends to $\widehat{\mathbf{v}}^{[j]}$ where $1 \leq j \leq p_z$. Since $\lambda_n = 2eC_0M_1^2\sqrt{\log p/n}$ is chosen such that $\mathbf{\Omega}_{1\cdot}$ belongs to the feasible set, we have

$$\frac{\|\widehat{\mathbf{v}}^{[1]}\|_2^2}{n} \leq \frac{\|\mathbf{W}\mathbf{\Omega}_{1\cdot}\|_2^2}{n}. \tag{99}$$

By Lemma 12 in Javanmard and Montanari [2014], we have

$$\frac{\|\widehat{\mathbf{v}}^{[1]}\|_2^2}{n} \geq \frac{(1-\lambda_n)^2}{\widehat{\mathbf{\Sigma}}_{11}}. \tag{100}$$

By the definition of $G_1$ and $G_3$, we establish (66). Let $\mathcal{I} = \{1, 2, \cdots, p_z\}$ and assume that $M \in \mathbb{R}^{p\times p_z}$ belongs to the feasible set $\|\widehat{\mathbf{\Sigma}}\mathbf{\Omega} - \mathbf{I}_{\cdot\mathcal{I}}\|_\infty \leq \lambda_n$, where $\mathbf{I}_{\cdot\mathcal{I}}$ denotes the sub-matrix of the identity matrix containing the column with index $i \in \mathcal{I}$, that is, $\|\widehat{\mathbf{\Sigma}}M - \mathbf{I}_{\cdot\mathcal{I}}\|_\infty \leq \lambda_n$, and hence

$$\|\widehat{\mathbf{\Sigma}}M\boldsymbol{\gamma} - \boldsymbol{\gamma}\|_\infty = \|\left(\widehat{\mathbf{\Sigma}}M - \mathbf{I}_{\cdot\mathcal{I}}\right)\boldsymbol{\gamma}\|_\infty \leq \|\widehat{\mathbf{\Sigma}}M - \mathbf{I}_{\cdot\mathcal{I}}\|_\infty\|\boldsymbol{\gamma}\|_1 \leq \lambda_n\|\boldsymbol{\gamma}\|_1. \tag{101}$$

Note that

$$\left|\boldsymbol{\gamma}^\intercal\widehat{\mathbf{\Sigma}}M\boldsymbol{\gamma} - \|\boldsymbol{\gamma}\|_2^2\right| = \left|\boldsymbol{\gamma}^\intercal\left(\widehat{\mathbf{\Sigma}}M\boldsymbol{\gamma} - \boldsymbol{\gamma}\right)\right| \leq \|\boldsymbol{\gamma}\|_1\|\widehat{\mathbf{\Sigma}}M\boldsymbol{\gamma} - \boldsymbol{\gamma}\|_\infty \leq \lambda_n\|\boldsymbol{\gamma}\|_1^2, \tag{102}$$

where the last inequality follows from (101). The inequality (102) informs that $M\boldsymbol{\gamma}$ is in the feasible set

$$\left|\boldsymbol{\gamma}\widehat{\mathbf{\Sigma}}\left(M\boldsymbol{\gamma}\right) - \|\boldsymbol{\gamma}\|_2^2\right| \leq \lambda_n\|\boldsymbol{\gamma}\|_1^2. \tag{103}$$

We define $\mu^*$ as

$$\begin{aligned}\mu^* = \arg\min_\mu \mu^\intercal\widehat{\mathbf{\Sigma}}\mu \\ \text{subject to} \quad \left|\boldsymbol{\gamma}^\intercal\widehat{\mathbf{\Sigma}}\mu - \|\boldsymbol{\gamma}\|_2^2\right| \leq \lambda_n\|\boldsymbol{\gamma}\|_1^2\end{aligned} \tag{104}$$

By (103), we have the following inequality,

$$\frac{1}{n}\left\|\sum_{j\in\mathcal{S}^*}\gamma_j\widehat{\mathbf{v}}^{[j]}\right\|_2^2 = \boldsymbol{\gamma}^\intercal M^\intercal\widehat{\mathbf{\Sigma}}M\boldsymbol{\gamma} \geq (\mu^*)^\intercal\widehat{\mathbf{\Sigma}}\mu^*. \tag{105}$$

In the following, we will show that $(\mu^*)^\intercal\widehat{\mathbf{\Sigma}}\mu^* = \langle\mu^*, \widehat{\mathbf{\Sigma}}\mu^*\rangle$ is further lower bounded. Since $\mu^*$ is feasible in the constrained set of (104), we have $\|\boldsymbol{\gamma}\|_2^2 - \boldsymbol{\gamma}^\intercal\widehat{\mathbf{\Sigma}}\mu^* - \lambda_n\|\boldsymbol{\gamma}\|_1^2 \leq 0$, and hence for any positive constant $c > 0$, we have

$$\begin{aligned}\langle\mu^*, \widehat{\mathbf{\Sigma}}\mu^*\rangle &\geq \langle\mu^*, \widehat{\mathbf{\Sigma}}\mu^*\rangle + c\left(\|\boldsymbol{\gamma}\|_2^2 - \boldsymbol{\gamma}^\intercal\widehat{\mathbf{\Sigma}}\mu^* - \lambda_n\|\boldsymbol{\gamma}\|_1^2\right)\\ &\geq \min_\mu\left(\langle\mu, \widehat{\mathbf{\Sigma}}\mu\rangle + c\left(\|\boldsymbol{\gamma}\|_2^2 - \boldsymbol{\gamma}^\intercal\widehat{\mathbf{\Sigma}}\mu - \lambda_n\|\boldsymbol{\gamma}\|_1^2\right)\right) = -\frac{c^2}{4}\langle\boldsymbol{\gamma}, \widehat{\mathbf{\Sigma}}\boldsymbol{\gamma}\rangle + c\left(\|\boldsymbol{\gamma}\|_2^2 - \lambda_n\|\boldsymbol{\gamma}\|_1^2\right).\end{aligned} \tag{106}$$

Note that $\|\gamma\|_1^2 \lambda_n \le s_{z1} \lambda_n \|\gamma\|_2^2 = C s_{z1} \sqrt{\log p/n} \|\gamma\|_2^2 \ll \|\gamma\|_2^2$, where the last inequality holds when $s_{z1} \sqrt{\log p/n} \to 0$. By (106), we have

$$
\begin{aligned}
\langle \mu^*, \widehat{\boldsymbol{\Sigma}} \mu^* \rangle &\ge \max_{c>0} -\frac{c^2}{4} \langle \gamma, \widehat{\boldsymbol{\Sigma}} \gamma \rangle + c \left( \|\gamma\|_2^2 - \lambda_n \|\gamma\|_1^2 \right) \\
&= \frac{\left( \|\gamma\|_2^2 - \lambda_n \|\gamma\|_1^2 \right)^2}{\langle \gamma, \widehat{\boldsymbol{\Sigma}} \gamma \rangle} \ge \frac{\|\gamma\|_2^4 \left( 1 - s_{z1} \lambda_n \right)^2}{\langle \gamma, \widehat{\boldsymbol{\Sigma}} \gamma \rangle}.
\end{aligned}
\tag{107}
$$

On the event $G_3$, we establish (67) for $\left\| \sum_{j \in \mathcal{S}^*} \gamma_j \widehat{\mathbf{v}}^{[j]} \right\|_2^2 / n$. The same argument holds for $\left\| \sum_{j \in \mathcal{V}^*} \gamma_j \widehat{\mathbf{v}}^{[j]} \right\|_2^2 / n$. Note that

$$
\frac{1}{M_2} \le \boldsymbol{\Theta}_{11} + \beta^2 \boldsymbol{\Theta}_{22} - 2\beta \boldsymbol{\Theta}_{12} = \begin{pmatrix} 1 & -\beta \end{pmatrix} \boldsymbol{\Theta} \begin{pmatrix} 1 \\ -\beta \end{pmatrix} \le M_2 \left( 1 + \beta^2 \right).
\tag{108}
$$

Combined with (67), we establish the first inequality of (68). Note that

$$
\frac{1}{n} \left\| \sum_{j \in \mathcal{S}^*} \gamma_j \widehat{\mathbf{v}}^{[j]} \right\|_2^2 \le \left( 2M_2 \sum_{j \in \mathcal{V}^*} |\gamma_j| \right)^2 \le s_{z1} \|\gamma\|_2^2.
\tag{109}
$$

Combined with (108), we establish the second inequality of (68). By the similar argument, we can establish (69).

### D.5. Proof of Lemma 6

In the following proof, we will use the shorthand $\langle a, b \rangle_J = \sum_{j \in J} a_j b_j$. We have the following decompositions for $\widehat{\|\gamma\|_2^2} - \|\gamma\|_2^2$ and $\widehat{\gamma^\intercal \boldsymbol{\Gamma}} - \gamma^\intercal \boldsymbol{\Gamma}$,

$$
\begin{aligned}
\widehat{\|\gamma\|_2^2} - \|\gamma\|_2^2 =& 2\langle \gamma, D^\gamma \rangle_{\widehat{\mathcal{S}}} + 2\langle \gamma, \Delta^\gamma \rangle_{\widehat{\mathcal{S}}} + \langle D^\gamma, D^\gamma \rangle_{\widehat{\mathcal{S}}} + \langle \Delta^\gamma, \Delta^\gamma \rangle_{\widehat{\mathcal{S}}} + 2\langle D^\gamma, \Delta^\gamma \rangle_{\widehat{\mathcal{S}}} \\
&- \left( \sum_{j \in \mathcal{S}^* \setminus \widehat{\mathcal{S}}} \gamma_j^2 - \sum_{j \in \widehat{\mathcal{S}} \setminus \mathcal{S}^*} \gamma_j^2 \right),
\end{aligned}
\tag{110}
$$

and

$$
\begin{aligned}
\widehat{\gamma^\intercal \boldsymbol{\Gamma}} - \gamma^\intercal \boldsymbol{\Gamma} =& \langle \gamma, D^{\boldsymbol{\Gamma}} \rangle_{\widehat{\mathcal{S}}} + \langle \boldsymbol{\Gamma}, D^\gamma \rangle_{\widehat{\mathcal{S}}} + \langle \gamma, \Delta^{\boldsymbol{\Gamma}} \rangle_{\widehat{\mathcal{S}}} + \langle \boldsymbol{\Gamma}, \Delta^\gamma \rangle_{\widehat{\mathcal{S}}} + \langle D^\gamma, D^{\boldsymbol{\Gamma}} \rangle_{\widehat{\mathcal{S}}} + \langle \Delta^\gamma, \Delta^{\boldsymbol{\Gamma}} \rangle_{\widehat{\mathcal{S}}} \\
&+ \langle D^\gamma, \Delta^{\boldsymbol{\Gamma}} \rangle_{\widehat{\mathcal{S}}} + \langle D^{\boldsymbol{\Gamma}}, \Delta^\gamma \rangle_{\widehat{\mathcal{S}}} - \left( \sum_{j \in \mathcal{S}^* \setminus \widehat{\mathcal{S}}} \gamma_j \boldsymbol{\Gamma}_j - \sum_{j \in \widehat{\mathcal{S}} \setminus \mathcal{S}^*} \gamma_j \boldsymbol{\Gamma}_j \right).
\end{aligned}
\tag{111}
$$

Recall that $\widehat{\mathbf{v}}^{[j]} = \mathbf{W}^\intercal \widehat{\mathbf{U}}_{\cdot j}$, then we have the following expression

$$
\langle \gamma, D^\gamma \rangle_{\widehat{\mathcal{S}}} = \frac{1}{n} \sum_{j \in \widehat{\mathcal{S}}} \gamma_j (\widehat{\mathbf{v}}^{[j]})^\intercal \boldsymbol{\Pi}_{\cdot 2},
$$

and

$$\langle \boldsymbol{\gamma}, D^{\boldsymbol{\Gamma}} \rangle_{\widehat{\mathcal{S}}} + \langle \boldsymbol{\Gamma}, D^{\boldsymbol{\gamma}} \rangle_{\widehat{\mathcal{S}}} = \frac{1}{n} \sum_{j \in \widehat{\mathcal{S}}} (\widehat{\mathbf{v}}^{[j]})^{\mathsf{T}} \left( \gamma_j \boldsymbol{\Pi}_{\cdot 1} + \boldsymbol{\Gamma}_j \boldsymbol{\Pi}_{\cdot 2} \right).$$

Note that $\widehat{\mathcal{S}}$ is correlated with the error $\boldsymbol{\Pi}_{\cdot 1}$ and $\boldsymbol{\Pi}_{\cdot 2}$. However, we can compare $\widehat{\mathcal{S}}$ with the true support $\mathcal{S}^*$,

$$\begin{aligned}
\langle \boldsymbol{\gamma}, D^{\boldsymbol{\gamma}} \rangle_{\widehat{\mathcal{S}}} - \langle \boldsymbol{\gamma}, D^{\boldsymbol{\gamma}} \rangle_{\mathcal{S}^*} &= \frac{1}{n} \sum_{j \in \widehat{\mathcal{S}}} \gamma_j (\widehat{\mathbf{v}}^{[j]})^{\mathsf{T}} \boldsymbol{\Pi}_{\cdot 2} - \frac{1}{n} \sum_{j \in \mathcal{S}^*} \gamma_j (\widehat{\mathbf{v}}^{[j]})^{\mathsf{T}} \boldsymbol{\Pi}_{\cdot 2} \\
&= \frac{1}{n} \sum_{j \in \widehat{\mathcal{S}} \setminus \mathcal{S}^*} \gamma_j (\widehat{\mathbf{v}}^{[j]})^{\mathsf{T}} \boldsymbol{\Pi}_{\cdot 2} - \frac{1}{n} \sum_{j \in \mathcal{S}^* \setminus \widehat{\mathcal{S}}} \gamma_j (\widehat{\mathbf{v}}^{[j]})^{\mathsf{T}} \boldsymbol{\Pi}_{\cdot 2},
\end{aligned} \tag{112}$$

and

$$\begin{aligned}
&\left( \langle \boldsymbol{\gamma}, D^{\boldsymbol{\Gamma}} \rangle_{\widehat{\mathcal{S}}} + \langle \boldsymbol{\Gamma}, D^{\boldsymbol{\gamma}} \rangle_{\widehat{\mathcal{S}}} \right) - \left( \langle \boldsymbol{\gamma}, D^{\boldsymbol{\Gamma}} \rangle_{\mathcal{S}^*} + \langle \boldsymbol{\Gamma}, D^{\boldsymbol{\gamma}} \rangle_{\mathcal{S}^*} \right) \\
=& \frac{1}{n} \sum_{j \in \widehat{\mathcal{S}}} (\widehat{\mathbf{v}}^{[j]})^{\mathsf{T}} \left( \gamma_j \boldsymbol{\Pi}_{\cdot 1} + \boldsymbol{\Gamma}_j \boldsymbol{\Pi}_{\cdot 2} \right) - \frac{1}{n} \sum_{j \in \mathcal{S}^*} (\widehat{\mathbf{v}}^{[j]})^{\mathsf{T}} \left( \gamma_j \boldsymbol{\Pi}_{\cdot 1} + \boldsymbol{\Gamma}_j \boldsymbol{\Pi}_{\cdot 2} \right) \\
=& \frac{1}{n} \sum_{j \in \widehat{\mathcal{S}} \setminus \mathcal{S}^*} (\widehat{\mathbf{v}}^{[j]})^{\mathsf{T}} \left( \gamma_j \boldsymbol{\Pi}_{\cdot 1} + \boldsymbol{\Gamma}_j \boldsymbol{\Pi}_{\cdot 2} \right) - \frac{1}{n} \sum_{j \in \mathcal{S}^* \setminus \widehat{\mathcal{S}}} (\widehat{\mathbf{v}}^{[j]})^{\mathsf{T}} \left( \gamma_j \boldsymbol{\Pi}_{\cdot 1} + \boldsymbol{\Gamma}_j \boldsymbol{\Pi}_{\cdot 2} \right).
\end{aligned} \tag{113}$$

Hence, the residual terms are

$$\begin{aligned}
R^{\boldsymbol{\gamma}} =& \sqrt{n} \left( 2 \langle \boldsymbol{\gamma}, \Delta^{\boldsymbol{\gamma}} \rangle_{\widehat{\mathcal{S}}} + \langle D^{\boldsymbol{\gamma}}, D^{\boldsymbol{\gamma}} \rangle_{\widehat{\mathcal{S}}} + \langle \Delta^{\boldsymbol{\gamma}}, \Delta^{\boldsymbol{\gamma}} \rangle_{\widehat{\mathcal{S}}} + 2 \langle D^{\boldsymbol{\gamma}}, \Delta^{\boldsymbol{\gamma}} \rangle_{\widehat{\mathcal{S}}} \right) \\
&+ \sqrt{n} \left( \frac{1}{n} \sum_{j \in \widehat{\mathcal{S}} \setminus \mathcal{S}^*} \gamma_j (\widehat{\mathbf{v}}^{[j]})^{\mathsf{T}} \boldsymbol{\Pi}_{\cdot 2} - \frac{1}{n} \sum_{j \in \mathcal{S}^* \setminus \widehat{\mathcal{S}}} \gamma_j (\widehat{\mathbf{v}}^{[j]})^{\mathsf{T}} \boldsymbol{\Pi}_{\cdot 2} \right) - \sqrt{n} \left( \sum_{j \in \mathcal{S}^* \setminus \widehat{\mathcal{S}}} \gamma_j^2 - \sum_{j \in \widehat{\mathcal{S}} \setminus \mathcal{S}^*} \gamma_j^2 \right),
\end{aligned} \tag{114}$$

and

$$\begin{aligned}
R^{\text{inter}} =& \sqrt{n} \left( \langle \boldsymbol{\gamma}, \Delta^{\boldsymbol{\Gamma}} \rangle_{\widehat{\mathcal{S}}} + \langle \boldsymbol{\Gamma}, \Delta^{\boldsymbol{\gamma}} \rangle_{\widehat{\mathcal{S}}} + \langle D^{\boldsymbol{\gamma}}, D^{\boldsymbol{\Gamma}} \rangle_{\widehat{\mathcal{S}}} + \langle \Delta^{\boldsymbol{\gamma}}, \Delta^{\boldsymbol{\Gamma}} \rangle_{\widehat{\mathcal{S}}} + \langle D^{\boldsymbol{\gamma}}, \Delta^{\boldsymbol{\Gamma}} \rangle_{\widehat{\mathcal{S}}} + \langle D^{\boldsymbol{\Gamma}}, \Delta^{\boldsymbol{\gamma}} \rangle_{\widehat{\mathcal{S}}} \right) \\
&+ \sqrt{n} \left( \frac{1}{n} \sum_{j \in \widehat{\mathcal{S}} \setminus \mathcal{S}^*} (\widehat{\mathbf{v}}^{[j]})^{\mathsf{T}} \left( \gamma_j \boldsymbol{\Pi}_{\cdot 1} + \boldsymbol{\Gamma}_j \boldsymbol{\Pi}_{\cdot 2} \right) - \frac{1}{n} \sum_{j \in \mathcal{S}^* \setminus \widehat{\mathcal{S}}} (\widehat{\mathbf{v}}^{[j]})^{\mathsf{T}} \left( \gamma_j \boldsymbol{\Pi}_{\cdot 1} + \boldsymbol{\Gamma}_j \boldsymbol{\Pi}_{\cdot 2} \right) \right) \\
&- \sqrt{n} \left( \sum_{j \in \mathcal{S}^* \setminus \widehat{\mathcal{S}}} \gamma_j \boldsymbol{\Gamma}_j - \sum_{j \in \widehat{\mathcal{S}} \setminus \mathcal{S}^*} \gamma_j \boldsymbol{\Gamma}_j \right).
\end{aligned} \tag{115}$$

Define $\mathcal{S}_0^* = \left\{ j : |\gamma_j| > \sqrt{2.05 \log p_z} \sqrt{\text{Var}\left( D_j^{\boldsymbol{\gamma}} \right)} \right\}$ to be the set of strong signals, on the event $A_2$, we have

$$\mathcal{S}_0^* \subset \widehat{\mathcal{S}} \subset \mathcal{S}^*, \quad \text{and} \quad \left| \widehat{\mathcal{S}} \right| \leq s_{z1}. \tag{116}$$

On the event $A_3$, we have

$$\max\left\{\left\|D^{\mathbf{\Gamma}}\right\|_\infty, \left\|D^\gamma\right\|_\infty\right\} \leq \left(1 + 12\sqrt{\frac{\log p}{n}}\right) M_1 \sqrt{\frac{2.05\log p_z}{n}} \max\{\sigma_1, \sigma_2\}. \quad (117)$$

On the event $S \cap G \cap A$,

$$\max\left\{\|\Delta^\gamma\|_\infty, \|\Delta^{\mathbf{\Gamma}}\|_\infty\right\} \leq Cs\frac{\log p}{n}\max\{\sigma_1, \sigma_2\}. \quad (118)$$

Combing (116), (117) and (118), we have on the event $S \cap G \cap A$,

$$\max\left\{\langle D^\gamma, D^\gamma\rangle_{\widehat{S}}, \langle D^{\mathbf{\Gamma}}, D^{\mathbf{\Gamma}}\rangle_{\widehat{S}}\right\} \leq Cs_{z1}\frac{\log p_z}{n}, \quad (119)$$

$$\max\left\{\langle \Delta^\gamma, \Delta^\gamma\rangle_{\widehat{S}}, \langle \Delta^{\mathbf{\Gamma}}, \Delta^{\mathbf{\Gamma}}\rangle_{\widehat{S}}\right\} \leq Cs_{z1}\left(s\frac{\log p}{n}\right)^2. \quad (120)$$

Note that

$$\left|\langle D^\gamma, \Delta^\gamma\rangle_{\widehat{S}}\right| \leq \sqrt{\langle D^\gamma, D^\gamma\rangle_{\widehat{S}}\langle\Delta^\gamma, \Delta^\gamma\rangle_{\widehat{S}}} \leq \frac{1}{2}\left(\langle D^\gamma, D^\gamma\rangle_{\widehat{S}} + \langle\Delta^\gamma, \Delta^\gamma\rangle_{\widehat{S}}\right).$$

Hence, we have

$$\left|\langle D^\gamma, D^\gamma\rangle_{\widehat{S}} + \langle\Delta^\gamma, \Delta^\gamma\rangle_{\widehat{S}} + 2\langle D^\gamma, \Delta^\gamma\rangle_{\widehat{S}}\right| \leq Cs_{z1}\frac{\log p_z}{n} + Cs_{z1}\left(s\frac{\log p}{n}\right)^2, \quad (121)$$

and

$$\left|\langle D^\gamma, D^{\mathbf{\Gamma}}\rangle_{\widehat{S}} + \langle\Delta^\gamma, \Delta^{\mathbf{\Gamma}}\rangle_{\widehat{S}} + \langle D^\gamma, \Delta^{\mathbf{\Gamma}}\rangle_{\widehat{S}} + \langle D^{\mathbf{\Gamma}}, \Delta^\gamma\rangle_{\widehat{S}}\right| \leq Cs_{z1}\frac{\log p_z}{n} + Cs_{z1}\left(s\frac{\log p}{n}\right)^2. \quad (122)$$

We also have the following control

$$2\left|\langle\gamma, \Delta^\gamma\rangle_{\widehat{S}}\right| \leq \|\gamma\|_2\sqrt{\langle\Delta^\gamma, \Delta^\gamma\rangle_{\widehat{S}}} \leq C\|\gamma\|_2\sqrt{s_{z1}}s\frac{\log p}{n}, \quad (123)$$

and

$$\left|\langle\gamma, \Delta^{\mathbf{\Gamma}}\rangle_{\widehat{S}} + \langle\mathbf{\Gamma}, \Delta^\gamma\rangle_{\widehat{S}}\right| \leq C\left(\|\gamma\|_2 + \|\mathbf{\Gamma}\|_2\right)\sqrt{s_{z1}}s\frac{\log p}{n}. \quad (124)$$

On the event $S \cap G \cap A$, we have $\widehat{S}\backslash S^* = \varnothing$ and hence $\frac{1}{n}\sum_{j\in\widehat{S}\backslash S^*}\gamma_j(\widehat{\mathbf{v}}^{[j]})^\mathsf{T}\mathbf{\Pi}_{\cdot 2} = 0$, $\sum_{j\in\widehat{S}\backslash S^*}\gamma_j\mathbf{\Gamma}_j = 0$ and $\sum_{j\in\widehat{S}\backslash S^*}\gamma_j^2 = 0$; On the event $S \cap G \cap A$, we also have

$$\left|\frac{1}{n}\sum_{j\in S^*\backslash\widehat{S}}\gamma_j(\widehat{\mathbf{v}}^{[j]})^\mathsf{T}\mathbf{\Pi}_{\cdot 2}\right| \leq \frac{1}{n}s_{z1}\max_{j\in S^*\backslash\widehat{S}}|\gamma_j|\left|(\widehat{\mathbf{v}}^{[j]})^\mathsf{T}\mathbf{\Pi}_{\cdot 2}\right|$$

$$\leq s_{z1}\sqrt{2.05\log p_z}\sqrt{\text{Var}\left(D_j^\gamma\right)}\left(1 + 12\sqrt{\frac{\log p}{n}}\right)M_1\sqrt{\frac{2.05\log p_z}{n}}\sigma_2 \leq \frac{s_{z1}\log p_z}{n}. \quad (125)$$

On the event $S \cap G \cap A$, we get

$$\left| \sum_{j \in \mathcal{S}^* \setminus \widehat{\mathcal{S}}} \boldsymbol{\gamma}_j^2 \right| \leq s_{z1} \frac{\log p_z}{n}. \tag{126}$$

By (114), (121), (123), (125) and (126), we establish that on the event $S \cap G \cap A$,

$$|R^\gamma| \leq C s_{z1} \frac{\log p_z}{\sqrt{n}} + C \|\boldsymbol{\gamma}\|_2 \sqrt{s_{z1}} s \frac{\log p}{\sqrt{n}}. \tag{127}$$

Similarly, we can establish that and

$$\left| R^{\text{inter}} \right| \leq C s_{z1} \frac{\log p_z}{\sqrt{n}} + C \left( \|\boldsymbol{\gamma}\|_2 + \|\boldsymbol{\Gamma}\|_2 \right) \sqrt{s_{z1}} s \frac{\log p}{\sqrt{n}}. \tag{128}$$

We can establish (71) and (74) by taking $\boldsymbol{\Gamma}_j = \beta \boldsymbol{\gamma}_j$. Note that

$$\frac{2\sqrt{\log p}}{\sqrt{n}} \left\| \sum_{j \in \mathcal{S}^*} \gamma_j \widehat{\mathbf{v}}^{[j]} \right\|_2 \sqrt{\boldsymbol{\Theta}_{22}} = 2\sqrt{\frac{\log p}{n}} \sqrt{\frac{\boldsymbol{\Theta}_{22}}{\boldsymbol{\Theta}_{11} + \beta^2 \boldsymbol{\Theta}_{22} - 2\beta \boldsymbol{\Theta}_{12}}} \sqrt{\text{Var}_H} \|\boldsymbol{\gamma}\|_2^2,$$

$$\frac{\sqrt{\log p}}{\sqrt{n}} \left\| \sum_{j \in \mathcal{S}^*} \gamma_j \widehat{\mathbf{v}}^{[j]} \right\|_2 \sqrt{\boldsymbol{\Theta}_{11} + \beta^2 \boldsymbol{\Theta}_{22} + 2\beta \boldsymbol{\Theta}_{12}} = \sqrt{\frac{\log p}{n}} \sqrt{\frac{\boldsymbol{\Theta}_{11} + \beta^2 \boldsymbol{\Theta}_{22} + 2\beta \boldsymbol{\Theta}_{12}}{\boldsymbol{\Theta}_{11} + \beta^2 \boldsymbol{\Theta}_{22} - 2\beta \boldsymbol{\Theta}_{12}}} \sqrt{\text{Var}_H} \|\boldsymbol{\gamma}\|_2^2. \tag{129}$$

By the definition of $A_4$ and $A_5$ in (61) and Lemma 5, we establish

$$\max \left\{ \left| \frac{2}{n} \sum_{j \in \mathcal{S}^*} \gamma_j (\widehat{\mathbf{v}}^{[j]})^{\mathsf{T}} \boldsymbol{\Pi}_{\cdot 2} \right|, \left| \frac{1}{n} \sum_{j \in \mathcal{S}^*} \gamma_j (\widehat{\mathbf{v}}^{[j]})^{\mathsf{T}} \left( \boldsymbol{\Pi}_{\cdot 1} + \beta \boldsymbol{\Pi}_{\cdot 2} \right) \right| \right\} \leq C s_{z1} \sqrt{\frac{\log p}{n}} \|\boldsymbol{\gamma}\|_2 \tag{130}$$

Combined with (74), we establish (75).

### D.6.  Proof of Lemma 2

Since $\min_{j \in \mathcal{S}^*} |\boldsymbol{\gamma}_j| \geq \delta_{\min} \gg \sqrt{\log p / n}$, on the event $A \cap S \cap G$, we have $\widehat{\mathcal{S}} = \mathcal{S}^*$. Without loss of generality, we assume $\widehat{\mathcal{S}} = \{1, 2, \cdots, q\}$ and analyze $\widehat{\mathcal{V}}^{[1]}$. We start with the analysis of $\widehat{\beta}^{[1]}$,

$$\sqrt{n} \left( \widehat{\beta}^{[1]} - \left( \beta + \frac{\pi_1}{\gamma_1} \right) \right) = T^{\beta, 1} + \Delta^{\beta, 1}, \tag{131}$$

where

$$T^{\beta, 1} = \frac{1}{\sqrt{n} \gamma_1} (\widehat{\mathbf{v}}^{[1]})^{\mathsf{T}} \left( \boldsymbol{\Pi}_{\cdot 1} - \left( \beta + \frac{\pi_1}{\gamma_1} \right) \boldsymbol{\Pi}_{\cdot 2} \right) \quad \text{and} \quad \Delta^{\beta, 1} = R_1 + R_2, \tag{132}$$

with

$$R_1 = \frac{\sqrt{n}}{\gamma_1} \left( \Delta_1^{\boldsymbol{\Gamma}} - \left( \beta + \frac{\pi_1}{\gamma_1} \right) \Delta_1^{\gamma} \right) \text{ and } R_2 = \frac{-\left( D_1^{\gamma} + \Delta_1^{\gamma} \right)}{\gamma_1 + \left( D_1^{\gamma} + \Delta_1^{\gamma} \right)} \left( T^{\beta, 1} + R_1 \right). \tag{133}$$

We focus on the case $j = 1$ and analyze the following estimator,

$$\widehat{\boldsymbol{\pi}}^{[1]} = \widehat{\boldsymbol{\Gamma}} - \widehat{\beta}^{[1]}\widehat{\boldsymbol{\gamma}}. \tag{134}$$

Note that

$$\widehat{\boldsymbol{\pi}}_k^{[1]} - \boldsymbol{\pi}_k = -\frac{\boldsymbol{\pi}_1}{\gamma_1}\gamma_k + \left(\widehat{\boldsymbol{\Gamma}}_k - \boldsymbol{\Gamma}_k\right) - \left(\beta + \frac{\boldsymbol{\pi}_1}{\gamma_1}\right)\left(\widehat{\gamma}_k - \gamma_k\right) - \gamma_k\left(\widehat{\beta}^{[1]} - \left(\beta + \frac{\boldsymbol{\pi}_1}{\gamma_1}\right)\right)$$
$$- \left(\widehat{\beta}^{[1]} - \left(\beta + \frac{\boldsymbol{\pi}_1}{\gamma_1}\right)\right)\left(\widehat{\gamma}_k - \gamma_k\right). \tag{135}$$

By (93) and (94), we have

$$\sqrt{n}\left(\widehat{\boldsymbol{\Gamma}}_k - \boldsymbol{\Gamma}_k\right) = \frac{1}{\sqrt{n}}(\widehat{\mathbf{v}}^{[k]})^{\intercal}\boldsymbol{\Pi}_{\cdot 1} + \sqrt{n}\Delta_k^{\boldsymbol{\Gamma}}. \tag{136}$$

By (95) and (96), we have

$$\sqrt{n}\left(\widehat{\gamma}_k - \gamma_k\right) = \frac{1}{\sqrt{n}}(\widehat{\mathbf{v}}^{[k]})^{\intercal}\boldsymbol{\Pi}_{\cdot 2} + \sqrt{n}\Delta_k^{\gamma}. \tag{137}$$

By plugging (132), (136) and (137) into (135), we have the following decomposition of $\widehat{\boldsymbol{\pi}}_k^{[1]} - \boldsymbol{\pi}_k$

$$\sqrt{n}\left(\widehat{\boldsymbol{\pi}}_k^{[1]} - \boldsymbol{\pi}_k\right) = -\sqrt{n}\frac{\boldsymbol{\pi}_1}{\gamma_1}\gamma_k + T^{\boldsymbol{\pi}_k} + \Delta^{\boldsymbol{\pi}_k}, \tag{138}$$

where

$$T^{\boldsymbol{\pi}_k} = \frac{1}{\sqrt{n}}\left((\widehat{\mathbf{v}}^{[k]})^{\intercal} - \frac{\gamma_k}{\gamma_1}(\widehat{\mathbf{v}}^{[1]})^{\intercal}\right)\left(\boldsymbol{\Pi}_{\cdot 1} - \left(\beta + \frac{\boldsymbol{\pi}_1}{\gamma_1}\right)\boldsymbol{\Pi}_{\cdot 2}\right),$$

and

$$\Delta^{\boldsymbol{\pi}_k} = \sqrt{n}\left(\Delta_k^{\boldsymbol{\Gamma}} - \left(\beta + \frac{\boldsymbol{\pi}_1}{\gamma_1}\right)\Delta_k^{\gamma} - \gamma_k\Delta^{\beta,1}\right) - \left(T^{\beta,1} + \Delta^{\beta,1}\right)\left(\widehat{\gamma}_k - \gamma_k\right). \tag{139}$$

Define the events for $i \in \mathcal{S}^*$,

$$F^i = \left\{\max_{k \in \mathcal{S}^*, k \neq i}|T^{\boldsymbol{\pi}_k}| \leq 2.02\sqrt{\log p_z}\sqrt{\boldsymbol{\Theta}_{11} + \left(\beta + \frac{\boldsymbol{\pi}_i}{\gamma_i}\right)^2\boldsymbol{\Theta}_{22} - 2\left(\beta + \frac{\boldsymbol{\pi}_i}{\gamma_i}\right)\boldsymbol{\Theta}_{12}}\frac{\left\|\widehat{\mathbf{v}}^{[k]} - \frac{\gamma_k}{\gamma_i}\widehat{v}^{[i]}\right\|_2}{\sqrt{n}}.\right\}$$

Then for $F = \cap_{i \in \mathcal{S}^*}F^i$, we have

$$\mathbf{P}\left(F\right) \geq 1 - Cs_{z1}^2p_z^{-2.04} \geq 1 - cp^{-c}. \tag{140}$$

The proof of Lemma 2 relies on the following lemmas. The following lemma provides upper bound and lower bound for the variance term and the proof of the following lemma can be found in Section D.9.

LEMMA 11. *On the event $A \cap S \cap G$, we have*

$$
\sqrt{\boldsymbol{\Theta}_{11} + \left(\beta + \left|\frac{\pi_1}{\gamma_1}\right|\right)^2 \boldsymbol{\Theta}_{22} - 2\left(\beta + \frac{\pi_1}{\gamma_1}\right)\boldsymbol{\Theta}_{12}} \frac{\left\|\widehat{\mathbf{v}}^{[k]} - \frac{\gamma_k}{\gamma_1}\widehat{\mathbf{v}}^{[1]}\right\|_2}{\sqrt{n}}
$$
$$
\leq 1.1\sqrt{M_1 M_2}\left(1 + \left|\frac{\gamma_k}{\gamma_1}\right|\right)\sqrt{1 + \left(\beta + \frac{\pi_1}{\gamma_1}\right)^2},
$$
(141)

*and*

$$
\sqrt{\boldsymbol{\Theta}_{11} + \left(\beta + \left|\frac{\pi_1}{\gamma_1}\right|\right)^2 \boldsymbol{\Theta}_{22} - 2\left(\beta + \frac{\pi_1}{\gamma_1}\right)\boldsymbol{\Theta}_{12}} \frac{\left\|\widehat{\mathbf{v}}^{[k]} - \frac{\gamma_k}{\gamma_1}\widehat{\mathbf{v}}^{[1]}\right\|_2}{\sqrt{n}}
$$
$$
\geq 0.45\sqrt{\frac{M_1}{M_2}}\left(1 + \left|\frac{\gamma_k}{\gamma_1}\right|\right)\sqrt{1 + \left(\beta + \frac{\pi_1}{\gamma_1}\right)^2}.
$$
(142)

LEMMA 12. *On the event $A \cap S \cap G \cap F^1$, for large $n$, we have*

$$
0.995 \leq \frac{\sqrt{\widehat{\boldsymbol{\Theta}}_{11} + \left(\widehat{\beta}^{[1]}\right)^2 \widehat{\boldsymbol{\Theta}}_{22} - 2\widehat{\beta}^{[1]}\widehat{\boldsymbol{\Theta}}_{12}}\left\|\widehat{\mathbf{v}}^{[k]} - \frac{\widehat{\gamma}_k}{\widehat{\gamma}_1}\widehat{\mathbf{v}}^{[1]}\right\|_2}{\sqrt{\boldsymbol{\Theta}_{11} + \left(\beta + \frac{\pi_1}{\gamma_1}\right)^2 \boldsymbol{\Theta}_{22} - 2\left(\beta + \frac{\pi_1}{\gamma_1}\right)\boldsymbol{\Theta}_{12}}\left\|\widehat{\mathbf{v}}^{[k]} - \frac{\gamma_k}{\gamma_1}\widehat{\mathbf{v}}^{[1]}\right\|_2} \leq 1.005.
$$
(143)

*On the event $A \cap S \cap G \cap F^1$, we have*

$$
\max_{k \in \mathcal{S}^*} \frac{1}{\sqrt{n}}|T^{\pi_k}| \leq 2.02\sqrt{\frac{\log p_z}{n}}\sqrt{\boldsymbol{\Theta}_{11} + \left(\beta + \frac{\pi_1}{\gamma_1}\right)^2 \boldsymbol{\Theta}_{22} - 2\left(\beta + \frac{\pi_1}{\gamma_1}\right)\boldsymbol{\Theta}_{12}} \frac{\left\|\widehat{\mathbf{v}}^{[k]} - \frac{\gamma_k}{\gamma_1}\widehat{\mathbf{v}}^{[1]}\right\|_2}{\sqrt{n}};
$$
(144)

*and*

$$
\max_{k \in \mathcal{S}^*} \frac{1}{\sqrt{n}}|\Delta^{\pi_k}| \leq \frac{1}{300}\sqrt{\frac{\log p_z}{n}}\sqrt{\boldsymbol{\Theta}_{11} + \left(\beta + \frac{\pi_1}{\gamma_1}\right)^2 \boldsymbol{\Theta}_{22} - 2\left(\beta + \frac{\pi_1}{\gamma_1}\right)\boldsymbol{\Theta}_{12}} \frac{\left\|\widehat{\mathbf{v}}^{[k]} - \frac{\gamma_k}{\gamma_1}\widehat{\mathbf{v}}^{[1]}\right\|_2}{\sqrt{n}}.
$$
(145)

To analyze the first column of the voting matrix, that is $\mathrm{VM}_{\cdot,1}$, we first establish the following observation, on the event $A \cap S \cap G \cap F$,

$$
\text{If } \frac{\pi_k}{\gamma_k} = \frac{\pi_1}{\gamma_1}, \text{ then } k \in \widehat{\mathcal{V}}^{[1]}.
$$
(146)

$$
\text{If } \left|\frac{\pi_k}{\gamma_k} - \frac{\pi_1}{\gamma_1}\right| \geq C_*(1/\delta_{\min})\sqrt{\log p_z/n}, \text{ then } k \notin \widehat{\mathcal{V}}^{[1]},
$$
(147)

with

$$
C_* = 12\left(1 + \max_{k \in \mathcal{S}^*}\left|\frac{\boldsymbol{\Gamma}_k}{\gamma_k}\right|\right)\sqrt{M_1/M_2}.
$$

<u>Proof of (147)</u> In this case, $\left|\pi_k - \frac{\pi_1}{\gamma_1}\gamma_k\right| \geq C_*|\gamma_k|(1/\delta_{\min})\sqrt{\log p_z/n}$. We rewrite (138) as

$$\widehat{\pi}_k^{[1]} - \left(\pi_k - \frac{\pi_1}{\gamma_1}\gamma_k\right) = \frac{1}{\sqrt{n}}\left(T^{\pi_k} + \Delta^{\pi_k}\right). \tag{148}$$

It is sufficient to show that $k \in \mathcal{S}^*$,

$$\left|\pi_k - \frac{\pi_1}{\gamma_1}\gamma_k + \frac{1}{\sqrt{n}}\left(T^{\pi_k} + \Delta^{\pi_k}\right)\right| \geq 2.05\sqrt{\widehat{\Theta}_{11} + \left(\widehat{\beta}^{[1]}\right)^2\widehat{\Theta}_{22} - 2\widehat{\beta}^{[1]}\widehat{\Theta}_{12}}\sqrt{\frac{\log p_z}{n}}\frac{\left\|\widehat{\mathbf{v}}^{[k]} - \frac{\widehat{\gamma}_k}{\widehat{\gamma}_1}\widehat{\mathbf{v}}^{[1]}\right\|_2}{\sqrt{n}}.$$

The above equation can be established by the following two results,

$$\max_{k \in \mathcal{S}^*}\frac{1}{\sqrt{n}}|T^{\pi_k} + \Delta^{\pi_k}| \leq 2.05\sqrt{\widehat{\Theta}_{11} + \left(\widehat{\beta}^{[1]}\right)^2\widehat{\Theta}_{22} - 2\widehat{\beta}^{[1]}\widehat{\Theta}_{12}}\sqrt{\frac{\log p_z}{n}}\frac{\left\|\widehat{\mathbf{v}}^{[k]} - \frac{\widehat{\gamma}_k}{\widehat{\gamma}_1}\widehat{\mathbf{v}}^{[1]}\right\|_2}{\sqrt{n}}, \tag{149}$$

and

$$\left|\pi_k - \frac{\pi_1}{\gamma_1}\gamma_k\right| \geq 4.1\sqrt{\widehat{\Theta}_{11} + \left(\widehat{\beta}^{[1]}\right)^2\widehat{\Theta}_{22} - 2\widehat{\beta}^{[1]}\widehat{\Theta}_{12}}\sqrt{\frac{\log p_z}{n}}\frac{\left\|\widehat{\mathbf{v}}^{[k]} - \frac{\widehat{\gamma}_k}{\widehat{\gamma}_1}\widehat{\mathbf{v}}^{[1]}\right\|_2}{\sqrt{n}}. \tag{150}$$

By (143), (144) and (145), we establish (149). To establish (150), it is sufficient to show that

$$\left|4.1\sqrt{\widehat{\Theta}_{11} + \left(\widehat{\beta}^{[1]}\right)^2\widehat{\Theta}_{22} - 2\widehat{\beta}^{[1]}\widehat{\Theta}_{12}}\sqrt{\frac{\log p_z}{n}}\frac{\left\|\widehat{\mathbf{v}}^{[k]} - \frac{\widehat{\gamma}_k}{\widehat{\gamma}_1}\widehat{\mathbf{v}}^{[1]}\right\|_2}{\sqrt{n}}\right| \leq C_*|\gamma_k|(1/\delta_{\min})\sqrt{\log p_z/n} \tag{151}$$

By (141) and (143), we have

$$\begin{aligned}
&4.1\sqrt{\widehat{\Theta}_{11} + \left(\widehat{\beta}^{[1]}\right)^2\widehat{\Theta}_{22} - 2\widehat{\beta}^{[1]}\widehat{\Theta}_{12}}\sqrt{\frac{\log p_z}{n}}\frac{\left\|\widehat{\mathbf{v}}^{[k]} - \frac{\widehat{\gamma}_k}{\widehat{\gamma}_1}\widehat{\mathbf{v}}^{[1]}\right\|_2}{\sqrt{n}} \\
\leq &4.1 \times 1.005 \times 1.1\sqrt{M_1 M_2}\left(1 + \left|\frac{\gamma_k}{\gamma_1}\right|\right)\sqrt{1 + \left(\beta + \frac{\pi_1}{\gamma_1}\right)^2}\sqrt{\frac{\log p_z}{n}} \\
\leq &|\gamma_k| 4.1 \times 1.005 \times 1.1\sqrt{M_1 M_2}\left(\left|\frac{1}{\gamma_k}\right| + \left|\frac{1}{\gamma_1}\right|\right)\left(1 + \left|\beta + \frac{\pi_1}{\gamma_1}\right|\right)\sqrt{\frac{\log p_z}{n}}.
\end{aligned} \tag{152}$$

The last term can be further upper bounded by

$$|\gamma_k|\frac{1}{\delta_{\min}}8.2 \times 1.005 \times 1.1\sqrt{M_1 M_2}\left(1 + \left|\beta + \frac{\pi_1}{\gamma_1}\right|\right)\sqrt{\frac{\log p_z}{n}} \leq \frac{C_*}{\delta_{\min}}|\gamma_k|\sqrt{\frac{\log p_z}{n}} \tag{153}$$

where the inequality follows from the definition of $C_*$. By (152) and (153), we conclude (151) and hence (150).

<u>Proof of (146)</u> In this case, $\left|\pi_k - \frac{\pi_1}{\gamma_1}\gamma_k\right| = 0$. For $k \in \mathcal{S}^*$, (138) can be re-expressed as

$$\sqrt{n}\left(\widehat{\pi}_k^{[1]} - 0\right) = T^{\pi_k} + \Delta^{\pi_k}. \tag{154}$$

By (143), (144) and (145), on the event $A \cap S \cap G \cap F^1$,

$$\max_{k \in \mathcal{V}^*} \frac{1}{\sqrt{n}} |T^{\pi_k} + \Delta^{\pi_k}| \leq 2.05 \sqrt{\widehat{\Theta}_{11} + \left(\widehat{\beta}^{[1]}\right)^2 \widehat{\Theta}_{22} - 2\widehat{\beta}^{[1]}\widehat{\Theta}_{12}} \sqrt{\frac{\log p_z}{n}} \frac{\left\| \widehat{\mathbf{v}}^{[k]} - \frac{\widehat{\gamma}_k}{\widehat{\gamma}_1} \widehat{\mathbf{v}}^{[1]} \right\|_2}{\sqrt{n}}.$$
(155)

Hence, we establish (146).

The results (146) and (147) can be generalized to (156) and (157) with the same proof. On the event $A \cap S \cap G \cap F$,

$$\text{If } \frac{\pi_k}{\gamma_k} = \frac{\pi_j}{\gamma_j}, \text{ then } k \in \widehat{\mathcal{V}}^{[j]}.$$
(156)

$$\text{If } \left| \frac{\pi_k}{\gamma_k} - \frac{\pi_j}{\gamma_j} \right| \geq C_*(1/\delta_{\min})\sqrt{\log p_z/n}, \text{ then } k \notin \widehat{\mathcal{V}}^{[j]}.$$
(157)

- If the assumption (IN1-M) holds, in this case, (IN1-P) automatically holds, then by (156) and (157), we have

$$\min_{k \in \mathcal{V}^*} \|\text{VM}_k\|_0 > \frac{q}{2} > \max_{k \in \mathcal{S}^* \setminus \mathcal{V}^*} \|\text{VM}_k\|_0$$
(158)

- If the assumption (IN1-M) does not hold, but (IN1-P) holds, then then by (156) and (157), we have

$$\min_{k \in \mathcal{V}^*} \|\text{VM}_k\|_0 = \max_{k \in \mathcal{S}^*} \|\text{VM}_k\|_0 > \max_{k \in \mathcal{S}^* \setminus \mathcal{V}^*} \|\text{VM}_k\|_0$$
(159)

Hence, on the event $A \cap S \cap G \cap F$, we have $\widehat{\mathcal{V}} = \mathcal{V}^*$.

### D.7.  Proof of Lemma 7

The proof of this lemma follows from the following results. Under the regularity assumptions $(R1) - (R3)$, as $\sqrt{s_{z1}}s \log p/\sqrt{n} \to 0$, we have

$$\max_{1 \leq i,j \leq 2} |\widehat{\Theta}_{ij} - \Theta_{ij}| \xrightarrow{p} 0;$$
(160)

and

$$\frac{\widehat{\|\gamma\|_2^2}}{\|\gamma\|_2^2} \xrightarrow{p} 1 \quad \text{and} \quad \frac{\left\| \sum_{j \in \widehat{\mathcal{S}}} \widehat{\gamma}_j \widehat{\mathbf{v}}^{[j]} \right\|_2}{\left\| \sum_{j \in \mathcal{S}^*} \gamma_j \widehat{\mathbf{v}}^{[j]} \right\|_2} \xrightarrow{p} 1.$$
(161)

By (27) and (160), we establish that

$$\frac{\sqrt{\widehat{\Theta}_{11} + \widehat{\beta}^2\widehat{\Theta}_{22} - 2\widehat{\beta}\widehat{\Theta}_{12}}}{\sqrt{\Theta_{11} + (\beta)^2\Theta_{22} - 2\beta\Theta_{12}}} \xrightarrow{p} 1.$$

Combined with (161), we establish (82).

**Proof of** (160) A stronger version of this proposition has already been proved in Ren et al. [2013], where part of it was already established in Sun and Zhang [2012]. To be self-contained, we will provide the sketch of the proof in the following.

The difference between $\widehat{\boldsymbol{\Theta}} - \boldsymbol{\Theta}$ can be decomposed as,

$$\widehat{\boldsymbol{\Theta}} - \boldsymbol{\Theta} = \boldsymbol{\Theta}^{\mathrm{ora}} - \boldsymbol{\Theta} + \widehat{\boldsymbol{\Theta}} - \boldsymbol{\Theta}^{\mathrm{ora}}, \tag{162}$$

where $\boldsymbol{\Theta}_{11}^{\mathrm{ora}} = \frac{1}{n}\|Y - \mathbf{Z}\boldsymbol{\Gamma} - \mathbf{X}\boldsymbol{\Psi}\|_2^2$, $\boldsymbol{\Theta}_{22}^{\mathrm{ora}} = \frac{1}{n}\|D - \mathbf{Z}\boldsymbol{\gamma} - \mathbf{X}\boldsymbol{\psi}\|_2^2$ and $\boldsymbol{\Theta}_{12}^{\mathrm{ora}} = \frac{1}{n}(Y - \mathbf{Z}\boldsymbol{\Gamma} - \mathbf{X}\boldsymbol{\Psi})^{\mathsf{T}}(D - \mathbf{Z}\boldsymbol{\gamma} - \mathbf{X}\boldsymbol{\psi})$. In the following, we only provide the detailed analysis of $\widehat{\boldsymbol{\Theta}}_{12} - \boldsymbol{\Theta}_{12}^{\mathrm{ora}}$. The other differences can be established in a similar way and the difference between $\boldsymbol{\Theta}^{\mathrm{ora}} - \boldsymbol{\Theta}$ can be established by central limit theorem.

$$\widehat{\boldsymbol{\Theta}}_{12} - \boldsymbol{\Theta}_{12}^{\mathrm{ora}} = \frac{1}{n}\begin{pmatrix} \widetilde{\boldsymbol{\gamma}} - \boldsymbol{\gamma} \\ \widetilde{\boldsymbol{\psi}} - \boldsymbol{\psi} \end{pmatrix}^{\mathsf{T}} \mathbf{W}^{\mathsf{T}}\mathbf{W}\begin{pmatrix} \widetilde{\boldsymbol{\Gamma}} - \boldsymbol{\Gamma} \\ \widetilde{\boldsymbol{\Psi}} - \boldsymbol{\Psi} \end{pmatrix} + \frac{1}{n}\boldsymbol{\Pi}_{\cdot 2}^{\mathsf{T}}\mathbf{W}\begin{pmatrix} \widetilde{\boldsymbol{\Gamma}} - \boldsymbol{\Gamma} \\ \widetilde{\boldsymbol{\Psi}} - \boldsymbol{\Psi} \end{pmatrix} + \frac{1}{n}\boldsymbol{\Pi}_{\cdot 1}^{\mathsf{T}}\mathbf{W}\begin{pmatrix} \widetilde{\boldsymbol{\gamma}} - \boldsymbol{\gamma} \\ \widetilde{\boldsymbol{\psi}} - \boldsymbol{\psi} \end{pmatrix}. \tag{163}$$

By (163), we have

$$\begin{aligned} \left|\widehat{\boldsymbol{\Theta}}_{12} - \boldsymbol{\Theta}_{12}^{\mathrm{ora}}\right| \leq & \frac{1}{\sqrt{n}}\|\mathbf{W}\begin{pmatrix} \widetilde{\boldsymbol{\gamma}} - \boldsymbol{\gamma} \\ \widetilde{\boldsymbol{\psi}} - \boldsymbol{\psi} \end{pmatrix}\|_2 \frac{1}{\sqrt{n}}\|\mathbf{W}\begin{pmatrix} \widetilde{\boldsymbol{\Gamma}} - \boldsymbol{\Gamma} \\ \widetilde{\boldsymbol{\Psi}} - \boldsymbol{\Psi} \end{pmatrix}\|_2 + \frac{1}{n}\|\boldsymbol{\Pi}_{\cdot 2}^{\mathsf{T}}\mathbf{W}\|_\infty\|\begin{pmatrix} \widetilde{\boldsymbol{\Gamma}} - \boldsymbol{\Gamma} \\ \widetilde{\boldsymbol{\Psi}} - \boldsymbol{\Psi} \end{pmatrix}\|_1 \\ & + \frac{1}{n}\|\boldsymbol{\Pi}_{\cdot 1}^{\mathsf{T}}\mathbf{W}\|_\infty\|\begin{pmatrix} \widetilde{\boldsymbol{\gamma}} - \boldsymbol{\gamma} \\ \widetilde{\boldsymbol{\psi}} - \boldsymbol{\psi} \end{pmatrix}\|_1. \end{aligned} \tag{164}$$

The following of the proof follows from Lemma 9 and definition of event $G$.

**Proof of** (161) For a given $0 < \epsilon_0 < 1$, we have

$$\mathbf{P}\left(\left|\frac{\|\gamma\|_2^2}{\widehat{\|\gamma\|_2^2}} - 1\right| \geq \epsilon_0\right) \leq \mathbf{P}\left(\left|\frac{\widehat{\|\gamma\|_2^2} - \|\gamma\|_2^2}{\|\gamma\|_2^2}\right| \geq \frac{\epsilon_0}{1 - \epsilon_0}\right).$$

By Lemma 6, on the event $A \cap S \cap G$, we have

$$\left|\frac{\widehat{\|\gamma\|_2^2} - \|\gamma\|_2^2}{\|\gamma\|_2^2}\right| \leq C\frac{1}{\|\gamma\|_2^2}\left(s_{\mathrm{z}1}\frac{\log p_{\mathrm{z}}}{n} + C\|\gamma\|_2\sqrt{\frac{2s_{\mathrm{z}1}\log p_{\mathrm{z}}}{n}}\right).$$

Since $\|\gamma\|_2^2 \gg (s\log p/\sqrt{n})^2$, we obtain that

$$\left|\frac{\widehat{\|\gamma\|_2^2} - \|\gamma\|_2^2}{\|\gamma\|_2^2}\right| \leq \frac{\epsilon_0}{1 - \epsilon_0} \text{ and } \mathbf{P}\left(\left|\frac{\|\gamma\|_2^2}{\widehat{\|\gamma\|_2^2}} - 1\right| \geq \epsilon_0\right) \leq \mathbf{P}\left((A \cap S \cap G)^c\right).$$

Combined with Lemma 4, we establish the first convergence result of (161). On the event $S \cap G \cap A$, we have

$$\left|\frac{\left\|\sum_{j\in\widehat{\mathcal{S}}}\widehat{\gamma}_j\widehat{\mathbf{v}}^{[j]}\right\|_2}{\left\|\sum_{j\in\mathcal{S}^*}\gamma_j\widehat{\mathbf{v}}^{[j]}\right\|_2} - 1\right| \leq \frac{\sum_{j\in\widehat{\mathcal{S}}}|\widehat{\gamma}_j - \gamma_j|\frac{\|\widehat{\mathbf{v}}^{[j]}\|}{\sqrt{n}} + \sum_{j\in\mathcal{S}^*\setminus\widehat{\mathcal{S}}}|\gamma_j|\frac{\|\widehat{\mathbf{v}}^{[j]}\|}{\sqrt{n}}}{\frac{1}{\sqrt{n}}\left\|\sum_{j\in\mathcal{S}^*}\gamma_j\widehat{\mathbf{v}}^{[j]}\right\|_2}$$

By Lemma 5, we have

$$\frac{\sum_{j\in\widehat{\mathcal{S}}}|\widehat{\gamma}_j - \gamma_j|\frac{\|\widehat{\mathbf{V}}^{[j]}\|}{\sqrt{n}} + \sum_{j\in\mathcal{S}^*\setminus\widehat{\mathcal{S}}}|\gamma_j|\frac{\|\widehat{\mathbf{V}}^{[j]}\|}{\sqrt{n}}}{\frac{1}{\sqrt{n}}\left\|\sum_{j\in\mathcal{S}^*}\gamma_j\widehat{\mathbf{v}}^{[j]}\right\|_2}$$

$$\leq \frac{\left(\sum_{j\in\widehat{\mathcal{S}}}|\widehat{\gamma}_j - \gamma_j| + \sum_{i\in\mathcal{S}^*\setminus\widehat{\mathcal{S}}}|\gamma_j|\right)\sqrt{\left(1+12\sqrt{\frac{\log p}{n}}\right)M_1}}{\sqrt{\frac{M_1\|\gamma\|_2^2(1-s_{z1}\lambda_n)^2}{1-12\sqrt{\frac{\log p}{n}}}}} \leq Cs_{z1}\sqrt{\frac{\log p}{n}} \leq \epsilon_0,$$

and hence the second convergence result of (161) follows from the following inequality.

## D.8. Proof of Lemma 8

Define $\widehat{\|\gamma\|}_2^2 = \sum_{j\in\widehat{\mathcal{V}}}\widehat{\gamma}_j^2$ and $\|\gamma_{\mathcal{V}^*}\|_2^2 = \sum_{j\in\mathcal{V}^*}\gamma_j^2$. The proof of this lemma is further based on the following results. Under the assumptions $(R1)-(R5)$ and (IN1)-(IN3). As $\sqrt{s_{z1}}s\log p/\sqrt{n}\to 0$, we have

$$\frac{\widehat{\|\gamma\|}_2^2}{\|\gamma_{\mathcal{V}^*}\|_2^2} \xrightarrow{p} 1 \quad \text{and} \quad \frac{\left\|\sum_{j\in\widehat{\mathcal{V}}}\widehat{\gamma}_j\widehat{\mathbf{v}}^{[j]}\right\|_2}{\left\|\sum_{j\in\mathcal{V}^*}\gamma_j\widehat{\mathbf{v}}^{[j]}\right\|_2} \xrightarrow{p} 1. \tag{165}$$

By (160) in the main paper, we establish that

$$\frac{\sqrt{\widehat{\boldsymbol{\Theta}}_{11} + \widehat{\beta}^2\widehat{\boldsymbol{\Theta}}_{22} - 2\widehat{\beta}\widehat{\boldsymbol{\Theta}}_{12}}}{\sqrt{\boldsymbol{\Theta}_{11} + (\beta)^2\boldsymbol{\Theta}_{22} - 2\beta\boldsymbol{\Theta}_{12}}} \xrightarrow{p} 1.$$

Combined with (165), we establish (84).

**Proof of** (165) For a given $0 < \epsilon_0 < 1$, we have

$$\mathbf{P}\left(\left|\frac{\|\gamma_{\mathcal{V}^*}\|_2^2}{\widehat{\|\gamma\|}_2^2} - 1\right| \geq \epsilon_0\right) \leq \mathbf{P}\left(\left|\frac{\widehat{\|\gamma\|}_2^2 - \|\gamma_{\mathcal{V}^*}\|_2^2}{\|\gamma_{\mathcal{V}^*}\|_2^2}\right| \geq \frac{\epsilon_0}{1-\epsilon_0}\right).$$

By Lemma 2, on the event $A \cap S \cap G \cap F$, we have $\widehat{\mathcal{V}} = \mathcal{V}^*$ and

$$\left|\frac{\widehat{\|\gamma\|}_2^2 - \|\gamma_{\mathcal{V}^*}\|_2^2}{\|\gamma_{\mathcal{V}^*}\|_2^2}\right| \leq C\frac{1}{\|\gamma_{\mathcal{V}^*}\|_2^2}\left(s_{z1}\frac{\log p_z}{n} + Cs_{z1}\left(s\frac{\log p}{n}\right)^2 + C\|\gamma_{\widehat{\mathcal{V}}}\|_2\sqrt{\frac{2|\widehat{\mathcal{V}}|\log p_z}{n}}\right). \tag{166}$$

Since $\|\gamma_{\mathcal{V}^*}\|_2^2 \gg (s\log p/\sqrt{n})^2$, we obtain that

$$\left|\frac{\widehat{\|\gamma\|}_2^2 - \|\gamma_{\mathcal{V}^*}\|_2^2}{\|\gamma_{\mathcal{V}^*}\|_2^2}\right| \leq \frac{\epsilon_0}{1-\epsilon_0} \text{ and } \mathbf{P}\left(\left|\frac{\|\gamma_{\mathcal{V}^*}\|_2^2}{\widehat{\|\gamma\|}_2^2} - 1\right| \geq \epsilon_0\right) \leq \mathbf{P}\left((A\cap S\cap G\cap F)^c\right).$$

Combined with Lemma 4 and (140), we establish the first converge result of (165).

On the event $S \cap G \cap A \cap F$, we have

$$\left| \frac{\left\| \sum_{j \in \widehat{\mathcal{V}}} \widehat{\gamma}_j \widehat{\mathbf{v}}^{[j]} \right\|_2}{\left\| \sum_{j \in \mathcal{V}^*} \boldsymbol{\gamma}_j \widehat{\mathbf{v}}^{[j]} \right\|_2} - 1 \right| \leq \frac{\sum_{j \in \mathcal{V}^*} |\widehat{\gamma}_j - \gamma_j| \frac{\|\widehat{\mathbf{v}}^{[j]}\|}{\sqrt{n}}}{\frac{1}{\sqrt{n}} \left\| \sum_{j \in \mathcal{V}^*} \boldsymbol{\gamma}_j \widehat{\mathbf{v}}^{[j]} \right\|_2}.$$

By Lemma 5, we have

$$\frac{\sum_{j \in \mathcal{V}^*} |\widehat{\gamma}_j - \gamma_j| \frac{\|\widehat{\mathbf{v}}^{[j]}\|}{\sqrt{n}}}{\frac{1}{\sqrt{n}} \left\| \sum_{j \in \mathcal{V}^*} \boldsymbol{\gamma}_j \widehat{\mathbf{v}}^{[j]} \right\|_2} \leq \frac{\left( \sum_{j \in \mathcal{V}^*} |\widehat{\gamma}_j - \gamma_j| \right) \sqrt{\left( 1 + 12\sqrt{\frac{\log p}{n}} \right) M_1}}{\sqrt{\frac{M_1 \|\boldsymbol{\gamma}\|_2^2 (1 - s_{z1}\lambda_n)^2}{1 - 12\sqrt{\frac{\log p}{n}}}}} \leq C s_{z1} \sqrt{\frac{\log p}{n}}.$$
$$(167)$$

Hence the second converge result of (165) follows from

$$\mathbf{P} \left( \left| \frac{\left\| \sum_{j \in \widehat{\mathcal{V}}} \widehat{\gamma}_j \widehat{\mathbf{v}}^{[j]} \right\|_2}{\left\| \sum_{j \in \mathcal{V}^*} \boldsymbol{\gamma}_j \widehat{\mathbf{v}}^{[j]} \right\|_2} - 1 \right| \geq \epsilon_0 \right) \leq \mathbf{P} \left( (S \cap G \cap A \cap F)^c \right).$$

### D.9. Proof of Lemmas 9, 10, 11 and 12

**Proof of Lemma 9** We only estalbish the first half of (90) and (91). The proof of the second half of (90) and (92) will be similar. The proof has been established in Sun and Zhang [2012] for fixed designs under certain assumptions for the design. In the following, we will check that the assumptions in Corollary 1 in Sun and Zhang [2012] are satisfied with high probability for the subgaussian random designs considered in this paper and then apply equation (23) in Sun and Zhang [2012]. By the definition of $\tau^*$ in Sun and Zhang [2012], we have $\tau^* \leq \tau$ where $\tau$ is defined in (30). Hence, on the event $S_1$, equation (23) in Sun and Zhang [2012] holds. By the relationship between $\ell_1$ cone invertibility factor and the restricted eigenvalue established in Lemma 13 of Cai and Guo [2016a], we obtain that on the event $S \cap G$,

$$\|\widetilde{\boldsymbol{\Gamma}} - \boldsymbol{\Gamma}\|_1 + \|\widetilde{\boldsymbol{\Psi}} - \boldsymbol{\Psi}\|_1 \leq C \frac{s\lambda_0 \sigma_1}{\kappa^2(\mathbf{H}, 4s, 1 + 2\epsilon_0)}. \tag{168}$$

Similar to the proof of Lemma 13 in Cai and Guo [2016a], we establish

$$\kappa^2 \left( \mathbf{H}, 4s, 1 + 2\epsilon_0 \right) \geq \frac{n}{\max \|\mathbf{W}_{.j}\|_2^2} \kappa^2 \left( \mathbf{W}, 4s, (1 + 2\epsilon_0) \left( \frac{\max \|\mathbf{W}_{.j}\|_2}{\min \|\mathbf{W}_{.j}\|_2} \right) \right). \tag{169}$$

Hence, on the event $G \cap S$, we establish the first half of (90). Since

$$\frac{1}{n} \|\mathbf{Z}(\widetilde{\boldsymbol{\Gamma}} - \boldsymbol{\Gamma}) + \mathbf{X}(\widetilde{\boldsymbol{\Psi}} - \boldsymbol{\Psi})\|_2^2 \leq \|\frac{1}{n} \mathbf{W}^{\mathsf{T}} \mathbf{W} \begin{pmatrix} \widetilde{\boldsymbol{\Gamma}} - \boldsymbol{\Gamma} \\ \widetilde{\boldsymbol{\Psi}} - \boldsymbol{\Psi} \end{pmatrix} \|_\infty \left( \|\widetilde{\boldsymbol{\Gamma}} - \boldsymbol{\Gamma}\|_1 + \|\widetilde{\boldsymbol{\Psi}} - \boldsymbol{\Psi}\|_1 \right),$$

we establish (91).

**Proof of Lemma 10** The decompositions (93) and (95) are established by the definitions of $\widehat{\boldsymbol{\Gamma}}$ and $\widehat{\boldsymbol{\gamma}}$. The error bound (97) follows from the following inequality

$$|\Delta_j^{\boldsymbol{\gamma}}| \leq \| \left( \frac{1}{n}(\widehat{\mathbf{U}}_{\cdot j})^{\intercal}\widehat{\boldsymbol{\Sigma}} - e_j^{\intercal} \right) \|_{\infty} \| \begin{pmatrix} \widetilde{\boldsymbol{\gamma}} - \boldsymbol{\gamma} \\ \widehat{\boldsymbol{\psi}} - \boldsymbol{\psi} \end{pmatrix} \|_1.$$

## Proof of Lemma 11

This lemma can be established by a similar argument with Lemma 5. On the event $A \cap S \cap G$, we have

$$\frac{1}{\sqrt{M_2}}\sqrt{1 + \left(\beta + \frac{\pi_1}{\gamma_1}\right)^2} \leq \sqrt{\boldsymbol{\Theta}_{11} + \left(\beta + \frac{\pi_1}{\gamma_1}\right)^2 \boldsymbol{\Theta}_{22} - 2\left(\beta + \frac{\pi_1}{\gamma_1}\right)\boldsymbol{\Theta}_{12}} \leq \sqrt{M_2}\sqrt{1 + \left(\beta + \frac{\pi_1}{\gamma_1}\right)^2},$$
$$(170)$$

and

$$\frac{\sqrt{M_1}\sqrt{1 + \left(\frac{\gamma_j}{\gamma_1}\right)^2}\,|1 - 2\lambda_n|}{\sqrt{1 - 12\sqrt{\frac{\log p}{n}}}} \leq \frac{\left\|\widehat{\mathbf{v}}^{[k]} - \frac{\gamma_j}{\gamma_1}\widehat{\mathbf{v}}^{[1]}\right\|_2}{\sqrt{n}} \leq \sqrt{M_1}\left(1 + \left|\frac{\gamma_j}{\gamma_1}\right|\right)\sqrt{1 + 12\sqrt{\frac{\log p}{n}}}.$$
$$(171)$$

Hence, (141) and (142) follow from the above inequalities (170) and (171).

## Proof of Lemma 12

(143) follows from the standard convergence analysis and (144) follows from high probability statement of Gaussian random variable. It remains to establish (145). We will analyze the expression (139) term by term. Note that on the event $A \cap S \cap G$,

$$\left|T^{\beta,1}\right| \leq \frac{\sqrt{\log p_z}}{\sqrt{n}\,|\gamma_1|}\|\widehat{\mathbf{v}}^{[1]}\|_2\sqrt{\boldsymbol{\Theta}_{11} + \left(\beta + \frac{\pi_1}{\gamma_1}\right)^2\boldsymbol{\Theta}_{22} - 2\left(\beta + \frac{\pi_1}{\gamma_1}\right)\boldsymbol{\Theta}_{12}}$$

$$\leq C\frac{1}{|\gamma_1|}\frac{\|\widehat{\mathbf{v}}^{[1]}\|_2}{\sqrt{n}}\sqrt{\log p_z}\left(1 + \left|\beta + \frac{\pi_1}{\gamma_1}\right|\right);$$

$$|R_1| \leq C\frac{1}{|\gamma_1|}\left(\left|\beta + \frac{\pi_1}{\gamma_1}\right| + 1\right)s\frac{\log p}{\sqrt{n}}; \qquad (172)$$

$$|R_2| \leq C\frac{1}{|\gamma_1|}\left(\sqrt{\frac{\log p}{n}} + s\frac{\log p}{n}\right)\left(\left|T^{\beta,1}\right| + |R_1|\right)$$

$$\leq C\frac{1}{|\gamma_1|}\sqrt{\frac{\log p}{n}}\frac{1}{|\gamma_1|}\frac{\|\widehat{\mathbf{v}}^{[1]}\|_2}{\sqrt{n}}\sqrt{\log p_z}\left(1 + \left|\beta + \frac{\pi_1}{\gamma_1}\right|\right).$$

Hence

$$\frac{1}{\sqrt{n}}\left|\gamma_j\Delta^{\beta,1}+\left(T^{\beta,1}+\Delta^{\beta,1}\right)(\widehat{\gamma}_j-\gamma_j)\right|\leq C\frac{|\gamma_j|}{\sqrt{n}}\left(|R_1|+|R_2|\right)+\sqrt{\frac{\log p_{\mathsf{z}}}{n}}\left|T^{\beta,1}\right|$$

$$\leq C\frac{|\gamma_j|}{|\gamma_1|}\left(\left|\beta+\frac{\pi_1}{\gamma_1}\right|+1\right)s\frac{\log p}{n}+C\frac{|\gamma_j|}{|\gamma_1|}\frac{\sqrt{\log p\log p_{\mathsf{z}}}}{n}\frac{1}{|\gamma_1|}\frac{\|\widehat{\mathbf{v}}^{[1]}\|_2}{\sqrt{n}}\left(1+\left|\beta+\frac{\pi_1}{\gamma_1}\right|\right)$$

$$+C\frac{1}{|\gamma_1|}\frac{\|\widehat{\mathbf{v}}^{[1]}\|_2}{\sqrt{n}}\frac{\sqrt{\log p\log p_{\mathsf{z}}}}{n}\left(1+\left|\beta+\frac{\pi_1}{\gamma_1}\right|\right)$$

$$\leq C\left(1+\left|\beta+\frac{\pi_1}{\gamma_1}\right|\right)\left(\frac{|\gamma_j|}{|\gamma_1|}s\frac{\log p}{n}+\left(1+\frac{|\gamma_j|}{|\gamma_1|}\right)\frac{1}{|\gamma_1|}\frac{\sqrt{\log p\log p_{\mathsf{z}}}}{n}\right).$$
(173)

Since

$$\left|\Delta_j^{\mathbf{\Gamma}}-\left(\beta+\frac{\pi_1}{\gamma_1}\right)\Delta_j^{\gamma}\right|\leq C\left(\left|\beta+\frac{\pi_1}{\gamma_1}\right|+1\right)s\frac{\log p}{n},$$

we have

$$\max_{j\in\mathcal{S}^*}\frac{1}{\sqrt{n}}|\Delta^{\pi_j}|\leq\left(1+\left|\beta+\frac{\pi_1}{\gamma_1}\right|\right)\left(1+\frac{|\gamma_j|}{|\gamma_1|}\right)\left(s\frac{\log p}{n}+\frac{1}{|\gamma_1|}\frac{\sqrt{\log p\log p_{\mathsf{z}}}}{n}\right)$$
(174)

By the assumption $\min_{j\in\mathcal{S}^*}|\gamma_j|\gg\sqrt{\log p/n}$ and (142), we establish (145).